\documentclass[11pt] {article} 
  \usepackage{amsmath}
    \usepackage{amssymb}
  \usepackage[dvipsnames]{xcolor}   
     \usepackage{pdfsync}

\newtheorem{example}{Example}[section]
\newtheorem{theorem}{Theorem}[section]
\newtheorem{lemma}{Lemma}[section]
\newtheorem{corollary}{Corollary}[section]

\newtheorem{remark}{Remark}[section]

\newcommand{\eqnsection}{
   \renewcommand{\theequation}{\thesection.\arabic{equation}}
   \makeatletter
   \csname @addtoreset\endcsname{equation}{section} 
   \makeatother}
\def \ov{\overline}

\def \be{\begin{equation}}
\def \ee{\end{equation}}
\def \bt{\begin{theorem}} 
\def \et{\end{theorem}}
\def \bl{\begin{lemma}} 
\def \el{\end{lemma}}
\def \bea{\begin{eqnarray}}
\def \eea{\end{eqnarray}}
\def \bas{\begin{eqnarray*}}
\def \eas{\end{eqnarray*}}

\def \al{\alpha}
\def \bb{\beta}
\def \ga{\gamma}

\def \de{\delta}
\def \De{\Delta}
\def \ep{\epsilon}
\def \vep{\varepsilon}
\def \la{\lambda}

\def \om{\omega}
\def \Om{\Omega}

\def \vf{\varphi}
\def \si{\sigma}

\def \th{\theta}

\def \ze{\zeta}

\def \ff{\infty}
\def \wh{\widehat}
\def \wt{\widetilde}

\def\stl{\stackrel{law}{=}}

\def \BB{{\cal B}}

\def \FF{{\cal F}}

\def \KK{{\cal K}}

\def \TT{{\cal T}}

\def \YY{{\cal Y}}

\def\b1{\mathbf 1}
\def \({\left(}
\def \){\right)}

\def \nn{\nonumber}
 
\def \Proof{\noindent{\bf Proof $\,$ }}

\def \bc{\begin{center} }
\def \ec{\end{center} }
\def \bs{\begin{slide} }
\def \es{\end{slide} }

\def\square{{\vcenter{\vbox{\hrule height.3pt
        \hbox{\vrule width.3pt height5pt \kern5pt
           \vrule width.3pt}
        \hrule height.3pt}}}}
\def\qed{{\hfill $\square$ \bigskip}}

\eqnsection

\begin{document}

\title{ Law of the iterated logarithm for $k/2$--permanental processes and the local times of   related Markov processes }
%related to the

 \author{  Michael B. Marcus\,\, \,\, Jay Rosen \thanks{Research of     Jay Rosen was partially supported by  grants from the Simons Foundation.   }}
\maketitle
 
 \begin{abstract}       Let $Y$ be a  symmetric   Borel right process  with locally compact state space $T\subseteq R^{1}$ and  potential    densities    $u(x,y)$ with respect to some $\sigma$-finite measure on $T$.
Let  $g$   and $f$ be   finite excessive functions    for  $ Y$. Set
 \be
  u_{g, f}(x,y)= u(x,y)+g(x)f(y),\qquad x,y\in  T. \nonumber
 \ee 

    In this paper we take $Y$ to be   a  symmetric    L\'evy process,   or a diffusion,   that is killed at the end of an independent exponential time or the first time it hits 0. 

 Under  general smoothness  conditions on   $g$, $f$,
  $u$ and points  $d\in T$, laws of the iterated logarithm are found for  $X_{k/2} =\{X_{k/2}(t), t\in T \}$, a $k/2-$permanental process with kernel   $ \{u_{g, f}(x,y),x,y\in T \}$, of the following form:
For   all integers $k\geq 1$,   
 \begin{equation}
 \limsup_{  x \to 0}\frac{|  X_{k/2}( d+x)-  X_{k/2}  (d)|}{ \(   2 \sigma^{2}\(x\)\log\log 1/x\)^{1/2}}=    \( 2 X  _{k/2} (d)\)^{1/2}, \qquad a.s. , \nonumber
 \ee 
 where,  
 \begin{equation} 
   \sigma^2(x)=u(d+x,d+x)+u(x,x)-2u(d+x,x).\nonumber
\end{equation}

Using these limit theorems and the Eisenbaum Kaspi Isomorphism Theorem, laws of the iterated logarithm are found for the local times  of  certain Markov processes with   potential    densities that have the form of  $ \{u_{g, f}(x,y),x,y\in T \}$  or are  slight modifications of it.  

\end{abstract}

\maketitle

 \section{Introduction}\label{sec-1}

Let $\{K(s,t),s,t\in T \}$ be a kernel with the property that for all ${\mathbf t_n} =(t_{1},\ldots,t_{n})$ in $T^n $, the matrix  $\mathbf{ \KK(      t_n)}=\{K(t_{i},t_{j}),i,j\in [1,n] \}$   determines an $n-$dimensional random variable random variable $ (X_{\al}(t_{1}),\ldots, X_{\al}( t_{n}))$ with  Laplace transform,   
\begin{equation}
  E\(e^{-\sum_{i=1}^{n}s_{i}  X_{\al }(t_i) }\) 
= \frac{1}{ |I+ \mathbf{ \KK(      t_n)}S(\mathbf{ s_n })|^{ \al}},   \label{int.1}
\end{equation}
  where $S( \mathbf{        s_n} )$ is a  diagonal matrix with positive entries $\mathbf{ s_n }=( s_{1},\ldots,s_{n} )$, and $\al>0$.     We refer to  $ \mathbf{ \KK(      t_n)}$ as the kernel of $(X_{\al }(t_1),\ldots,X_{\al }(t_n))$.     
It follows from the Kolmogorov Extension Theorem that $\{K(s,t),s,t\in T \}$ determines  a stochastic process which we denote by   $X_\al =\{X_{\al}(t),t\in T \}$ and refer to as an $\al-$permanental process.

There are well known examples of $\al-$permanental processes when $\al=k/2$.
Let $\eta=\{\eta(t);t\in  R^1\}$ be a mean zero   Gaussian process with covariance $U=\{U(s,t),s,t\in R^1\}.$   Let $\{\eta_{i};i=1,\ldots, k\}$ be independent copies of $\eta$. Define,
\begin{equation} \label{1.9mm}
 Y_{k/2}(t)=\sum_{i=1}^{k}\frac{\eta^2_{i}(t)}{2},\qquad t\in \TT.\end{equation}
The positive stochastic process  $\{Y_{k/2}(t),t\in R^1 \}$ is referred to as a chi--square process of order  $k$  with covariance $U$.
It is well known that for $\{t_1 ,\ldots,t_n\}\subset \TT$  
\begin{equation}
  E\(e^{ -\sum_{j=1}^{n}s_{j} Y_{k/2}(t_{j})   }\)
= \frac{1}{ |I+  \ov U(   \mathbf{    t_n})S(\mathbf{        s_n })|^{k/2}}, \label{int.1pp} 
\end{equation}
 where $   \ov U(\mathbf{t_n})$ is the $n\times n$ matrix, $\{
U(t_{j},t_{k}); 1\leq j,k\leq n\}$.  
 Therefore $Y_{k/2}=\{Y_{k/2}(t);t\in \TT\}$ is a  $k/2$--permanental process with  kernel $\ov U(n)$.  

 It is important to note that whereas (\ref{int.1pp}) holds for all integers $k$ whenever $U$ is the covariance of a Gaussian process, it may or may not hold when $k/2$
   is replaced by some $\al>0$ that is not a multiple of 1/2. (The question of when 
(\ref{int.1pp}) holds for all $\al>0$ was raised by Paul L\'evy when he asked when does a Gaussian process have infinitely divisible squares. This was answered by 
Griffiths and Bapat. An extensive treatment of these results is given in \cite[Chapter 13]{book}.) Nevertheless,
it follows from Eisenbaum and Kaspi, \cite[Theorem 3.1]{EK} that  if   $\wt U=\{\wt U(s,t),s,t\in T \}$ is 
a    potential density of a transient Markov process on   $T$ which is   finite for all $s,t\in T$,   then there exist     $\al-$permanental processes  $\{Y_{\al,t},t\in T \}$ with kernel      $\wt U$ for all $\al>0$.  

It is at this point that this really becomes interesting because in the Eisenbaum  Kaspi    theorem   the  kernel   $\wt U$ is not necessarily  symmetric,  or equivalent to a symmetric  matrix;  (see \cite{MRnonsym}).    When it is not,   the corresponding permanental processes are not   squares of a Gaussian process. They are a new kind of positive stochastic process defined by kernels that need not be symmetric. It is intriguing  to study  sample path properties of these processes. They can not be analyzed using standard Gaussian process techniques.  

\medskip In this paper we obtain LILs  for the local behavior of   a large class of $k/2-$permanental processes.
We can  get an idea of what they should look like by considering what   they are  for $\{\eta^2 (t)/2,t\in R^1 \}$,    where $\eta=\{\eta(t);t\in  R^1\}$ is a mean zero   Gaussian process   with covariance  $ \{u(x,y),x,y\in R^1 \}$. (As we point out above,  $\{\eta^2 (t)/2,t\in R^1 \}$,  is a $1/2-$permanental process  with kernel   $  u$.)  
 
 Let $\phi(t)$ be a positive function  on $[0,\ep]$ for some $\ep>0$.
Suppose that for $d\in R^1$,    where $u(d,d)=E(\eta^2(d))>0$,   
\begin{equation}
 \limsup_{  t \to 0}\frac{|  \eta(   d+t)-  \eta  (d)|}{ \( 2\phi(t)\log\log 1/t\)^{1/2}}=  1,\qquad a.s., \label{1.3aaq}
 \end{equation}
 or, equivalently,  
 \begin{equation}
 \lim_{\de\to 0}\sup_{  |t| \leq \de}\frac{|  \eta(   d+t)-  \eta  (d)|}{ \( 2\phi(t)\log\log 1/t\)^{1/2}}=  1,\qquad a.s.\label{1.3aaqe}
 \end{equation} 
  (These local moduli of continuity hold  for a very large class of Gaussian processes.     See e.g.,  \cite[Theorem 7.2.15]{book}.) Then since,
\begin{equation} 
 \frac{\eta ^{2}(d+t)}{2}-\frac{\eta ^{2}(d)}{2}=\frac{1}{2}(\eta( d+ t)-\eta( d))(\eta( d+t)+\eta( d)),
\end{equation}
we see that,   
\bea
 &&\label{1.3bb}\limsup_{  t \to 0}\frac{|  \eta^2(  d+t)/2-  \eta^2  (d)/2|}{ \( 2\phi(t)\log\log 1/t\)^{1/2}}\\ &&\quad= \limsup_{  t  \to  0}\frac{|  \eta ( d+ t) -  \eta   ( d) |}{  \( 2\phi(t)\log\log 1/t\)^{1/2}} |\eta (d)| \nn=  \sqrt2 \( \frac{  \eta^2 (d)  }{2}\)^{1/2}, \qquad a.s. 
 \eea
 In \cite{MRchi}  we generalize this to  Chi-square  processes,  $X_{k/2}(t)=\sum_{i=1}^{k} \eta^2_{i}(t)/2$,   $k\ge 1$, where  $\{\eta_{i},i=1,\ldots, k\}$ are independent copies of $\eta$. We find that  
for all   $k\ge 1$, 
\begin{equation}
  \limsup_{  t\to d}\frac{|  X_{k/2}(t)-  X_{k/2}  (d)|}{ ( 2 \phi(t) \log\log 1/t)^{1/2}}= \sqrt {2 }  X ^{1/2}_{k/2} (d) \qquad a.s.\label{1.34slr}
  \end{equation} 
   In this paper we extend this LIL to a wide class of $k/2$ permanental processes with kernels that are not symmetric and consequently are not Chi--square processes.

	\medskip	Let $T$  denote a locally compact space with a countable base and let $Y$ be a symmetric   transient Borel right process \index{Borel right process} with state space $ T$,  and   potential    densities  $u(x,y)$ with respect to some reference measure $m$.
   If there exist  functions    $u^{\bb}= \{u^{\bb}(x,y),x,y\in T\}$ for which,  
   \begin{equation}
 E^{ x} \(   \int_{0}^{ \ff}e^{- \bb t}h\(     Y_{t}\)\,dt\)=\int_T   u^{\bb}(x,y) h(y) \,dm(y),\label{potdef7}
   \end{equation}
 for all $x\in T$ and all  positive measurable functions $h$  on $T$, we refer to $u^\bb$
 as the   $\bb-$potential densities of $Y$  with respect to  the reference measure $m$.  (In this case we say that $Y$ is strongly symmetric, \cite[p. 75]{book}.)   When    (\ref{potdef7}) holds for $\bb=0$,    we say that $0-$potential densities exist     and refer to them  simply as the potential  densities of $Y$.    In addition  we use the convention that   the domain of  $h$ is extended to include $\De$,  the cemetery state for $Y$, and   $h(\De)=0$. (This is necessary for (\ref{potdef7}) to make sense.)
 
% When $\bb>0$ the $\bb-$potential densities of $Y$ are the   0--potential  densities of a  Borel right process that is obtained by killing   $Y$ at an independent exponential time with mean $1/\bb$. 
 
\medskip Let  $\{f(x)$, $x\in T\}$ be an excessive function\index{excessive function} for   $ Y$.  (A function $f$ is   excessive for $Y$ if $  E^{x}\(f( Y_{t})\)\uparrow  f(x)$ as $t\to 0$ for all  $x\in T$.)     
  It is easy to check that for any positive measurable function $h$,  
  \begin{equation}
   f(x)=\int_T   u(x,y) h(y) \,dm(y)=E^{ x} \(   \int_{0}^{ \ff}h\(     Y_{t}\)\,dt\),\label{potdef}
   \end{equation}
  is excessive for $Y$. 
Excessive functions  that can be expressed in this way are called  potential functions for $Y$.

  Let   $\{f(x)$, $x\in T \}$ and $\{g(x)$, $x\in T \}$ be   excessive functions\index{excessive function} for   $ Y$. Set  
\be
   u _{g, f}(x,y)= u(x,y)+g(x)f(y),\qquad x,y\in  T.\label{1.3qssj}
  \ee 
	 The fact that $u_{g, f}:=\{u _{g, f}(x,y),x,y\in  T \}$ is the kernel of $\al-$permanental processes    is given by the next theorem.  When $g\equiv 1$ this result is     in  \cite[Theorem 1.11]{MRejp}.  

    \bt\label{theo-borelN} 
    
Let ${ Y}\!=\!
(\Om,  \FF_{t}, Y_t,\th_{t},P^x
)$ be a  symmetric   transient Borel right process \index{Borel right process} with state space $ T$, and  continuous  potential    densities  $u(x,y)$ with respect to some   locally bounded reference measure $m$. Then for any   finite excessive functions   $g$ and $f$ for  $   Y$ and $\al>0$, $u_{g, f}$  
is the kernel\index{kernel} of an $\al$--permanental process. 
\et
In \cite[Theorem 1.11]{MRejp} we require  that  $u(x,y)>0$ for all $x,y\in T$. This is  used, along with \cite[(2.5)]{FR},  to show that the reference measure $m$   is locally bounded. In the statement of Theorem \ref{theo-borelN} we simply include the statement that the reference measure $m$   is locally bounded  in the hypothesis.   
 
\medskip
When $ Y $  in Theorem \ref{theo-borelN} has    continuous  potential densities $u=\{u(x,y),\newline x,y\in T\}$ with respect to some reference measure on $T$, 
   it follows from  \cite[Lemma 3.3.3]{book} that  $\{u(x,y);x,y\in  T\}$    is positive definite and therefore is  the covariance  of a Gaussian process. We refer to Gaussian processes with covariances that are the potential densities of transient Markov processes as associated Gaussian processes. By   \cite[Lemma 3.4.3]{book} these covariances have the property that,   
    \begin{equation} \label{smp.1}
  u(x,y)\le u(x,x)\wedge u(y,y), \hspace{.2 in}\forall x,y\in T.
\end{equation}
   (This is a significant condition since, in general,   covariances of Gaussian processes only satisfy $u(x,y)\le (u(x,x)  u(y,y))^{1/2}$.)
    It follows from \cite[Lemma 3.4.6]{book}  that (\ref{smp.1}) implies that,
  \begin{equation}
 u(x,x)>0, \hspace{.2 in}\forall x\in T.\label{poscond.8}
  \end{equation}

    \medskip  
 To obtain LILs for permanental processes with kernels of the form of (\ref{1.3qssj}) we must identify the symmetric potential densities $u(x,y)$. In   Theorems \ref{theo-1.2zz}  we take them to be the potential densities of   L\'evy processes killed at    a time that is an independent exponential random variable.  

Let $Z=\{Z_{t};t\in R^1 \}$ be a real valued symmetric 
L\'evy process with characteristic exponent $\psi$, i.e., 
\begin{equation} \label{1.27nn}
  E\(e^{i\la Z_{t}}\)=e^{-t\psi(\la)},\qquad \forall\,t\ge 0,
\end{equation}
where,
\begin{equation} \label{1.20nn}
  \frac{1}{ 1+\psi(\la)}\in L^{1}(R^1).
\end{equation}  For $\bb>0$, the $\bb-$potential density of $Z$ with respect to Lebesgue measure  is,  
\begin{equation} \label{1.21nn}
u^{\bb}:=u^{\bb}(x,y)=\frac{1}{\pi}\int_{0}^{\ff}\frac{ \cos\la (x-y)}{\bb+ \psi(\la)}\,d\la,\qquad x,y\in R^{1};  
\end{equation}
see e.g.,   \cite[(4.84)]{book}.   Note that the potential density $u^{\bb}$ is the  $0-$potential density of the process   $Y=\{Y_t,t\in R^1 \}$ that is obtained by killing  $Z$ at the end of an independent exponential time with mean $1/\bb$. These are the processes   that we consider in   Theorem  \ref{theo-1.2zz}.  

\medskip  Note that the function $\{u^{\bb}(x,y),x,y\in R^1 \}$
is   translation invariant and consequently it is  the covariance of a stationary Gaussian process.   The proof of Theorem \ref{theo-1.2zz} uses    sample path properties of these Gaussian processes. 
When we consider translation invariant potential densities like $u^{\bb}(x,y)$ we alternately write them as,
\begin{equation} \label{ti}
  u^\bb(x,y)\qquad\text{and}\qquad  u^\bb(x-y).
\end{equation}
 
 Throughout this paper we  use the notation,
  \begin{equation} \label{deld}
  \De_d(\de)=[d-\de,d+\de] . 
\end{equation} \vspace{-.2in}
 
 \begin{theorem}\label{theo-1.2zz}  {\rm I}.
  Let $Y$ be a   symmetric  L\'evy process in $R^1$   killed at an exponentially distributed time with mean $1/\bb$, with    continuous       potential densities $u^\bb=\{u^{\bb}(x-y) ,x,y\in R^1\}$    with respect to Lebesgue measure.    
     For any integer $k\ge 1$, let $  X_{k/2}=\{  X_{k/2}(t), t\in  R^{ 1}\}$ be a    $k/2-$permanental process with kernel 
 \be
    u^\bb _{g, f}(x,y)=  u^\bb (x-y)+g(x)f(y),\qquad x,y\in R^{ 1}.\label{1.wwq}
 \ee 
  where    $g,f$ are    excessive functions for   $Y$.

  Set 
\begin{equation} \label{sigti.1}
 (\si^\bb)^{2}\(x\)=2\(u^\bb(0)-u^\bb(x)\) . 
 \end{equation} 
 Assume that,  
 \begin{itemize} 
  \item[(i)]$(\si^\bb)^{2}\(x\)$ is  regularly varying  at $0$ with index $0<r\leq1$;
   \item[(ii)] 
   for some   $p>1$ and $\de>0$,
\begin{equation} \label{1.38nna}
   (\si^\bb)^{2}\(x\) \ge x\( \log 1/x\)^{ p},\qquad x\in [0,\de].
\end{equation}
\end{itemize}
Then 
 for any   $d\in R^1$  and  $g, f\in C^{1}(\De_d(\de))$    for some $\de>0$,    \begin{equation}
  \limsup_{  x \to 0}\frac{|  X_{k/2}( d+x)-  X_{k/2}  (d)|}{ \(2   (\si^\bb)^{2}\(x\)\log\log 1/|x|\)^{1/2}}=    \( 2X  _{k/2} (d)\)^{1/2}, \qquad a.s.\label{1.34szqa}
  \end{equation}

  {\rm II.} This LIL   continues to hold   without requiring condition $(ii)$   when $g, f\in C^{2}(\De_d(\de))$,   for some $\de>0$, and   $g'(d)=f'(d)=0$.   \end{theorem}

  The random variable $X_{k/2}(d)$ is a  chi--square random variable of order $k$  which is the sum of $k$ independent squares of a normal random variable with mean zero and variance $u^\bb(0) + g(d)f (d)$.

   \medskip   We now explain the significance of   Theorem \ref{theo-1.2zz}, II.  To begin we note that $\psi(\la)$, the characteristic exponent of a L\'evy process can be written as,
   \begin{equation} \label{1.32nn}
  \psi(\la)=C\la^2+\psi_1(\la),\qquad\text{where $\psi_1(\la)=o(\la^2)$ as $\la\to\ff$}.
\end{equation}
When $C>0$ the L\'evy process has a Gaussian component and,  
\begin{equation} \label{1.29nn}
   (\si ^{  \bb  })^2(x)\sim \frac{|x|}{C},\qquad \text{as $x\to 0$.}
\end{equation} 
This follows from (\ref{1.26mmw2})  below because  $C_{2}=1$;  (see \cite[3.741.3]{GR}).    
 Consequently,  the condition in (\ref{1.38nna})   implies that  the L\'evy processes $Y$ do  not have a Gaussian component. 
 
 We also want to consider L\'evy processes with a Gaussian component.  
 Taking the Fourier transform we see that,  
\begin{equation} \label{1.21nnqq}
{e^{-a |y-x|}\over a} =\frac{1}{ \pi}\int_{-\ff}^{\ff}\frac{ \cos\la (x-y)}{a^2+ \la^2}\,d\la.  
\end{equation}
 Therefore, for $\bb>0$, 
  \begin{equation}      u^{\bb}_C(x-y) :=\frac{1}{ \pi}\int_{0}^{\ff}\frac{ \cos\la (x-y)}{\bb+ C\la^2}\,d\la=\frac{1}{2 }{e^{- ({ \bb/C})^{1/2}\,| x-y|}\over  ({ \bb C})^{1/2}},\qquad x,y\in  R^{ 1}.  \label{exe.13ja}
\end{equation} 
When $C=1/2$ this is the potential   density of     Brownian motion  killed after an independent exponential time   with mean $1/\bb$.  
  In these cases the   potential densities $\{u^{\bb}(x,y);x,y\in R^{1} \}$ are  also the covariances of  Ornstein--Uhlenbeck processes on the positive half line,  which are stationary Gaussian process on $R^1$. 
  
  Note that for all $\bb>0$,
  \begin{equation} \label{1.27mm}
   (\si_C^{\bb})^2\(x\)=2(u_C^{\bb}(x,x)-u_C^{\bb}(x,0))\sim \frac{|x|}{C},\qquad \text{as $x\to 0$},
\end{equation}
  does not depend on $\bb$.
   Theorem \ref{theo-1.2zz}, II  enables us to consider 
L\'evy processes $Y$ with a Gaussian component.  
   When $ (\si^\bb)^{2}\(x\)$ is  regularly varying  at $0$ with index $0<r<1$, (\ref{1.38nna}) is satisfied and we don't need the restrictive condition that $f'(d)=g'(d)=0$.

  \medskip 

 For all $\bb\ge 0$ set,  
\begin{equation}
 (\si^{\bb})^2\(x\)=\frac{ 2}{\pi}\int_{0}^{\ff}\frac{1- \cos(\la x)}{\bb+ \psi(\la)}\,d\la.\label{lv.21}
\end{equation}
(This is (\ref{sigti.1}) when $\bb>0$.)
   In Section \ref{sec-ex}  we show that    there exist characteristic exponents  $\psi$ that are asymptotic to every regularly varying   function at infinity with index $1<r<  2$, and to a large class of regularly varying   functions at infinity with index 2 subject to the condition that   $\psi(\la)=o(\la^{2}) $.
 When $\psi$ is asymptotic to a regularly varying   function at infinity with index  $1<r\le  2$,  
      \begin{equation} \label{1.26mmw2}
  (\si ^{  \bb  })^2(x)\sim C_r \frac{1}{|x|\psi(1/x)},\qquad  \text{as $x\to 0$}, \, \forall\,  \bb\ge 0 ,
\end{equation}
 where,   
    \begin{equation} 
  \label{1.26mmpa}
  C_r=\frac{4}{\pi}\int_{0}^{\ff}  \frac{\sin^2 s/2}{s^r}\,ds. 
\end{equation}
Consequently, the denominator in (\ref{1.34szqa}) does not depend on $\bb$, although it does depend on $r$.

 \medskip  We now consider  a class of symmetric potential  densities that are not translation invariant.    As in 
 Theorem \ref{theo-1.2zz} let   $Y $ be a   symmetric   L\'evy process in   $R^1$ that is killed at the end of an exponential time with mean $1/\bb$, with  continuous translation invariant      potential densities $u^\bb=\{u^\bb(x,y),x,y\in R^1\}$ with respect to Lebesgue measure.   
  Let $Y' $ be the symmetric  Borel right  process with state space   $T=R^1-\{0\}$ obtained by killing  $Y$ the first time it hits $0$. The process   $Y'$ has  continuous  potential densities with respect to Lebesgue measure restricted to $T$, which we denote by $v^{\bb}=\{v^{\bb}(x,y),x,y\in T\}$,  where   
  \bea\label{1.35nn}
 v^{\bb} (x,y)&=&u^\bb(x,y)-\frac{u^\bb(x,0)u^\bb(0,y)}{u^\bb(0,0)} \label{80.4a}\\
 &=&u^\bb(x-y)-\frac{u^\bb(x)u^\bb(y)}{u^\bb(0)};\nn
 \eea
see \cite[(4.165)]{book}.

   \begin{theorem}\label{theo-1.2zzk}   {\rm I.} Let  $Y' $  be a  symmetric  Borel right  process with state space   $T=R^1-\{0\}$ as described above with continuous  potential densities   $v^\bb=\{v^\bb(x,y),x,y\in T\}$ with respect to Lebesgue measure restricted to $T$.
   Let $(\si^{\bb})^2(x)$ be as in (\ref{sigti.1}) 
 and assume that;
    \begin{itemize}
  \item[(i)]
  $  (\si^{\bb})^2(x)$ is regularly varying at $0$ with index $0<r\le 1$;
  \item[(ii)]  for some   $p>1$,   and $\de>0,$
\begin{equation} \label{1.38nnb}
  (\si^{\bb})^2(x) \ge x\( \log 1/x\)^{ p},\qquad x\in [0,\de];
\end{equation}   \item[(iii)]   $ (\si^{\bb})^2(x)\in C^{1}\(T\).$ 
 
  \end{itemize}

For integers $k\ge 1$, let $  X_{k/2}=\{  X_{k/2}(t), t\in  T\}$ be a    $k/2$-permanental process with kernel  
  \be
     v^{\bb} _{g, f}(x,y)= v^{\bb} (x,y)+ g (x)f(y),\qquad x,y\in T,\label{80.5}
  \ee  
  where $g$ and $f $ are excessive functions for  $Y'$.    Then 
  for any  $d\neq 0$,  when  $g, f\in C^{1}(\De_d(\de))$  for some $\de>0$,  
  \begin{equation}
   \limsup_{  x \to 0}\frac{|  X_{k/2}( d+x)-  X_{k/2}  (d)|}{\( 2(\si^{\bb})^2(x) \log\log 1/|x|\)^{1/2}}=    \( 2X  _{k/2}  (d)\)^{1/2}, \qquad a.s.\label{80.6}
   \end{equation}

  {\rm II.}  This LIL continues to hold without requiring condition $(ii)$  when $g, f\in C^{2}(\De_d(\de))$,  for some $\de>0$,   $g'(d)=f'(d)=0$ and $\psi(\la)$ is a stable mixture, (defined immediately below).
   \end{theorem}

  \medskip In Section \ref{sec-ex} we show that condition $(iii)$ of Theorem \ref{theo-1.2zzk}, that $ (\si^{\bb})^2(x)\in C^{1}\(T\)$, is satisfied when   $\psi(\la)$ is the   characteristic exponent of L\'evy processes  that are stable mixtures. That is when,   
\begin{equation} \label{4.13q}
 \psi\(\la\)=\int_{1}^{2}  \la ^{s}\,d\mu(s),
\end{equation}
where $\mu$ is a finite positive measure on $(1,2]$. See Remark \ref{rem-sm}.   

 \begin{remark}\label{rem-BM} {\rm   There is a further restriction that has to be added to Theorem \ref{theo-1.2zzk}, II. When $Y$ is exponentially killed Brownian motion starting at some point $x>0$, then because its  paths are continuous, the process $Y' $  obtained by killing  $Y$ the first time it hits $0$ has  state space   $(0, \ff)$ not $R^1-\{0 \}$. Similarly, if it starts at    $x<0$}.
\end{remark}

   We now consider  (\ref{1.21nn}) with $\bb=0$,   in which case, 
\begin{equation} \label{1.21nn0}
 u^{0}(0,0)=\frac{1}{\pi}\int_{0}^{\ff}\frac{ 1}{  \psi(\la)}\,d\la .  
\end{equation} 
Note that when
\begin{equation} 
  \frac{\psi(\la)}{\la}=O(1),\qquad \text{as $\la\to 0$},
\end{equation} 
 $u^{0}(0,0)=\ff.$   Nevertheless, by modifying  $u^{0}(x,y)$   we can obtain the potential densities of    Markov processes. 
  
  Let  
 \begin{equation} \label{lv.22}
\phi(x)=\frac{1}{\pi}\int_{0}^{\ff}\frac{1-\cos\la x}{ \psi(\la)}\,d\la ,\qquad x \in R^{1}.  
\end{equation}
(The integral near $\la=0$ is finite and continuous in $x$ because $\psi$ is the characteristic exponent of  a L\'evy process.  See   \cite[(4.77)]{book}. It is finite  and continuous in $x$ on the rest of its domain by (\ref{1.20nn}).)  
Set,  
\begin{equation} \label{u0}
u^{(0)}=u^{(0)} (x,y)=\phi(x)+\phi(y)-\phi(x-y) ,\qquad x,y\in R^1-\{0\}.  
\end{equation}
It follows from \cite[Theorem 4.2.4]{book}  that $u^{(0)}$  is the potential density with respect to Lebesgue measure of the L\'evy process with characteristic exponent $\psi(\la)$, say $\ov Z$, killed the first time it hits 0. Denote this process  by   $\ov Y $.    
For example if $\ov Z$ is a symmetric stable process $u^{(0)}$ is the covariance of  fractional Brownian motion.

 Since $\{u^{(0)}(x,y)$, $x,y\in R^1-\{0\} \}$ is the potential density  of a  symmetric  Borel right  process it is positive definite and therefore it is the covariance of a Gaussian process. Furthermore, since  $\{u^{(0)}(x,y)$, $x,y\in R^1 \}$ is continuous and  $u^{(0)}(0,0)=0$, it is
also the covariance of a Gaussian process, say  $\eta=\{\eta(x),x\in R^{1} \}$, with $\eta(0)=0$. We define,
\bea \label{sig0}
  (\si^{(0)})^2(x,y)&=&E\(\eta(x)-\eta(y)\)^2 \\
  &=&\nn u^{(0)} (x,x)+u^{(0)} (y,y) -2u^{(0)} (x,y)=2\phi(x-y).
\eea
Note that $\eta$ has stationary increments;   
see Remark \ref{rem-3.2}. Therefore we write,
  \label{ref-3.2a}
\begin{equation} 
(\si^{(0)})^2(x,y)=(\si^{(0)})^2(x-y).
\end{equation}
  We have,
  \begin{equation} \label{1.44mm}
   (\si^{(0)})^2(x)=2\phi(x)= \frac{2}{\pi}\int_{0}^{\ff}\frac{1-\cos\la x}{ \psi(\la)}\,d\la =(\si^{0})^2(x),
\end{equation}
where $(\si^{0})^2(x)$ is given in (\ref{lv.21}).

  \begin{theorem}\label{theo-1.2aa}  {\rm I.} Let   $\ov Y$ be a  symmetric  Borel right  process with state space   $T=R^1-\{0\}$ and  continuous   potential densities  $u^{(0)} $,   (defined in  (\ref{u0})).
     Let $(\si^{0})^2(x)$ be as in (\ref{1.44mm}) 
 and assume that;
   \begin{itemize}
  \item[(i)]
  $(\si^{0})^2(x)$ is regularly varying at $0$ with index $0<r\le 1$; 
  \item[(ii)]  for some   $p>1$,   and $\de>0,$
\begin{equation} \label{1.38nnqq}
(\si^{0})^2(x) \ge x\( \log 1/x\)^{ p},\qquad x\in [0,\de];
\end{equation} 
 \item[(iii)]     $(\si^{0})^2(x)\in C^{1}(T)$.  
  \end{itemize}

For integers $k\ge 1$, let $  X_{k/2}=\{  X_{k/2}(t), t\in  R^{ 1}-\{0\}\}$ be a    $k/2$-permanental process with kernel 
  \be
     u _{g, f}(x,y)= u^{(0)} (x,y)+g(x)f(y),\qquad x,y\in R^{ 1}-\{0\},\label{1.wwqa}
  \ee  
 where   $g, f $ are excessive functions for   $\ov Y$.   
  Then,    for any $d\neq 0$,  when   $g, f\in C^{1}(\De_d(\de))$  for some $\de>0$,
  \begin{equation}
   \limsup_{  x \to 0}\frac{|  X_{k/2}( d+x)-  X_{k/2}  (d)|}{ \(2(\si^{0})^2 (x) \log\log 1/|x|\)^{1/2}}=     \(2 X  _{k/2} (d)\)^{1/2}, \qquad a.s.\label{1.34szw}
   \end{equation}

 {\rm II.}  
  This LIL continues to hold without requiring condition $(ii)$   when    $g, f\in C^{2}(\De_d(\de))$, for some $\de>0$,          $g'(d)=f'(d)=0$   and $ (\si^{(0)})^2(x) $ is increasing in some neighborhood of   $d$. 
    \end{theorem}
   It follows from \cite[Lemma 9.6.3]{book}  that  $ (\si^{(0)})^2(x) $ is increasing when $\psi(\la)$ is a stable mixture.

\medskip  We now consider permanental processes for which the symmetric part of the kernel is the potential density of a diffusion. 
 Let  $\mathcal{Z}$ be a  diffusion in $R^1$ that is regular and without traps  and  is symmetric with respect to some $\si$--finite measure $m$.  
(A diffusion is   regular and without traps when
$P^{x}\( T_{y}<\ff\)>0,   \forall x,y\in  R^1$.) We   assume that the second order generator for $\mathcal{Z}$ has continuous coefficients.  

Analogous to the study of the L\'evy process $Z$ 
 on page \pageref{1.27nn} we consider three transient symmetric diffusions determined by $\cal Z$:
\begin{itemize} 
\item[\bf a.]     $Z$ is   $\cal Z$  killed  at the end of an independent exponential time with mean $1/\bb$. The potential of $Z$ is
 \begin{equation} \label{diff.1}
\wt u^\bb (x,y)= \left\{
 \begin{array} {cc}
 p(x)q(y),& \quad x\leq y  
 \\
 q(x)p(y),& \quad y\leq x   
\end{array}  \right. ,
\end{equation}
where $p,q\in C^{2}(R^{1})$  are positive and  strictly convex   and $p$ is  increasing and  
$q$ is decreasing.   (For the dependence of $p$ and $q$ on $\bb$ see  \cite[(4.14)]{book}.)
\item[\bf b.]      $Z'$ is the  process with state space  $T=(0,\ff) $ that is obtained by starting   $  Z$ in $T$ and then  killing it  the first time it hits 0. The potential of $Z'$ is   
\begin{equation}
\wt v^\bb  (x,y)=\wt u^\bb (x,y)-\frac{\wt u^\bb (x,0)\wt u^\bb (0,y)}{\wt u^\bb(0,0)},\qquad  x,y>0, \label{diff.22}  
\end{equation}

so that
 \begin{equation} \label{diff.1v}
\wt v^\bb (x,y)= \left\{
 \begin{array} {cc}
 p(x)q(y)-\displaystyle\frac{p(0)}{q(0)}q(x)q(y),& \quad 0<x\leq y  
 \\\\
 q(x)p(y)-\displaystyle\frac{p(0)}{q(0)}q(x)q(y),& \quad 0<y\leq x.   
\end{array}  \right. 
\end{equation}
 
\item[\bf c.]       $\ov Z$ is the process with state space  $T=(0,\ff)$ that is obtained by starting  $\cal Z$ in $T$ and then  killing it  the first time it hits 0.  
$\ov Z$ has potential densities
  \begin{equation} \label{diff.5}
 \wt u_{T_{0}}(x,y)=s(x)\wedge s(y),\qquad  x,y>0,
\end{equation}
where   $s\in C^{2}(0,\ff)$  is   strictly positive and  strictly increasing.

 \end{itemize}

   \begin{theorem}\label{theo-1.5} {\bf a}. 
       For any integer $k\ge 1$, let $  X_{k/2}=\{  X_{k/2}(t), t\in  R^{ 1}\}$ be a    $k/2-$permanental process with kernel 
 \be
    \wt u^\bb _{g, f}(x,y)=  \wt u^\bb (x,y)+g(x)f(y),\qquad x,y\in R^{ 1}.\label{diffa.1x}
 \ee  
where   $f,g$ are excessive functions for   $Z$.
 Then for any   $d\in R^1$ for which  $g, f\in C^{2}(\De_{d}(\de))$ for some $\de>0$, and $g'(d)=f'(d)=0$,  
 \begin{equation}
  \limsup_{  x \to 0}\frac{|  X_{k/2}( d+x)-  X_{k/2}  (d)|}{ \(2 \tau(d)|x|\log\log 1/|x|\)^{1/2}}=    \( 2X  _{k/2} (d)\)^{1/2}, \qquad a.s.\label{1.34szq}
  \end{equation}
  where,
  \begin{equation} \label{8.51q}
 \tau(d)=\wt u^\bb(d,d)\frac{d}{dx}\(\log \frac{ p(x)}{q(x)}\)\bigg|_{x=d} \,.
\end{equation} 
{\bf b}.  
       For any integer $k\ge 1$, let $  X_{k/2}=\{  X_{k/2}(t), t\in (0,\ff)\}$ be a    $k/2-$permanental process with kernel 
 \be
    \wt v^\bb _{g, f}(x,y)=  \wt v^\bb (x,y)+g(x)f(y),\qquad x,y>0.\label{xkk}
 \ee  
where   $f,g$ are excessive functions for   $Z'$.
 Then for any   $d>0$ for which  $g, f\in C^{2}(\De_{d}(\de))$ for some $\de>0$,   (\ref{1.34szq})
   continues to hold for $x\to 0$.  
 
 \medskip\noindent {\bf c}. 
       For any integer $k\ge 1$, let $  X_{k/2}=\{  X_{k/2}(t), t\in (0,\ff)\}$ be a    $k/2-$ permanental process with kernel 
 \be
    \wt u_{T_0, g, f}(x,y)=  \wt u_{T_0} (x,y)+g(x)f(y),\qquad x,y>0.\label{diffa.1}
 \ee  
     where $ f,g $  are excessive functions for    $\ov Z$.  Then for any   $d>0$ for which  $g, f\in C^{2}(\De_{d}(\de))$ for some $\de>0$, and $g'(d)=f'(d)=0$,  
 \begin{equation}
  \limsup_{  x \to 0}\frac{|  X_{k/2}( d+x)-  X_{k/2}  (d)|}{ \(2 s'(d)|x|\log\log 1/|x|\)^{1/2}}=    \( 2X  _{k/2} (d)\)^{1/2}, \qquad a.s.\label{diffa.3q}
  \end{equation}
  
     \end{theorem}
     
      \begin{remark} {\rm \label{rem-tau} To reconcile the terms in the denominators of (\ref{1.34szqa}) and (\ref{1.34szq}) we show in Lemma \ref{lem-4.9} that
  $(\wt \si^\bb )  ( x-y )\sim
 \tau(x)|x-y |$  as $|x-y|\to 0$.
}\end{remark}

     \begin{remark} {\rm Killing a diffusion on $R^{1}$ the first time it hits $0$ gives rise to three  different processes: One on $(0,\ff)$, a second one on $(-\ff, 0)$, and a third process on $R^{1}-\{0\}$. In the third case, if the process  starts in $(0,\ff)$ it remains in $(0,\ff)$ throughout its lifetime, and similarly if it in starts in $(-\ff,0)$.  In this paper we consider the process on $(0,\ff)$.
}\end{remark}

      \begin{remark} {\rm 
  In all the LILs in Theorems \ref{theo-1.2zz}--\ref{theo-1.5} we take $\limsup_{x\to 0}$. The statements also hold for $\limsup_{x\uparrow 0} $ and   $\limsup_{x\downarrow  0}$. This observation is particularly relevant to the LIL in  Theorem \ref{theo-1.5}, \textbf{c.}, since we show in Section \ref{sec-9.9}  that   when $g, f\in C^{2}(0,\ff)$ are excessive functions for   $\ov Z$ and $g'(d)=f'(d)=0$ for some $d>0$, then $g(x)$ and $f(x) $ are constant for 
 $x\geq d$. Consequently $ \wt u_{T_0, g, f}(x,y)$ in (\ref{diffa.1}) is symmetric for  $x, y\geq d$.  Therefore we want to consider the analogue of  (\ref{diffa.3q}) with $\limsup_{x\to 0}$ replaced by  $\limsup_{x\uparrow 0} $.  
}\end{remark}

   We give two specific examples of LILs for permanental processes. In the first we consider the situation in Theorem \ref {theo-1.2aa} when,   
 \begin{equation} 
  \label{1.26mmpaz}
\phi(x)=\frac{4|x|^{\rho}}{\pi}\int_{0}^{\ff}  \frac{\sin^2 s/2}{s^{\rho+1}}\,ds:= |x|^{\rho}  C_{\rho+1},\qquad \rho\le 1. 
\end{equation}
Here $  \phi(x)+\phi(y)-\phi(x-y)$ is the potential of a symmetric stable process with index $\rho+1$, when $\rho<1$,  or a constant times Brownian motion  when $\rho=1$, killed the first time it hits 0.  Label this process $\wt Y$. The function $\phi(x)$ is given by (\ref{lv.22}) with $\psi(\la)=|\la|^{\rho+1}$. 
 
 For integers $k\ge 1$, let $  Z_{k/2}=\{  Z_{k/2}(t), t\in  R^{ 1}-\{0\}\}$ be a    $k/2$-permanental process with kernel 
    \be
   C_{\rho+1}\( |x|^{\rho}  +|y|^{\rho}    -|x-y|^\rho\) +g(x)f(y),\qquad x,y\in R^{ 1}-\{0\},\label{1.wwqb}
  \ee  
 where $g, f$ are excessive functions for $\wt Y$, with the property that   $g, f\in C^{1}$ on $[d,d+\de]$, $0  \notin[d,d+\de] $.   
 Then when $\rho<1$,  
  \begin{equation}
   \limsup_{  x \to 0}\frac{|  Z_{k/2}( d+x)-  Z_{k/2}  (d)|}{ \( 4C_{\rho+1}\,|x|^\rho \log\log 1/|x|\)^{1/2}}=    \(2 Z  _{k/2} (d)\)^{1/2}, \qquad a.s.\label{1.34szwq}
   \end{equation}
  If in addition    $g, f\in C^{2}(T)$ and    $g'(d)=f'(d)=0$  this LIL continues to hold   when $\rho=1$.  (In which case $|x|   +|y|     -|x-y| =|x|\wedge |y|$ is the covariance of    Brownian motion.)

   \medskip   Here $Z_{k/2}(d)$ is a  chi--square random variable of order $k$, as described in the paragraph following Theorem \ref{theo-1.2zz}, in which the  mean zero normal random  variable has  variance   $ 2C_{\rho+1}|d|^{\rho}+ g(d)f (d)$.

  \medskip  We next consider the situation in  Theorem \ref{theo-1.2zz}  where $u^\bb(x,y)$ is given in  (\ref{exe.13ja}) with $C=1/2$. This is the potential  density of Brownian motion,   killed after an independent exponential time     with mean $1/\bb$.
For any integer $k\ge 1$ let $  \wt Z_{k/2}=\{  \wt Z_{k/2}(t), t\in  R^{ 1}\}$ be a    $k/2$-permanental process with kernel  $\{u^\bb(x,y)+g(x)f(y),x,y\in R^1 \}$,  
where   $f,g$ are $C^{2}(R^1)$ excessive function for   $\ov B$.
Then   for any   $d\in R^{ 1}$ with $f'(d),g'(d)=0$, for all integers $k\ge 1  $,   
\begin{equation}
   \limsup_{  x \to 0}\frac{|  \wt Z_{k/2}(d+x)-  \wt Z_{k/2}  ( d)|}{ \(  4|x|\log\log 1/|x|\)^{1/2}} =   \(2   \wt Z _{k/2} ( d)\)^{1/2}, \qquad a.s.\label{1.34aa}
   \end{equation}
 
     Here $\wt Z_{k /2} (d)$ is a chi--square random variable of order k, as described in the paragraph following Theorem \ref{theo-1.2zz},  in which the   mean zero normal random variable has variance $ (2\bb)^{-1/2} +g(d) f (d) $.

 \medskip  Our initial interest in permanental processes was that their sample path properties can be used to obtain sample path properties of local times  of related Markov processes.  (See \cite{MRsuf}.) In  Sections 
    \ref{sec-emodMP}--\ref{sec-Borel} we give LILs for the local times of    Markov processes with potential densities that are the same as or closely related to kernels of the form of $u_{g,f}(x,y)$ in (\ref{1.3qssj}), and consequently, in different examples,  are equal to the kernels of the permanental processes in Theorems \ref{theo-1.2zz}--\ref{theo-1.5}.  
     These Markov processes are obtained by rebirthing a symmetric   transient Borel right process,   say \index{Borel right process} $Y$   with state space $ T$,  and   potential    densities  $u(x,y)$ with respect to some reference measure $m$, to obtain a new Markov process $\wt Y$. I.e.,  $\wt Y$ starts out as $Y$, but instead of going to a cemetery  state at the end of its lifetime it  is   ``reborn'' with    measure $\nu$, $| \nu|\le 1$, and remains in the cemetery  state with probability $1-|\nu|$.  If   the process is reborn it goes to the set $B\subset  T$ with measure $\nu(B) $, after which  it continues to evolve repeatedly   as we just described. When $| \nu|< 1$ we say the process  is partially reborn. In this case the process eventually dies out. When $| \nu|= 1$
    we say the process  is fully reborn.
    In Sections \ref{sec-LILpr} and \ref{sec-Borel} we use Theorems \ref{theo-1.2zz}--\ref{theo-1.5}
   to obtain LILs for the local times of partial and  fully   rebirthed   Markov processes.
    
  The following LIL is an example of the results the are obtained. It is a special case of LILs   in Example \ref{ex-7.1}:
  Let   $Y$ be a symmetric   L\'evy process killed at the end of an exponential time with mean $1/\bb$, as in Theorem \ref{theo-1.2zz} and let $\wt Y$ be a fully rebirthed Markov process determined by $Y$ and the probability measure $\mu$.  The process $\wt Y$ has potential densities of the form 
  \be
   u^{\bb} _{C, f}(x,y)= u^{\bb}(x-y)+Cf(y),\qquad x,y\in   R^1,\label{1.3qssjr}
  \ee 
for some  known constant constant $C$, where
\begin{equation}  
 f(y)=\int_{R^1} u^{\bb}(y-z)\,d\mu(z).  \label{1.3qssjr37}
\end{equation}
 (Such functions  are excessive for $Y$.)
    
     Let 
$   L=\{    L^y_{t}\,,\,(y,t) \in R^1\times R_+\}$
denote the local   times for
$\wt   Y$, normalized so that  
\begin{equation}
   E^{ x}\(\int_{0}^{\ff} e^{-\bb t} d_{t}  L^{y}_t\)=    u^{\bb} _{C, f}(x,y),\qquad x,y\in R^{ 1}. \label{az.1cq}
\end{equation}
 When   $f\in C^1(R^1)$   and $(i)$ and $(ii)$ in Theorem \ref{theo-1.2zz} hold, 
   $     L$ is jointly continuous on $R^1\times R_+$, and for all $v\in R^1$, 
\begin{equation}
 \limsup_{  u \to v}
\frac{ | L^{u}_t-  L^{v}_t|}{\(2   (\si^\bb)^{2}\(u-v\)\log\log 1/|u-v|\)^{1/2}} =
   \(  L^{v}_t\)^{1/2},\qquad \text{for a.e. $t\in R_+$}\label{az.3q},
\end{equation}
  $ P^{x}$ almost surely, for all $x\in R^1$.

 \medskip  Here is a more explicit  example of   LILs for the local times of a partially rebirthed diffusion.  
  Let ${\cal X}=\{{\cal X}(t),{t\in R_{+}^1 }\}$ be a diffusion
  with potential densities 
 \begin{equation} \label{page15}
  \wt u_{T_0}(x,y)=s(x) \wedge s(y), \qquad x,y>0,\end{equation}
 described in (\ref{diff.5}),   starting in $(0,\ff)$ and killed the first time it hits 0. 
 Let ${\wt{ \cal X} }=\{{\wt{ \cal X}(t)},t\in R_{+}^1 \}$ be a partial rebirthing of ${\cal X}$   determined by a measure $\mu$ as described in  Remark \ref{rem-6.1}. We take $\mu=h(x)\,dx$
 where $h(x)$ is positive and continuous on $[0,x_0]$ with   $h(x )=0$ on $[x_0,\ff)$.  The partially rebirthed process $\wt{ \cal X} $ has state space   $S=(0,\ff)\cup\{* \} $, where $*$ is an isolated point that $\cal X$ and subsequently  $\wt{ \cal X} $ goes to when it dies and awaits a possible rebirth. Let $\{\mathbf{\wt u}_{T_0}(x,y),x,y\in  S\}$ denote the potential density of $\wt{ \cal X} $. It is a function of $  \wt u_{T_0}(x,y)$ 
 and   $\mu$  as described in Theorem \ref{theo-borelpr}.   

\medskip Let  $\wt L=  \{ \wt L_{t}^{y}, (y,t)\in   S \times (0,\ff)\}$ be the local times   of $  \wt{ \cal X} $ normalized so that
 \begin{equation} \label{xzxq}
  E^x\wt L_{\ff}^{y}=\mathbf{\wt u}_{T_0}(x,y).
  \end{equation}
 Then for any $x_0 \in (0,\ff)$
\begin{equation}
   \limsup_{  x\to 0}\frac{|\wt   L_{t}^{x_0+x}- \wt  L_{t}^{x_0}|}{  (  2 s'(x_0)|x|\log\log 1/x)^{1/2}}=  \({2\wt  L_{t}^{x_0}}\)^{1/2} \qquad a.s.  \mbox{ for a.e. } t<\ze, \label{1.66mm}
   \end{equation} and 
  \begin{equation}
   \limsup_{  x\to 0}\frac{|   \wt L_{\ff}^{x_0+x}- \wt L_{\ff}^{x_0}|}{  (  2s'(x_0)  |x|\log\log 1/x)^{1/2}}=  \({2 \wt L_{\ff}^{x_0}}\)^{1/2} \qquad a.s.,    \label{1.67mm}
   \end{equation}
 $ P^{v}$ almost surely, for all $v\in (0,\ff)$.      For details see  Remarks \ref{rem-6.3} and \ref{rem-9.1}.  (This also holds for $ \limsup_{  x\uparrow 0}$).
 
  To see the dependence  of $\wt L_{\ff}^{x_0}$ on $\mu$ we note that,
\begin{equation} \label{xzx}
  E^v\wt L_{\ff}^{x_0}= s(v)\wedge s(x_0) + \int_0^{x_0} s(y) h(y)\,dy. 
  \end{equation}
  This integral can be expressed in a more probabilistic fashion  since it is equal to $\|h\|_1E(s(\th))$ where $\th$ is a  random variable with probability density $h/\|h\|_1$.

\medskip  Results on sample path properties of permanental processes are given in \cite{MRsuf, MRnec,MRejp}. The upper bounds in Theorems \ref{theo-1.2zz}--\ref{theo-1.5} are given in    \cite[Theorems  1.13]{MRejp}, and in fact hold for all $\al>0$, not just $\al=k/2$ for integers $k\ge 1$.
However, finding lower bounds for the local and uniform moduli of continuity of permanental processes  seems to be more difficult.  
     In \cite[Section 7]{MRejp} we  use a Slepian type lemma to get lower bounds for certain sequences of permanental process like $\{  X_{\al}(t_{j}),j\in \mathbb N \}$. However, we can not apply this lemma to increments  $\{  X_{\al}(t_{j})-  X_{\al}(0),j\in \mathbb N \}$. 	
    (The classical Slepian's lemma for Gaussian processes applies to increments because the increment  of a Gaussian process is a Gaussian process. We can not apply the  Slepian type lemma for permanental processes to increments of    a permanental process  because   increments of a permanental process  are generally not a permanental   process. This is easy to see. Let $\{X(t), t\in T\}$ be a permanental process. Let $t_{0}\in T$ and assume that $X(t_{0})\ne 0$. Then $\{X(t)-X(t_{0}), t\in T\}$ is not a permanental process because the terms of  a permanental process are all positive.)

  \medskip	   To prove Theorems  \ref{theo-1.2zz}--\ref{theo-1.5}  we  use a   comparison theorem between general permanental processes and permanental processes with symmetric kernels,  Theorem \ref{cor-7.1jq}, and show in Lemma \ref{theo-4.1} how to use it to relate $\{X_{k/2}(t),t\in[0,1] \}$ to chi-square processes so that we can use Gaussian techniques to study their sample path properties.    In Section \ref{sec-LIL} we prove Theorem \ref{theo-1.2z} which gives   lower bounds  for the LILs of a wide class of permanental processes with kernels that are not symmetric.  Theorems  \ref{theo-1.2zz}--\ref{theo-1.5} are proved in Section \ref{sec-4}.
 In Section \ref{sec-ex} and  \ref{sec-ex1} we give examples of potentials and excessive functions that satisfy the hypotheses of Theorems \ref{theo-1.2zz}--\ref{theo-1.5}.   

\medskip In Section    \ref{sec-emodMP} we consider Markov processes with potential densities that are  the same as  the kernels of  1/2--permanental processes.  In  Theorems \ref{theo-LILMP} and \ref{theo-5.2newa} we show that the local and uniform moduli of continuity 
of the local times of the Markov processes are the same as the local and uniform moduli of continuity 
of the 1/2--permanental processes. These theorems are used in Section \ref{sec-LILpr} to give LILs for the local times of partially rebirthed Markov processes with potential densities that are  the same as  the kernels of  the permanental processes in Theorems   \ref{theo-1.2zz}--\ref{theo-1.5}. These theorems are also used in Section \ref{sec-Borel} to give LILs for the local times of certain fully rebirthed Markov processes.

The hypotheses of Theorems \ref{theo-1.2zz}--\ref{theo-1.5} contain conditions on kernels of permanental processes and   excessive functions of related Markov processes. In Section \ref{sec-ex} we give many examples   when the conditions on the kernels are satisfied. We do the same in Section \ref{sec-ex1} for the excessive functions.

 \medskip  
 We thank Pat Fitzsimmons for his help on Sections   \ref{sec-emodMP} and \ref{sec-Borel}.

\section{Comparing general permanental sequences to symmetric  permanental sequences}

    Finding lower bounds for stochastic processes is a difficult problem. The only effective method   we  know is to use   Slepian's Lemma, which applies only to  Gaussian processes. Therefore  we need a relationship between permanental processes that are not Gaussian processes and Gaussian processes. In this paper we accomplish this by using  Theorem \ref{cor-7.1jq} below.

\medskip    Let  $C=\{ c_{ i,j}\}_{ 1\leq i,j\leq n}$ be an
$n\times n$ matrix. We call $C$ a positive matrix,
and write 
$C\geq 0$, if
$c_{ i,j}\geq 0$ for all
$i,j$.  

\medskip A matrix
$A$  is said to be a nonsingular  $M$-matrix   if
\begin{enumerate}
\item[(1)] $a_{ i,j}\leq 0$ for all $i\neq j$.
\item[(2)] $A$ is nonsingular and $A^{ -1}\geq 0$.
\end{enumerate}

 It follows from \cite[Lemma 4.2]{EK} that the right-hand side of   (\ref{int.1}) is a Laplace transform for all $\al>0$ if
  $A=K^{-1}$ exists and  is a  nonsingular $M$-matrix.   We refer to   $A$ as the  $M$-matrix  corresponding to $\wt X$. We also use the terminology that $K$ is an inverse $M$-matrix.

 \medskip	When $A$ is an $M$-matrix, we define 
\begin{equation}
   A_{sym}= \left\{
   \begin{array}{cl}
   A_{j,j}&j=1,\ldots,n\\
  -(A_{i,j}A_{j,i})^{1/2}&i,j=1,\ldots,n, i\ne j
   \end{array}  \right. .\label{2.15}
   \end{equation}
In addition, when $A=K^{-1}$ we define
\begin{equation}
     K_{isymi}= (A_{sym})^{-1}.\label{con7}
   \end{equation}
  (The notation $isymi$ stands for: take the inverse, symmetrize and take the inverse again.)    
 
  \begin{lemma}\label{lem-2.2} If $K$  is an inverse $M$-matrix, then so is     $K_{ isymi}$.    
\end{lemma}

\Proof By hypothesis $A=K^{-1}$ is a non-singular $M$-matrix. Therefore, by \cite[Lemma 3.3]{MRHD},  $A_{sym}$ is a non-singular  $M$-matrix. We have denoted it's inverse by $K_{isymi}$.   \qed

   \medskip	 We now show how we can relate probabilistic properties of $\al$-permanental processes that have kernels that are  not equivalent to a  symmetric kernel  to $\al$-permanental processes that have kernels that are    symmetric, and therefore, when $\al=k/2$, to   Chi--square processes of order $k$. The next  lemma combines  \cite[Corollaries 3.1 and 3.2]{MRHD}.

    \begin{theorem}\label{cor-7.1jq}  For any $\al>0$ let $Z_{\al}(n)=(Z_{\al,1},\ldots,Z_{\al,n})$ be the $\al$-perman-ental random variable determined by an $n\times n$ kernel $K(n)$ that is an inverse   $M$-matrix and set $A(n)=K(n)^{-1}$. Let $\wt Z_\al(n) $   be the $\al$-permanental random  variable determined by   $K(n)_{isymi}$.  Then for all functions $g_{n}$ of $Z_\al(n)$ and $\wt Z_\al(n)$ and sets $\BB_{n}$ in the range of $g_{n}$
\bea
\lefteqn{ \frac{|A(n)|^{\al}}{|  A(n)_{sym}|^{\al}}  P\(g_{n}( \wt Z_\al(n))\in \BB_{n}\)  \le  P\( g_{n}(   Z_\al(n))\in\BB_{n}\)\label{3.16w}}\\
 & & \hspace{1in}\le\(1- \frac{|A(n)|^{\al}}{|  A(n)_{sym}|^{\al}} \) +\frac{|A(n)|^{\al}}{|  A(n)_{sym}|^{\al}}P\(g_{n}( \wt Z_\al(n))\in\BB_{n}\).\nn
 \eea
 \end{theorem}
 We use this theorem with $\al=1/2$, in which case   $\wt Z_\al(n)$	is a  process of Gaussian squares.
  
  \medskip We now develop the material need to apply Theorem \ref{cor-7.1jq} to the $\al$-permanental processes considered in Theorems \ref{theo-1.2zz}--\ref{theo-1.2aa}.  We begin with the proof of 
  Theorem \ref{theo-borelN}.
     
\medskip \noindent \textbf{Proof of  Theorem \ref{theo-borelN} } By  the Kolmogorov extension theorem it suffices to show that for any finite subset  $T_{m}=\{t_{0}, t_{1},\ldots, t_{m}\}$ in $T$ there exists an $\al$--permanental vector, say $X_{\al}$, with kernel  
\be
  u_{g,f}(t_{i},t_{j})=u(t_{i},t_{j})+g(t_{i})f(t_{j}),\qquad t_{i},t_{j}\in  T_{m}.\label{bor.1}
  \ee 
  To do this we introduce another point $t_{-1}\notin T_m$ and  show that  there exists an $\al$-permanental vector with kernel,
 \bea  \label{bor.2a}
\begin{array}{ lcll }
 u_{g,f} (t_{-1},t_{-1})&=&1,& \\
  u_{g,f}(t_{i},t_{-1})&=&g( t_{i}),&      i= 0,\ldots,m,\\
  u_{g,f}(t_{-1},t_{j})&=& f(t_{j}),&   j=0,\ldots,m, \\
  u_{g,f}(t_{i},t_{j})&=&u(t_{i},t_{j})+g(t_{i})f(t_{j}), &  i,j=0,\ldots,m.
\end{array} 
\eea 
It is clear from the definition of permanental vectors in (\ref{int.1}) by its Laplace transform that if, say   $\wt Y$, is a permanental vector on a set $S$ then it is a permanental vector on subsets of $S$. Therefore, $X_{\al}$ is a permanental vector. 
  
\medskip Let  $\{U_{i,j},i,j=0,1,\ldots,m \}$ be an  $m+1\times m+1$ matrix with entries $U_{i,j}=u(t_{i},t_{j})$.   It follows from \cite[Lemma 4.2]{EK} that to show that (\ref{bor.2a}) is the kernel of an an $\al$--permanental vector, it suffices to show that the following $(m+2)\times (m+2)$
matrix is invertible and its inverse is an M--matrix:
  \bea 
 && K=\(
\begin{array}{ cccc } 1 &f( t_{0})&\ldots&f(t_{m})  \\
g( t_{0})   &U_{ 0,0}+g( t_{0}) f( t_{0})&\ldots&U_{0,m}+ g( t_{0}) f( t_{m})  \\
\vdots& \vdots &\ddots &\vdots  \\
g( t_{m})  &U_{m,0}+g( t_{m})f(t_{o} )&\ldots&U_{m,m}+ g( t_{m})f( t_{m}  )
\end{array}\)\!. \label{bor.3}
  \eea
 It is clear from (\ref{bor.3}), by subtracting $g( t_j)$ times  the first row from the $j+1-$st row  that,
  \begin{equation}
 |K |=  |U|.\label{bor.6ax}
  \end{equation}
  Since it follows from \cite[Lemma A.1]{MRnec} that $U$ is invertible, $|U|>0$.
Therefore $K$ is invertible.  

Let  
\be
A=K^{-1}.\label{bor.6ay}
\ee 
By multiplying the following matrix on the right by $K$ one can check that,   
  
\be  A = \left(\!\!\!
\begin{array}{ cccc  }  1+\rho&- v_{0}&\dots&-v_{m} \\
-r_{0}& U^{0,0} & \dots &  U^{ 0 ,m}    \\
\vdots&\vdots& \ddots&\vdots  \\
- r_{m}& U^{m ,1 }&  \dots & U^{m,m }   \end{array}\!\!\!\right)  \label{bor.7},
  \ee 
   where,
  \begin{equation} \label{bor.8}
  r_{i}=\sum_{ j=0}^{m} U^{i,j}g(x_j),\quad  v_{j}=\sum_{ i=0}^{m} U^{i,j} f( x_i),\quad i,j=1,\ldots,m,  
\end{equation}
and
\begin{equation}
\rho=\sum_{i, j=0}^{m} U^ {i, j} f(x_i) g(x_j).\label{bor.9}
\end{equation}
(We use the notation $U^{i,j}=(U^{-1})_{i,j}.$)

 It follows from \cite[Lemma A.1]{MRnec} that $U^{-1}$ is an M--matrix, therefore to show that $A $ is an M--matrix, it suffices to show that for all $ j=0,\ldots,m$,  $ r_{j}\ge 0$ and $v_{j}\geq 0$.   Since $U$ is symmetric   to show this it suffices to show that for any   finite excessive function   $h$ 
 for  $ \ov Y$
 \begin{equation}
 \sum_{ i=0}^{m} U^{i,j} h( x_i)\geq 0, \quad i,j=1,\ldots,m.\label{bor.9a}
 \end{equation}
 
  We obtain (\ref{bor.9a}) by considering an earlier version of this theorem in which $g\equiv 1$ and $f=h$ where 
 \begin{equation}
h(y)=\int u(x,y)\,d\mu (x),\label{bor.9x}
\end{equation}
 and $\mu$ is a finite measure on T.
 We show in  \cite[Theorem A.1]{MRejp} that   (\ref{bor.2a}),   with $g\equiv 1$, is the kernel of a transient Borel right process.  We show in \cite[Lemma 7.1]{MRnec} this implies that $A$ is an M--matrix. This implies (\ref{bor.9a}) for functions $h$ expressible as in (\ref{bor.9x}).
 Since, as explained in the proof of  \cite[Theorem 1.11]{MRejp},     all excessive functions   are    limits of such  functions, it follows that (\ref{bor.9a}) holds for all excessive functions $h$.\qed

 \medskip    Let   \label{calY}${\cal Y}\!=\!
(\Om,  \FF_{t}, {\cal Y}(t),\th_{t},P^x
)$ be a  symmetric   Borel right  process with state space    $T=R^1$ or $R^1-\{0\}$ and    continuous  potential densities $u=\{u(x,y),x,y\in T\}$ with respect to some locally bounded reference  measure   $m$ on $T$.  In either case we assume that $u$ has a continuous extension $u=\{u(x,y),x,y\in R^1\}$ and  take 
 $\eta=\{\eta_{x},x\in R^1\}$ to be a    Gaussian process with covariance $\{u(x,y);x,y \in R^1\}$.
Set
\begin{equation} \label{3.5nn}
\si^2(x,y):=  E(\eta_x-\eta_y)^2=u(x,x)+u(y,y)-2u(x,y).
\end{equation}
By (\ref{smp.1})
\begin{equation} \label{cond1}
  u(x,y)\le u(x,x)\wedge  u(y,y),\qquad\forall x,y\in R^1.
\end{equation}
   Then obviously,
   \begin{equation} \label{4.3new}
0\le   u(x,x)-u(x,y)\le \si^2(x,y) ,\qquad\forall x,y\in R^1.
\end{equation}
 
 Fix $d\in T$.   We point out in (\ref{poscond.8}) that  $u(d,d)>0$. Therefore,  it  follows from the continuity of $u$ that for some $ \de>0$  and function $b$, 
\begin{equation} \label{4.4new}
  \min_{x,y\in \De_d(\de)}u(x,y)=b(d, \de)>0.
\end{equation}
(See (\ref{deld}).) 

\medskip  
 
 We impose the following regularity conditions on the covariance $u$. We assume that for each $d\in T$ there exists  a continuous increasing function  $   \vf_d(x)$, $x\ge 0$, such that $\vf_d(0)=0$, and there exists a   $\de >0$, which can depend on $d$, such that,
  \bea   \label{3.7q}
 \si^2(x,y) &\asymp &
  \vf_d\( {| x-y|}\),\qquad\forall x,y\in  \De_d(\de),\,\, \text{and}, \\
    u(x,x)-u(x,y) &\ge &
 c_1 \vf_d\( {| x-y|}\),\qquad\forall x,y\in  \De_d(\de),\,\, x\ge y, \nn
\eea 
 for some constant $c_1>0$,  and
  \begin{equation} \label{2.149}
 \limsup_{x\to 0} \frac{|x|}{\vf_d(x)}=O(1) .
\end{equation}(The notation in (\ref{3.7q}) indicates that there exist constants, $0<c_2\le c_3<\ff$, such that, $ c_2
  \vf\( {| x-y|}\)\le \si^2(x,y) \le c_3
  \vf\( {| x-y|}\)$.)  
  
  We   generally simply write $\De_d$ or simply $\De$ instead of $\De_d(\de)$ and $\vf$ instead of $\vf_d$.

   \medskip Let $\{f(x),x\in T \}$ and $\{g(x),x\in T \}$ be    finite excessive functions    for  $ \cal Y$.   We show in  Theorem \ref{theo-borelN} that
\be
  u_{g,f}(x,y)= u(x,y)+g(x)f(y),\qquad x,y\in T,\label{1.33ssj}
  \ee     
is the kernel\index{kernel} of an   $\al-$permanental process. We now show how  to   relate path properties of this permenental process to path properties of a mean zero Gaussian process     with a covariance that is a minor modification of $u(x,y)$.

 \medskip   For all $n$ sufficiently   let  $\{t'_j \}_{j=0}^{m(n)}$ be  an increasing sequence with
 \begin{equation} 
  t'_0=d\ge 0,\qquad \text{and}\qquad \lim_{n\to\ff}t'_{m(n)}=t'_0.
\end{equation}   Consider the  Gaussian sequence   
  \be
  \eta(m(n), n)=\{\eta (t'_j); j= 0,1 ,\ldots,m(n) \},\label{gs.1j}
  \ee 
  with covariance,  
  \be
  u( m(n),n):=\{u(x,y);x,y\in t'_0 ,t'_1,\ldots, t'_{m(n)}\},
  \ee
 We set,
  \be u(  m(n),n)=  G(m(n), n)=\{G(m(n), n)_{j,k};j,k=0,1,\ldots, m(n)\}.\label{3.19}
\ee 
and assume that \begin{equation} 
 \sup_{n,m(n)} G(m(n), n)_{j,k}>0,\qquad \forall\,j,k=0,1,\ldots, m(n).\label{assume}
\end{equation}

  We show in    \cite[Lemma A.1]{MRnec}   that  $(G(m(n), n))^{-1}$
 exists and is an $M-$matrix. That is, it has positive diagonal entries and off diagonal entries that are less than or equal to 0.

To simplify the notation we write,
\begin{equation} \label{2.37mm}
  t_{j}'=d+t_j,\qquad j =0,1,\ldots, m(n),
\end{equation}
and to indicate where we are going we point out that in what follows 
 we choose $m(n)$ so that,  \begin{equation} \label{4.3w}
  \lim_{n\to\ff}m\(n\)=\ff,\qquad\text{and}\qquad\lim_{n\to\ff} m \(n\) t_{m(n)}=0. 
\end{equation}

 We also require that $f$ and $g$ have a continuous derivative on $[t' _{0}, t' _{m(n)}]=[d,d+t _{m(n)}]$.   If  $l$ is a function with this property then,  \begin{equation} \label{3.27nn}
(l')^{*}|  t_{j} - t_{k} | \le  |l(t'_{j})-l(t'_{k})|\le (l')^{**}|  t_{j} - t_{k} |,\end{equation}
where, \begin{equation} \label{2.30nn}
   (l')^{* }=\inf_{u\in [t' _{0}, t' _{m(n)}]}  l'(u)    \qquad\text{and} \qquad   (l')^{**}=\sup_{u\in [t' _{0}, t' _{m(n)}]}  l'(u)  .
\end{equation}
 In particular,  
  \begin{equation} \label{3.27rr}
   l(t'_{j})=l(d) +O(  t_{m(n)}).\end{equation}
  Therefore, both $f$ and  $g$ satisfy   (\ref{3.27rr}).

 \medskip Consider the $(m(n)+2)\times (m(n)+2)$ matrix given  by (\ref{bor.3}) with $t_j$ replaced by $t_j'$, $m$ replaced by $m(n)$ and $U_{j,k}$ replaced by $G(m(n),n)_{j,k}$. and denote it by $K(m(n), n)$.
    (To simplify the notation we  often denote $m(n)$ by $m$.) 
  Written out, 
  \bea 
 \lefteqn{K(m, n)\label{15.5}}\\
 && =\(
\begin{array}{ cccc } 1 &f( t' _{0})&\ldots&f( t'_{m })  \\
g( t'_0)   &G(m, n)_{ 0,0}+g( t'_0) f( t' _{0})&\ldots&G(m, n)_{0,m}+ g( t'_0) f( t'_{m })  \\
\vdots& \vdots &\ddots &\vdots  \\
g( t'_n)  &G(m, n)_{m,0}+g( t'_n)f(t' _{0} )&\ldots&G(m, n)_{m,m}+ g( t'_n)f( t'_{m })  
\end{array}\).  \nn
  \eea

As we point out in (\ref{bor.6ax}),
  \begin{equation}
 |K(m, n) |=  |G(m, n)|>0.\label{15.6}
  \end{equation}
   Therefore $K(m, n) $ is invertible.  
Define.  
\be
A(m, n)=K(m, n)^{-1},\label{15.6a}
\ee 
and note that,
\begin{equation} \label{3.11mm}
  |A(m, n)|=|G(m, n)|^{-1}.
\end{equation}
It follows from (\ref{bor.7}) that,

\be  A(m, n) = \left(\!\!\!
\begin{array}{ cccc  }  1+\rho_{m,n}&- v_{0}&\dots&-v_{m} \\
-r_{0}& G(m, n)^{0,0 } & \dots &  G(m, n)^{ 0 ,m}    \\
\vdots&\vdots& \ddots&\vdots  \\
- r_{m}& G(m, n)^{m ,0 }&  \dots & G(m, n)^{m ,m }   \end{array}\!\!\!\right)  \label{15.8},
  \ee 
  where,
  \begin{equation} \label{4.16}
  r_{j}=\sum_{ k=0}^{m}G(m, n)^{j,k}g(t'_k),\quad  v_{j}=\sum_{ k=0}^{m}G(m, n)^{j,k} f( t'_k),\quad j=0,\ldots,m,  
\end{equation}
and
\begin{equation}
\rho_{m,n}=\sum_{ j,k=0}^{m}G(m, n)^ {j,k} f(t'_j) g(t'_k).\label{rhomn}
\end{equation}
  It is clear from (\ref{3.11mm}) that $|A(m , n)|$  does not depend on $\{f(t'_j), g(t'_j) \}_{j=0}^{m}$.   
 
\medskip   By following  the proof of Theorem \ref{theo-borelN}, beginning with the second paragraph, we   see that
  $A(m, n)$ is an $M$--matrix and      $K(m, n)$ is the kernel of a   permanental process.

Set,\be
  \begin{array}{ll }
   \mathbf r=(r_0,\ldots,r_{m}), &\quad \mathbf v=(v_0,\ldots,v_{m})\\  \mathbf   f=(f(t'_0),f (t'_1)\ldots,f (t'_m)),  
  &  \quad \mathbf g=(g(t'_0),g (t'_1)\ldots,g (t'_m)) 
  \end{array}  
\ee  
  By   definition   $ f(t'_0)=f(d)$ and $ g(t'_0)=g(d) $.   Consequently, 
\begin{equation} 
  \mathbf f=(f(d), f (t'_1)\ldots,f (t'_m)),\qquad\text{and}\qquad       \mathbf g=(g(d), g (t'_1)\ldots,g (t'_m)). 
 \end{equation}
  We assume that    $g(d)$ and $f(d)$ are  strictly greater than 0.  
   
   In this notation, using (\ref{3.27rr}),    we can write (\ref{rhomn}) as  \begin{equation} \label{3.17a}
   \rho(m,n)=\sum_{k=0}^m v_kg (t'_k)  =(\mathbf{v}\cdot \mathbf{g})=(\mathbf{v}\cdot \mathbf{1})\(g(d)+O(t_{m(n)})\).
\end{equation}
 where $\mathbf{1}$ is a vector with all its components equal to 1.
  Similarly,   
 \begin{equation} \label{2.41re}
   \rho(m,n)=\sum_{j=0}^m r_jf(t'_j)=(\mathbf{r}\cdot \mathbf{f})=(\mathbf{r}\cdot \mathbf{1})\(f(d)+O(t_{m(n)})\).
\end{equation}
 
 \begin{lemma}  \label{lem-rowsum} 
\begin{equation}  
    \sum_{j=0}^m r_j= \frac{g(d)}{G(m,n)_{0,0}}+ O\(\vf (    t_{m(n)} )\),\label{3.23a}
\end{equation}
\begin{equation}  
    \sum_{j=0}^m v_j= \frac{f(d)}{G(m,n)_{0,0}}+ O\(\vf (    t_{m(n)} )\),\label{3.23b}
\end{equation}
and  
\be\label{3.23}
 \rho(m,n)= \sum_{j=0}^m v_jg(t'_j)= \frac{g(d)f (d)}{G(m,n)_{0,0}}+O\(\vf (    t_{m(n)} )\).
\end{equation} 
  \end{lemma}

\Proof   
Consider (\ref{15.8}) as a block matrix. We have,
 \begin{equation} \label{2.45mm}
  |A(m, n)|=|G(m, n)|^{-1}\(1+\rho_{m,n}-\mathbf vG(m, n)\mathbf r\).
\end{equation}
Therefore, by  (\ref{3.11mm}), 
\begin{equation} 
  \rho_{m,n}=\mathbf vG(m, n)\mathbf r,
\end{equation}
or, equivalently, by (\ref{3.17a})
\begin{equation} \label{3.18}
(\mathbf{v}\cdot \mathbf{g})= \sum_{j=0}^m g (t'_j) v_j= \mathbf vG(m, n)\mathbf r .
\end{equation}   
(This is  also easy to verify directly because, since $G(m,n)$ is symmetric, $\mathbf vG(m, n)\mathbf r=\mathbf rG(m, n)\mathbf v$ and $\mathbf rG(m, n)=\mathbf g$.) 

Let,
\begin{equation} 
  G(m,n)_{*}=\min_{j,k \in[0,m]} G(m,n) _{j,k}  \text{\,\, and \,\,}  G(m,n)_{**}=\max_{j,k \in[0,m]} G(m,n)_{j,k}.
\end{equation} 
  It follows from (\ref{3.18}) that 
  \begin{equation} \label{4.27}
G(m,n)_{*} \sum_{j=0}^m v_j \sum_{j=0}^m r_j\le   \sum_{j=0}^m v_j g (t'_j)\le  G(m,n)_{**} \sum_{j=0}^m v_j \sum_{j=0}^m r_j.
\end{equation}
  Using (\ref{3.27rr}) we see that,
 \begin{equation} \label{3.43nn}
\frac{\(g(d)+O(t_{m(n)})\)}{G(m,n)_{**}}\le  \sum_{j=0}^m r_j\le \frac{\(g(d)+O(t_{m(n)})\)}{G(m,n)_{*}}.
\end{equation}
The statement in (\ref{3.23a}) now follows from Lemma \ref{lem-5.1} and  (\ref{2.149}). The proof of the statement in (\ref{3.23b}) is similar. Using (\ref{2.41re}) and  (\ref{3.23a})      we get (\ref{3.23}).\qed

 \medskip  Using the fact that  $G(m, n)$ is symmetric we see that, 
\medskip\be \hspace{-4in} A(m, n)_{sym} 
\ee
\vspace{-.1in}
 \be = \left (
  \begin{array}{ ccccc  }   1+\rho_{m,n}&-h_0&-h_{1 } &\dots&-h_{m  } \\
- h_{0}& G(m, n)^{0 ,0 }& G(m, n)^{0 ,1 } & \dots &  G(m, n)^{0 ,m }   \\
\vdots&\vdots&\vdots& \ddots&\vdots  \\
- h_{m} & G(m, n)^{m ,0 }& G(m, n)^{m ,1 }&  \dots & G(m, n)^{m ,m}  \end{array}\right ) ,\nn
  \ee 
where,  
\be
h_{j  }=\(v_{j  }r_{j  }\)^{1/2},\qquad j=0,\ldots,m .
  \label{15.14}
\ee  
We  write $A(m, n)_{sym}$ in block form,   
\be A(m, n)_{sym} =\left (
\begin{array}{ ccccc }1+\rho_{m,n} &-{\bf  h} \\
-{\bf  h}^{T}&G(m, n)^{-1}  \end{array}\right ),\label{15.15}
  \ee 
where ${\bf  h} =(h_{0 }, h_{1  },\ldots,h_{m  })$.  Then,
\begin{equation}
|A(m, n)_{sym} |=|G(m, n)^{-1} |\,\,\(  1+\rho_{m,n}   -{\bf  h}G(m, n){\bf  h}^{T}\).\label{15.16}
\end{equation}
Using this and (\ref{3.11mm})  we see that,  
\begin{equation}
\nu_{m,n}={ |A(m, n)_{sym} | \over  |A(m, n)  |}=  1+\rho_{m,n} - {\bf  h}G(m, n){\bf  h}^{T}.\label{15.17}
\end{equation}

  To apply Theorem \ref{cor-7.1jq} we    give conditions under which  $ \lim_{n\to\ff}\nu_{m,n}=1$.
We do this in the next   two theorems.

\begin{theorem} \label{lem-3.3}     If $f,g\in C^1(\De_{t_0'} )$ and,  \begin{equation} \label{5.10bb}
\lim_{n\to\ff} \sup_{0\le j\neq k\le  m(n) }  m(n)  |t_k-t_j| \(\vf\( |t_k-t_j|\)\)^{-1}  = 0,
\end{equation}
and,
\begin{equation} \label{3.57qz}
  t_{m(n)}=o\(\frac{1}{m(n)}\),\qquad  \text{as $n\to\ff$,}
\end{equation}
  then
    \begin{equation} \label{4.37}
  \lim_{n\to\ff}\nu_{m,n}=1.
  \ee 
 
 \end{theorem}

\Proof   To   obtain (\ref{4.37}),  by  (\ref{15.17}), we must show that,
\begin{equation} \label{4.39}
 \lim_{n\to\ff} \(\rho_{m,n} - {\bf  h}G(m, n){\bf  h}^{T}\)=0,
\end{equation}
or equivalently, by (\ref{3.17a}), that,
\begin{equation} \label{4.40}
 \lim_{n\to\ff} \(\sum_{j=0}^{m}v_{j}g(t'_j)- {\bf  h}G(m, n){\bf  h}^{T}\)=0.
\end{equation}
(Recall that      $m=m(n)$. Also, we note that it follows from  \cite[Lemma 3.3]{MRHD}  that $\(\sum_{j=0}^{m}v_{j}g(t'_j)- {\bf  h}G(m, n){\bf  h}^{T}\)\ge 0$. 
 
      We now find a lower bound for ${\bf  h}G(m, n){\bf  h}^{T}.$ Using the fact that $G(m,n)^{-1}$ is an $M$  matrix we  see that,  
 \bea \label{3.67}
  g(t'_j)v_{j}&=&g(t'_j)\sum_{ k=0}^{m}G(m, n)^{j,k} f( t'_{k})\\
  &=&
  g(t'_j)\sum_{ k=0}^{m}G(m, n)^{j,k} f(t'_{j})-  g(t'_j)\sum_{ k=0}^{m}|G(m, n)^{j,k} |( f(t'_{k})-f(t'_{j})  )\nn\\
    &=&
  \sum_{ k=0}^{m}G(m, n)^{j,k}g(t'_k) f(t'_{j})-  g(t'_j)\sum_{ k=0}^{m}|G(m, n)^{j,k} |( f(t'_{k})-f(t'_{j})  )\nn\\
  &-&f(t'_{j}) \sum_{ k=0}^{m}|G(m, n)^{j,k} |( g(t'_{k})-g(t'_{j})  )\nn\\
 &\ge & r_j f(t'_{j})  - \( g^{**}(f')^{**}+f^{**}(g')^{**}\)  \sum_{ k=0}^{m} |t_{k} - t_{j}  | |G(m, n)^{j,k}|,  \nn
  \eea
   where we use (\ref{3.27nn}) for the last inequality and $f^{**}$ and $ g ^{**}$ are defined similarly to $(f')^{**}$ and $ (g') ^{**}$ in (\ref{2.30nn}).    Using  (\ref{3.43nn}), which, obviously, shows that the  terms $r_j $ are bounded, we have,\be
r_j f(t'_{j})\ge r_j\( f(d)- t_{j}(f')^{*} \)\ge r_j f(d)- Ct_{m(n)} (f')^{*}  ,\label{3.57}
\ee
  for some constant $C>0$. 
  
  It follows from (\ref{3.7q}) that
  \be \label{lower} G(m(n), n)_{j,j}-G(m(n), n)_{j,k}  \ge   
C \vf\( |t_j-t_k|\),\qquad\forall  j\ge k,\, 
\ee  
for some constant $C>0$.   Therefore we can use  
   Lemma \ref{lem-5.2}   below  to see that,  
   \be   \sum_{ k=0 }^{ m }  |t_k-t_j||G(m, n)^{j,k}|=o\(\frac{1}{ m(n)}\).\label{5.11b}
\ee
Consequently,  
 \begin{equation} \label{4.45a}
    g(t'_j)v_{j}\ge  r_j f (d)+o\(\frac{1}{m(n)}\)- Ct_{m(n)} (f')^{*}  .
\end{equation}
  Using (\ref{3.57qz}) we get,
 \begin{equation} \label{4.45a1}
    g(t'_j)v_{j}\ge  r_j f (d)+o\(\frac{1}{m(n)}\).
\end{equation}
Finally, using (\ref{3.27rr}) and the hypothesis that $g (d)>0 $ we have,
 \begin{equation} \label{4.45a2}
    v_{j}\ge  r_j \frac{f (d)}{g(d)}+o\(\frac{1}{m(n)}\).
\end{equation}

  We need a lower bound   for $ v_j^{1/2}$. This is a subtle point so it is useful to get a more precise description of the term $o\( {1}/{m(n)}\)$ in   (\ref{4.45a2}). We write (\ref{4.45a}) as,  
  \begin{equation} \label{4.45aa}
   v_{j}\ge r_j \frac{f (d)}{g(d)}+\frac{\ep(n,j)}{ m(n)},
\end{equation} 
where $\lim_{n\to\ff}\ep(n,j)=0.$   If $\ep(n,j)\ge 0$, $  v_j^{1/2}\ge (r_j f(d)/g(d))^{1/2}$. This is the case when $r_j=0$ because $v_j\ge 0.$
If $\ep(n,j)< 0$ and  $r_j>0$,   
 \begin{equation} \label{4.45aam}
  v_j^{1/2}\ge \(r_j \frac{f (d)}{g(d)}- \frac{|\ep(n,j)|}{m(n)}\)^{1/2}\ge  \(r_j \frac{f (d)}{g(d)} \)^{1/2} -\( \frac{|\ep(n,j)|}{m (n)} \)^{1/2}.
\end{equation}
where we use the fact that      $v_{j}\ge 0$.   
This implies that,
 \begin{equation} \label{4.49q}
  \sum_{j=0}^{m(n) } \(r_jv_j\)^{1/2}\ge  \sum_{j=0}^{m(n)} r_j \,\(\frac{f   (d)}{g (d)}\)^{1/2}  -  \( \frac{|\ep(n,j)|}{m (n)} \)^{1/2}\,\,\sum_{j=0}^{m(n)} r_j^{1/2}.
\end{equation}
Using  the inequality,
 \be 
  \sum_{j=0}^{m(n)} r_j^{1/2}  \le  m^{1/2}(n)\(\sum_{j=0}^{m(n)} r_j \)^{1/2} ,
\ee
and Lemma \ref{lem-rowsum}, which shows that $\sum_{j=0}^{m(n)} r_j<2g(d)/G(m,n)_{0,0}:=C$, for all $n$ sufficiently large, we see that, 
\begin{equation} \label{4.49}
  \sum_{j=0}^{m(n)} \(r_jv_j\)^{1/2}\ge  \sum_{j=0}^{m(n)} r_j \, \(\frac{f   (d)}{g (d)}\)^{1/2}  -C|\ep^*(n)|^{1/2},
\end{equation}
where $\ep^*(n)=\sup_j \ep(n,j)$.   Note that $\lim_{n\to\ff}\ep^*(n)=0$.

  \medskip We use (\ref{4.49}) to find a lower bound for ${\bf  h}G(m, n){\bf  h}^{T}$ in (\ref{15.16}).
 It follows from 
 Lemma \ref{lem-5.0} that,  
 \begin{equation} 
  G(m,n)_{j,k}=G(m, n)_{0,0}+O\(    \vf\( t_{m(n)}  \)\).
\end{equation}
 Therefore, by Lemma \ref{lem-5.0}
 \bea \label{3.64}
   {\bf  h}G(m, n) {\bf  h}^{T}&=&\sum_{j,k =0}^{m(n)}(r_jv_j)^{1/2}G(m, n)_{j,k}(r_kv_k)^{1/2}\\
   && =\(\sum_{j =0}^{m(n)}(r_jv_j)^{1/2}\)^{2} \( G(m, n)_{0,0}  +O\(    \vf\( t_{m(n)}  \)\)\).\nn
\eea  
It follows from (\ref{4.49}) that,  
\begin{equation} 
  \(\sum_{j =0}^{m(n)}(r_jv_j)^{1/2}\)^{2} \ge 
\(\sum_{j =0}^{m(n)} r_j  \)^{2}\frac{f (d)}{g(d)}+o_{n}(1).
\end{equation}
 Consequently,
\bea 
   {\bf  h}G(m, n) {\bf  h}^{T} &=& G(m, n)_{0,0}\(\sum_{j =0}^{m(n)} r_j  \)^{2}\frac{f (d)}{g(d)}+o_{n}(1)\\
   &=&   \frac{g(d)f (d)}{G(m,n)_{0,0}}+O\(\vf (    t_{m(n)} )\)+o_{n}(1),
  \nn
\eea
where we use (\ref{3.23a}). Therefore, by (\ref{3.23}),
\begin{equation} 
   {\bf  h}G(m, n) {\bf  h}^{T}=\rho_{m,n}+O\(\vf (    t_{m(n)} )\)+o_{n}(1),
\end{equation}
which gives (\ref{4.39}).\qed

    In the next theorem we  do not require   (\ref{5.10bb}),  but we    
  impose  more stringent conditions on the excessive functions $f$  and $g$.    
  
    \begin{theorem} \label{lem-3.3q}   When $f,g\in C^2(\De_d(\de))$, for some $\de>0$ and   $f'(d)=g'(d) =0$,  
 \begin{equation} \label{4.37ww}
  \lim_{n\to\ff}\nu_{m,n}=1.
  \ee 
  \end{theorem}
 
  \Proof      
Follow the proof of Theorem \ref{lem-3.3} through   (\ref{4.40}). We now modify   (\ref{3.67}) to  take account of the fact that $f'(d)=0$.  By Taylor's Theorem,   
\be
f(t'_{j})-f(t'_{k})=f'(t'_{k}) \( t_{j} - t_{k}\)  +\frac{1}{2} \int_{t'_{j}}^{t'_{k}}f''(s) (  t'_{k}-s )\,ds.
\ee 
 Since $f'(d)=0$, 
\begin{equation} 
 f'(t'_{k})=  \int_{d}^{t'_{k}}f''(s)  \,ds,  
\end{equation}
so that
 \begin{equation} \label{3.27nnq}
  |f(t'_{j})-f(t'_{k})|\le 2(f'')^{*}\,\,t_{m(n)} |  t_{j} - t_{k} |,\quad\text{where}\quad   (f'')^{*}=\sup_{u\in [d,d+t _{m(n)}]} |f''(u)| ,
\end{equation}
   where we take $n$ large enough so that  $t_{m(n)}<\de$.

Using Lemma \ref{lem-5.2ww}    we see     that, 
\bea 
  \sum_{ k=0}^{m}|G(m, n)^{j,k} |  f(t'_{k}) -f(t'_{j})   |&\le& O\(t_{m(n)} \sum_{ k=0 }^{m(n)}  |  t_{j} - t_{k} ||G(m, n)^{j,k}|  \)\nn \\
  &\le& O\(t_{m(n)}   \)\nn.
\eea  
 Similarly,
 \bea 
  \sum_{ k=0}^{m}|G(m, n)^{j,k} |  g(t'_{k}) -g(t'_{j})   |\le  O\(t_{m(n)}   \) .
\eea  
Using these last two inequalities in (\ref{3.67}) and then using (\ref{3.57})
we get,  
\be  \label{3.61z}
 g(t'_{j}) v_{j}\ge   r_j f(d)+O(t_m(n)).    
  \ee 
   It follows from (\ref{4.3w}) that (\ref{4.45a1}) holds.
  The rest of the proof proceeds exactly as in the proof of Theorem \ref{lem-3.3}.\qed

 To apply Theorem  \ref{cor-7.1jq} we must find $ K(m, n)_{isymi}$. By (\ref{15.15}) and the formula for  the inverse of  $A(m, n)_{sym}$  we see that,
\bea 
&& \hspace{-.4in}K(m, n)_{isymi}\label{8.29z}\\
&&=\left (
\begin{array}{ ccccc } \nu^{-1}_{m,n} &\nu^{-1}_{m,n} {\bf  h}G(m, n) \\
\nu^{-1}_{m,n} G(m, n){\bf  h}^{T}  &G(m, n)+\nu^{-1}_{m,n} G(m, n){\bf  h}^{T}{\bf  h} G(m, n) \end{array}\right ). \nn
  \eea 
For $i,j=0,1,\ldots,m$  set,
\begin{equation} \label{3.36mmz}
  (\nu^{-1}_{m,n} G(m, n){\bf  h}^{T}{\bf  h}G(m, n))_{i,j}=a(m,n)_{i } a(m,n)_{j },
\end{equation}
where,
\begin{equation} 
  a(m,n)_{j } =\nu^{-1/2}_{m,n}\({\bf  h}G(m, n)\)_{j} \label{15.24z}.
\end{equation} 
(We   use  the fact that $G(m, n)$
 is symmetric.) 
 
 Note that,  
 \begin{equation} 
 \{ (K(m, n)_{isymi})_{j,k};0\le j,k\le m\}\label{defrep}
\end{equation} 
 is the covariance of the Gaussian sequence,  
 \begin{equation} \label{defrepe}
 \wt\eta= \{\eta( t'_{j})+a(m,n)_{j }\, \xi;0\le j,k\le m\},
\end{equation}
where   $\{\eta (t'_{j})  ;j=0,\ldots,m\}$ is a mean zero Gaussian process with covariance $G(m,n)$ and $\xi$ is an independent standard normal random variable.  Furthermore,  
 \begin{equation} \label{3.46}
 0<  a(m,n)_{j}\le   f^{1/2} (d)g^{1/2} (d) +O\(\vf^{1/2}\( t_{m(n)}\)\),
\end{equation}
and,  
\be  \label{3.47z}
   |a(m,n)_{j }-a(m,n)_{k }|  \le  \frac{ f^{1/2} (d)g^{1/2} (d)}{G(m, n)_{0,0}}  O\(\vf^{1/2}\( |t_{j}-t_{k}|\) \).  
\ee 
 To verify (\ref{3.46})   note that,
by the Schwartz Inequality, Lemma \ref{lem-rowsum} and (\ref{4.4new}),
\begin{equation} \label{3.80z}
  \sum_{k=0}^{m(n)}h_k\le \(\sum_{k=0}^{m(n)}r_k \sum_{k=0}^{m(n)}v_k\)^{1/2} \le   \frac{ f^{1/2} (d)g^{1/2} (d)}{G(m, n)_{0,0}} +O\(\vf^{1/2}\(t_{m(n)}\) \).
\end{equation}
Using this, Lemma \ref{lem-5.1} and the fact that $\nu_{m,n}\ge 1$ we get (\ref{3.46}). 

\medskip  For (\ref{3.47z}) we have,
 \be  
    |\({\bf  h}G(m, n)\)_{j}-\({\bf  h}G(m, n)\)_{k}| \le \sum_{i=0}^{m}h_i|G(m, n)_{i,j}-G(m, n)_{i,k}|.     \ee 
  Note that by (\ref{3.7q}), 
  \bea 
  &&|G(m, n)_{i,j}-G(m, n)_{i,k}| =E(\eta_i(\eta_j-\eta_k)|\\
  &&\qquad \le \(  E(\eta^2_i)E(\eta_j-\eta_k)^{2}\) ^{1/2}\nn\le  \sup_{0\le i\le m}\(E(\eta^2_i)\)^{1/2}\vf^{1/2}(|t_j-t_k|).
\eea
Consequently, using (\ref{3.80z}), we see that,
    \bea \! |\({\bf  h}G(m, n)\)_{j}-\({\bf  h}G(m, n)\)_{k}| &\le&\!\! \sup_{0\le i\le m}\(E(\eta^2_i)\)^{1/2} \(\sum_{i=0}^{m}h_i\)  \vf^{1/2}(|t_j-t_k|) \nn\\
     &\le &\!\! \sup_{0\le i\le m}\(E(\eta^2_i)\)^{1/2} \frac{f^{1/2} (d)g^{1/2} (d)}{G(m, n)_{0,0}}     \vf^{1/2}(|t_j-t_k|). \nn\\
\eea
Since $ \sup_{0\le i\le m}\(E(\eta^2_i)\)^{1/2}$ is bounded we get  (\ref{3.47z}).

\medskip   The matrix $ K(m, n)_{isymi} $ is an $(m+2)\times (m+2) $ matrix which we can write as  
   $\{(K(m, n)_{isymi})_{i,j};	 i,j=-1,0,1,\ldots,m\}$. Define,   
  \be\label{3.41}
\ov\KK (m,n) =  \{(K(m, n)_{isymi})_{i,j};	 i,j=0,1,\ldots,m\}.
  \ee
  As we point out in (\ref{defrep}) and (\ref{defrepe}) this is the covariance of the Gaussian process,  
\begin{equation} \label{3.42}
 \{\eta (t_{j})+a(m,n)_{j }\,\xi;j=0,\ldots,m\},
\end{equation}
where   $\{\eta (t_{j})  ;j=0,\ldots,m\}$ is a mean zero Gaussian process with covariance $G(m,n)$ and $\xi$ is an independent standard normal random variable.

We also consider,
  \be
\KK (m,n) =  \{ K(m, n) _{i,j};	 i,j=0,1,\ldots,m\},
  \ee
in (\ref{15.5}).  
The next lemma is a restatement of  Theorem    \ref{cor-7.1jq} for the $\al$--permanental random variables with kernels $\KK_{(m,n)}$ and $\ov \KK_{(m,n)}$.

 	 \bl\label{theo-4.1}    Let  $\wh  X_{\al,(m,n)}= \{\wh X_{\al,(m,n)}(t'_j), j=0,1,\ldots m\}$  be an $\al$--permanental random variable with kernel $ \KK (m,n)   $  and   $  Y_{\al,(m,n)}= \{Y_{\al,(m,n)}(t'_j), j=0,1,\ldots m\}$   be an $\al$-permanental random variable with kernel $   \ov\KK (m,n)    $. 	Then for all functions  $g_{n} $ of $   \wh  X_{\al,(m,n)}$ and $  Y_{\al,(m,n)}$ and sets $\BB_{n} $ in the range of $g_{n} $,    
\bea
\lefteqn{   \nu_{n,m} ^{-\al}    P\(g_{n}(      Y_{\al,(m,n)})\in \BB_{n}\) \le P\( g_{n}(   \wh  X_{\al,(m,n)})\in\BB_{n}\)\label{3.16wq}} \qquad\\
 && \le  1-\nu_{n,m} ^{-\al}  \nn +   \nu_{n,m} ^{-\al}  P\(g_{n}(     Y_{\al,(m,n)})\in\BB_{n}\).\nn
 \eea
\el
 
 \Proof We use   Theorem \ref{cor-7.1jq} and (\ref{15.17}) with $\al$--permanental processes determined by the kernels $K(m, n)$ and $K(m, n)_{isymi}$.\qed

   \begin{remark}\label{rem-3.1}
{\rm It is clear from the properties of $\KK (m,n) $ and the Kolmogorov extension theorem  that there exists an $\al$-permanental process $\wh  X_{\al }= \{\wh X_{\al }(t), t\in [d,d+\de]\}$ such that $\{\wh X_{\al}(t'_j), j=0,1,\ldots m\}= \{\wh X_{\al,(m,n)}(t'_j), j=0,1,\ldots m\}$.  The process $\wh  X_{\al }$ is the $\al$--permanental process mentioned in the paragraph containing (\ref{1.33ssj}), where it is referred to as $X_{\al}$, with kernel
    $\{u_{g,f}(x,y)$, $x,y\in T\}$. In that paragraph we said that we would relate    path properties of $X_{\al}$ to path properties of a mean zero Gaussian process     with a covariance that is a minor modification of $u(x,y)$.  This is what is done in Lemma \ref{theo-4.1}. }
\end{remark}

  \subsection{Some basic calculations} 
  Recall that,  
\begin{equation} \label{5.1}
  G(m,n)_{*}=\min_{j,k \in[0,m]} G(m,n)_{j,k}  \,\text{ and } \, G(m,n)_{**}=\max_{j,k \in[0,m]} G(m,n)_{j,k}.
\end{equation} We often abbreviate $G(m,n)$ by just writing $G$. 

  \begin{lemma} \label{lem-5.0}
  \begin{equation} \label{5.1nn}
  G_{**}-G_{*}=  O\(\vf\(t_{m}\)\).
\end{equation} 
\end{lemma}
 \Proof The matrix $\{G_{j,k} \}_{j,k=0}^{m}$ has a finite number of terms, consequently the values $G_{**}$ and $G_{*}$ are achieved by some of the terms, and since $G_{j,k}\le G_{j,j}\wedge G_{k,k}$, the value $G_{**}$ is achieved by at least one of the diagonal terms. Therefore, there are terms such that, 
 \begin{equation} \label{5.2nn}
 0\leq   G_{**}-G_{*}=G_{j,j}-G_{k,l}.
\end{equation}
We note that, 
\begin{equation} 
  G_{j,j}-G_{k,l}=\( G_{j,j}-G_{j,k}\)+\( G_{j,k}-G_{k,l}\).
\end{equation}
  It follows from (\ref{3.7q}) that for all $t_{m(n)}$ sufficiently small,
\begin{equation} 
  G_{j,j}-G_{j,k}\le c_2\vf(|t_j-t_k|)\le c_2\vf(t_{m(n)}),
\end{equation} 
for some constant $c_2$.
 In addition, by  (\ref{smp.1}), 
\begin{equation} 
   G_{j,k}-G_{k,l} \le G_{k,k}-G_{k,l}\le c_2\vf(|t_{k}-t_{l}|)\le c_2\vf(t_{m}),\end{equation} 
 so we get (\ref{5.1nn}).\qed  

\begin{lemma} \label{lem-5.1}For all $\{ t_j \}_{j=0}^{{m(n)}}$,
  \begin{equation} \label{5.3new}
  \frac{1}{G(m,n)_{j,k}}=\frac{1}{G(m,n)_{0,0}}+ O\(\vf\(t_{m(n)}\)\) ,\qquad\text{as}\quad  n\to\ff. 
\end{equation}
\end{lemma}

\Proof   Write,
\begin{equation} 
  \frac{1}{G_{j,k}}=\frac{1}{G_{0,0}}+\frac{G_{0,0}-G_{j,k}}{G_{j,k}G_{0,0}}.
\end{equation}
This implies that, 
\begin{equation} 
\frac{1}{G^{}_{0,0}}-\frac{G_{**}-G_{*}}{G_{j,k}G_{0,0}}\le   \frac{1}{G_{j,k}}\le \frac{1}{G^{}_{0,0}}+\frac{G_{**}-G_{*}}{G_{j,k}G_{0,0}},
\end{equation}
 Using    (\ref{5.1nn}) and (\ref{4.4new}) we get  (\ref{5.3new}).  \qed

     \begin{lemma} \label{lem-5.2a} When 
\begin{equation}\label{3.lby} 
  \(G\(m ,n\)_{j,j}-G\(m ,n\)_{j,k}\)\ge \frac{2}{C}\vf\(|t_j-t_k|\),\quad  \forall    \, k<j,\,\,  j,k=0,\ldots, m(n) 
\end{equation}
for some $C>0$ and all $n$ sufficiently large it follows that,
   \begin{equation} \label{5.aa} 
   \sum_{k=0 }^{m(n)}   \big |G\(m ,n\)^{k,j}\big | \, \vf\(|t_j-t_k|\) \le C,\quad  \forall    \,\,  j,k=0,\ldots, m(n), 
\end{equation}
and all $n$ sufficiently large.
\end{lemma}

\Proof We have,  
\bea \label{3.66}
1&=& \sum_{k=0}^{m} G\(m ,n\)_{j,k} G\(m ,n\)^{k,j}\\
&=&G\(m ,n\)_{j,j} G\(m ,n\)^{j,j}-\sum_{k=0,k\ne j}^{m} G\(m ,n\)_{j,k} |G\(m ,n\)^{k,j}|\nn 
\eea
Continuing  this, 
\bea \label{5.4}
&=&G\(m ,n\)_{j,j} \(G\(m ,n\)^{j,j}-\sum_{k=0,k\ne j}^{m}   |G\(m ,n\)^{k,j}|\) \\
&&\qquad + \sum_{k=0 }^{m}   |G\(m ,n\)^{k,j}| \(G\(m ,n\)_{j,j}-G\(m ,n\)_{j,k}\)\nn\\
&=&G\(m ,n\)_{j,j} r_j+  \sum_{k=0 }^{m}   |G\(m ,n\)^{k,j}| \(G\(m ,n\)_{j,j}-G\(m ,n\)_{j,k}\).\nn 
\eea
 Consequently,
\begin{equation} \label{5.11t}
  \sum_{k=0 }^{m}   |G\(m ,n\)^{k,j}| \(G\(m ,n\)_{j,j}-G\(m ,n\)_{j,k}\)\le 1.
\end{equation}

 Furthermore, since $\(G\(m ,n\)_{j,j}-G\(m ,n\)_{j,k}\)\geq 0$ for all $j,k$ it follows that 
\begin{equation} \label{5.11ty}
  \sum_{k<j }   |G\(m ,n\)^{k,j}| \(G\(m ,n\)_{j,j}-G\(m ,n\)_{j,k}\)\le 1.
\end{equation}
Using  (\ref{3.lby}) we see  that  
 \begin{equation} \label{511}
 \sum_{k<j }  \big|G\(m ,n\)^{k,j}\big| \,\vf\(|t_j-t_k|\) \le C/2.
\end{equation}
This is true for all   $j$, so in particular, setting $j=k$, we have,
 \begin{equation} \label{512}
 \sum_{i<k}  \big|G\(m ,n\)^{i,k}\big| \,\vf\(|t_k-t_i|\) \le C/2,
\end{equation}
or, equivalently,
 \begin{equation} \label{512j}
 \sum_{j<k}  \big|G\(m ,n\)^{j,k}\big| \,\vf\(|t_k-t_j|\) \le  C/2.
\end{equation}
Since  $\{G\(m ,n\)^{j,k} \}_{j,k=0,\ldots,m(n)}$ is symmetric this shows that,
 \begin{equation} \label{513}
 \sum_{j<k}  \big|G\(m ,n\)^{k,j}\big| \,\vf\(|t_j-t_k|\) \le C/2.
\end{equation}
 combining (\ref{511}) and (\ref{513}) gives (\ref{5.aa}).\qed

 Recall the condition in (\ref{5.10bb}),  
\begin{equation} \label{5.10bbz}
\lim_{n\to\ff} \sup_{0\le j\neq k\le  m(n) }  m(n)  |t_k-t_j| \(\vf\( |t_k-t_j|\)\)^{-1}  = 0.
\end{equation}

 \begin{lemma} \label{lem-5.2}   When   (\ref{3.lby}) and (\ref{5.10bbz}) hold,  
\be   \sum_{ k=0 }^{ m(n)}  |t_k-t_j||G(m, n)^{j,k}|=o\(\frac{1}{ m(n)}\).\label{5.11}
\ee
\end{lemma}

\Proof    Using Lemma \ref{lem-5.2a} and (\ref{5.10bbz}) we get,
\bea
\lefteqn{    \sum_{ k=0 }^{m(n)}  |t_k-t_j||G(m, n)^{j,k}|}\\
&& \le  \frac{1}{ m(n)} \sum_{ k=0 }^{m(n)}   |G(m, n)^{j,k}|\vf\( |t_k-t_j|\)\left\{ m(n)|t_k-t_j| \( \vf\( |t_k-t_j|\)\)^{-1}\right\}\nn \\
&& \le  \frac{C}{m(n)} \sup_{0\le j,k\le m(n) } m(n)  |t_k-t_j| \(\vf\( |t_k-t_j|\)\)^{-1} = o\(\frac{1}{m(n)}\).\nn
\eea
 \qed

The next lemma gives an estimate that is weaker than (\ref{5.11}) but, significantly,  it does not require   (\ref{5.10bbz}):  
  \begin{lemma} \label{lem-5.2ww}   \be   \sum_{ k=0}^{ m(n)}  |t_k-t_j||G(m, n)^{j,k}|=O\(1\).\label{5.11a}
\ee

\end{lemma}

\Proof We have,
\bea 
    &&   \sum_{ k=0}^{ m(n)}  |t _k-t _j||G(m, n)^{j,k}| \\
    &&\qquad  \le   \sum_{ k=0}^{m(n)}   |G(m, n)^{j,k}|\vf\( |t_k-t_j|\)\left\{  |t_k-t_j| \( \vf\( |t_k-t_j|\)\)^{-1}\right\}.\nn
\eea
 It follows from  (\ref{5.aa}) and (\ref{2.149}) that this last term is $O(1)$.\qed

   \section{Lower bounds for the LIL for $\mathbf {k/2}$ permanental processes}\label{sec-LIL}  
   
  Upper bounds for the LILs in Theorems \ref{theo-1.2zz}--\ref{theo-1.5} are given in  \cite[Theorem 1.2]{MRejp}. In this paper we show that they are best possible by finding the lower bounds. We then say that the moduli are exact. I.e., Theorems \ref{theo-1.2zz}--\ref{theo-1.5} give exact local moduli of continuity for the permanental processes they consider. 

The next theorem is a general result for the lower bound of the local modulus of continuity for certain permanental processes. It is used in Section   \ref{sec-4} to find lower bounds in Theorems \ref{theo-1.2zz}--\ref{theo-1.5}.

    \begin{theorem}\label{theo-1.2z}        Let $\cal Y$, $u$, $\eta$ and $\si^2(x,y)$ be as described on page \pageref{calY} 
 and consider the permanental processes $X_{k/2}=\{X_{k/2}(t),t\in T\}$, $k\ge 1$,  determined by   kernels 
\be
  u_{g, f}(x,y)= u(x,y)+g(x)f(y),\qquad x,y \in T,\label{4.1nn}
  \ee 
  where   $g, f $ are excessive functions for   $\cal Y$.
 Fix $d\in T$ and assume that   $u(d,d)>0 $  and  $g, f\in C^{1}(\De_d(\de))$ for some $\de>0$; see (\ref{deld}).
 Furthermore, assume there exists a function that   $\vf_d$ that satisfies (\ref{3.7q}) and (\ref{2.149}) with respect to $\si^2(x,y)$
 
For   $\th<1$ let,  
\begin{equation} 
t_{0}=0\hspace{.1 in}\mbox{ and }\hspace{.1 in}  t_j=\th^{n+1-j},\,\, j=1,\ldots,m(n),
\end{equation}
and,
\begin{equation} 
    m(n)=n+1-[n^{q}]   \text{ for some $0<q<1$}. 
\end{equation}
Set 
\begin{equation}
 t'_j=d+t_j,\,\, j=0,1,\ldots,m(n).\label{d4.1}
\end{equation}
 Assume that    for all $\ep>0$,   
 \be  \label{3.5mm}P\(\sup_{1\le j\le m(n)} \frac{\eta(t'_j)-\eta(d)}{(2\si^2(t'_j,d)\log\log 1/t_j)^{1/2} }\ge  (1-\ep) \)  \ge  (1-\ep)   ,
\ee
for all $n$ sufficiently large, and that one, or both, of the following conditions are satisfied:
\begin{itemize} 
  \item[(i)]  
  \begin{equation} \label{5.10bbqv}
\lim_{n\to\ff} \sup_{0\le j\neq k\le  m(n) }  m(n)  |t_k-t_j| \(\vf_d\( |t_k-t_j|\)\)^{-1}  = 0;
\end{equation}
\ 
  \item[(ii)]     $g, f\in C^{2}(\De_d )$  and $g'(d)=f'(d)=0$.    
\end{itemize}
Then   
      for all integers $k\ge 1  $,  
\begin{equation}
   \limsup_{  t \downarrow 0}\frac{|  X_{k/2}(d+t)-  X_{k/2}  (d)|}{\( 2\si^{2}(d+t,d)\log\log 1/t\)^{1/2}}\ge   {\sqrt 2 }  X ^{1/2}_{k/2} (d) \qquad a.s.\label{1.34sz}
   \end{equation}
 This also holds for $\limsup_{  t \uparrow 0}$.\end{theorem}

  \noindent \textbf{Proof of  Theorem \ref{theo-1.2z}  }   We use Theorem \ref{cor-7.1jq} to find sample path properties of permanental processes with kernels that are not symmetric in terms of permanental processes with symmetric kernels. Since the later are chi-square processes they can be analyzed as squares of  Gaussian processes using well known techniques. To get almost sure results we must show that,
   \begin{equation} \label{nu}
   \lim_{n\to\ff}\nu_{m,n}=1.
\end{equation}
  (See (\ref{15.17}).)  When condition $(i)$ is satisfied Theorem \ref{lem-3.3} shows that (\ref{nu}) holds. When condition  $(ii)$ is satisfied Theorem \ref{lem-3.3q} shows that (\ref{nu}) holds.   (These are the only uses of conditions $(i)$ and $(ii)$ in the proof.)
    
    We proceed with the proof of this theorem: 
     Set 
\begin{equation} 
  T_n=\{t_0,t_1,t_2,\ldots, t_{m(n)}\}=\{0,  \th^n,\th ^{n-1},\ldots,\th^{[n^{q}]}\},\label{4.62ms}
\end{equation} 
and  
\begin{equation}
T'_n=d+T_n:=\{d, t'_1,t'_2,\ldots, t'_{m(n)}\}.\label{3.10q}
\end{equation}
  Also,  excluding  the initial element $d$, we set, 
\begin{equation}
T_{n,1}=\{t_1,t_2,\ldots, t_{m(n)}\}, \qquad \text{and}\qquad T'_{n,1}=\{t'_1,t'_2,\ldots, t'_{m(n)}\}.
\end{equation} 
We  use Lemma \ref{theo-4.1} to get probability estimates for $X_{k/2,(m,n)}=\{X_{k/2}(t),t  \in T'_{n}\}$,  in terms of a chi-square process of order   $k$  that is determined by a Gaussian sequence with covariance  $\ov\KK_{(m,n)}$    in (\ref{3.41}). Specifically,
    \begin{equation} \label{consistent.1}
  (\ov\KK_{(m,n)})_{j,k}=G(m,n)_{j,k} +a(m,n)_{ {j}}\,\,a(m,n)_{k }:=\ov u(t'_j, t'_{k}),
\end{equation}
where,
\begin{equation}
G(m(n), n)_{j,k}= u(t'_j, t'_{k}).\label{4.3wxyz}
\end{equation}
(See (\ref{3.19}).)
This is the covariance of the Gaussian sequence   $\ov\eta=\{ \ov\eta(t'_j); t'_j\in T'_{n}\}$ where,  
\begin{equation} \label{d3.42}
\ov\eta({t'_{j}}) =  \eta({t'_j})+a(m,n)_{j }\,\xi,\qquad j= 0,1,\ldots,m ,
\end{equation}
where   
$\xi$ is an independent standard normal random variable.

\medskip Let $\ov u$ denote the covariance of $\ov \eta$.   It follows from  (\ref{3.46}) that    $\ov u(d,d)>0.$  Using this,   (\ref{3.7q}), (\ref{3.46}) and (\ref{3.47z}) we see that,
\be 
  \frac{\ov u (d,d)-\ov u (d,t'_j)}{\ov u (d,d)} \le  C\vf^{1/2}\(t_j\) ,\label{4.17}
\ee 
  for some constant $C.$   Let  
  \begin{equation} \label{3.15q}
\psi_d(t)=\( 2\si^{2}(d+t,d)\log\log 1/t\)^{1/2}.\end{equation}
It  follows from    (\ref{2.149})  that, 
 \begin{equation}
 \frac{\vf_d^{1/2}\(t\)}{ \psi_d(t) }\to 0,\qquad\mbox{ as }t\to 0 .\label{link1}
\end{equation}
(We often   write $\vf_d$ and $\psi_d$ simply as $\vf$ and $\psi$.)

\medskip Let $\{\ov\eta_i,i=1,\ldots,k\}$ be independent copies of the Gaussian sequence $\ov\eta$.  
 Consider the chi-square process of order $k$ 
 \begin{equation} 
  \sum_{i=1}^{k} \ov \eta^{2}_{i }(t'),\qquad t'\in T'_n. 
\end{equation}
Note that, 
\begin{eqnarray}
 \ov \eta_{i}^{2}(t')-\ov  \eta_{i}^{2}(d)&=&( \ov \eta_{i}(t')- \ov \eta_{i}(d))( \ov \eta_{i}(t')+\ov  \eta_{i}(d))
\label{16.40mn}\\
&=& (\ov  \eta_{i}(t')-\ov  \eta_{i}(d))( 2\ov  \eta_{i}(d)+( \ov \eta_{i}(t')- \ov \eta_{i}(d)))   \nonumber\\
&=& 2( \ov \eta_{i}(t')- \ov \eta_{i}(d)) \ov  \eta_{i}(d)+(\ov  \eta_{i}(t')-\ov  \eta_{i}(d))^{2}.  \nonumber
\end{eqnarray}
Consequently,  
\begin{equation}
 \sup_{t'\in T'_{n,1}}{\sum_{i=1}^{k}(\ov \eta^{2}_{i }(t')-\ov\eta^{2}_{i }(d)) \over   2 \psi (t)}\ge \sup_ {t'\in T'_{n,1}}{\sum_{i=1}^{k}( \ov\eta _{i }(t')-\ov\eta _{i }(d))\ov\eta _{i }(d) \over    \psi (t)}\label{5.4mm3},
   \end{equation}
  where, by (\ref{3.10q}),
\begin{equation} \label{3.20q}
  t=t'-d.
\end{equation}     We write this as, 
 \bea
 \sup_{t'\in T'_{n,1}}{\sum_{i=1}^{k}(\ov \eta^{2}_{i }(t')-\ov\eta^{2}_{i }(d)) \over   2 \psi (t)}&\ge&  \nn \sup_{t'\in T'_{n,1}}{\sum_{i=1}^{k}\( \ov\eta_{i}(t')-\frac{\ov u(d,t')}{\ov u(d,d)} \ov\eta_{i}(d)\)  \ov\eta_{i}(d) \over  \psi (t)}\\
 &&\,\, -\sup_{t'\in T'_{n,1}}{ \(\frac{\ov u(d,d)-\ov u(d,t')}{\ov u(d,d)}\)\sum_{i=1}^{k} \ov\eta^{2}_{i}(d) \over    \psi (t)} . \label{4.20} \eea
Set,
\begin{equation} 
\de_{m,n} = \sup_{t'\in T'_{n,1}} \frac{\ov u(d,d)-\ov u(d,t')}{\ov u(d,d)\psi(t)},
\end{equation}
and note that by (\ref{4.17}) and (\ref{link1}), 
 \begin{equation} 
\lim_{n\to\ff}\de_{m,n}=0 .
\end{equation}
  With this notation we rewrite (\ref{4.20}) as,  
\bea
 && \sup_{t'\in T'_{n,1}}{ \sum_{i=1}^{k}( \ov\eta^{2}_{i }(t')-\ov\eta^{2}_{i }(d)) \over    2\psi(t) \| \vec{\ov\eta}(d)\|_2 }+\de_{m,n} \| \vec{\ov\eta}(d)\|_2  \label{3.43mm}\\
 &&\qquad \ge\nn \sup_{t'\in T'_{n,1}}{\sum_{i=1}^{k}\( \ov\eta_{i}(t')-\frac{\ov u(d,t')}{\ov u(d,d)} \ov\eta_{i}(d)\)   \over  \psi (t)}\frac{\ov\eta_{i}(d) }{\| \vec{\ov\eta}(d)\|_2} ,\nn 
 \eea
  where for $k$ dimensional Gaussian vectors we use the notation,    $\vec{\ov \eta}(t)=\vec{\ov \eta}_k(t)=(\ov\eta_{1}(t),\ldots,\ov\eta_{k}(t))$, for $t\in R^1.$ 
 
\medskip  The next step in the proof is the critical observation. We give it as a lemma.

 \begin{lemma}   For all $t'\in R^1$,
 \be 
   \sum_{i=1}^{k}\( \ov\eta_{i}(t' )-\frac{\ov u (d,t')}{\ov u (d,d)} \ov\eta_{i}(d)\)    \frac{\ov\eta_{i}(d) }{\| \vec{\ov\eta}(d)\|_2} \stl   \(\ov\eta (t')-\frac{\ov u(d,t')}{\ov u(d,d)} \ov\eta (d)\).  \label{3.24} 
\ee
\end{lemma}

\Proof  Let $\{\xi_{i}(t'),t'\in T_{n}\}$, $i=1,\ldots,k$, be independent copies of a mean zero Gaussian process  $\{\xi (t'),t'\in T'_{n}\}$, and   set  $\vec \xi(t')=(\xi_{1}(t'),\ldots,\xi_{k}(t'))$. Let $\vec v\in R^{k}$ with $\|\vec v\|_{2}=1$. Computing the covariances we see that,
\begin{equation}
  \{ \vec \xi(t')\cdot \vec v),t'\in T'_{n}\}\stl \{\xi(t'),t'\in T'_{n}\}\label{4.2dmm}.
   \end{equation}
Therefore, since 
   $\( \ov\eta_{i}(t')-\frac{\ov u (d,t')}{\ov u (d,d)} \ov\eta_{i}(d)\) $  and $\ov \eta_{i}(d)$ are independent for $i=1,\ldots,k$, we see that,
   \bea
&& \Big  \{   \(\vec{\ov \eta} (t')-\frac{\ov u (d,t')}{\ov u (d,d)}\vec{\ov \eta} (d)\)\cdot  \frac{\vec{\ov \eta} (d)}{ \|\vec{\ov \eta} (d)\|_{2}},t'\in T'_{n}\Big\}\label{4.2dmn}\\
&&\qquad\stl  \Big  \{   \( \ov\eta (t')-\frac{\ov u (d,t')}{\ov u (d,d)} \ov\eta (d)\),t'\in T'_{n}\Big\}.\nn
      \eea
  This gives(\ref{3.24})  
      \qed

It follows from (\ref{3.24})  that the second line of (\ref{3.43mm}) is equal in law to  
\bea 
   \sup_{t'\in T'_{n,1}}{ \(\ov\eta (t')-\frac{\ov u(d,t')}{\ov u (d,d)} \ov\eta (d)\)  \over  \psi (t)}&=&\sup_{t'\in T'_{n,1}}\({ \( \ov\eta (t')-  \ov\eta (d)\)  \over  \psi (t)}+{    \ov u (d,d)-\ov u (d,t')    \over \ov u (d,d) \psi (t)}\ov\eta(d)\)\nn\\
 &\ge& \sup_{t'\in T'_{n,1}}{ \( \ov\eta (t')-  \ov\eta (d)\)  \over  \psi (t)}- \de_{m,n}|\ov\eta(d)|\label{3.46mm}.
\eea
 Therefore, using  (\ref{3.43mm}) and (\ref{3.46mm}) we see that for any $  \ep>0$,  
    \bea
 &&P\(  \sup_{t'\in T'_{n,1}}{ \sum_{i=1}^{k}( \ov\eta^{2}_{i }(t')-\ov\eta^{2}_{i }(d)) \over   2 \psi(t) \| \vec{\ov \eta} (d)\|_2 }\ge  1-  \ep-\de_{m,n} \| \vec{\ov \eta} (d)\|_2    \) \label{5.14mmz}\\
 &&\qquad \ge P\(  \sup_{t'\in T'_{n,1}}{ \(\ov \eta (t')- \ov \eta (d)\)   \over  \psi (t)} \ge 1-  \ep- \de_{m,n}  |\ov\eta(d)| \)\nn .
   \eea
 Set, 
\begin{equation} 
  P\(\de_{m,n} \|    \vec{\ov \eta} (d)\|_2   \le \ep/2\)= 1-\ep(n). 
\end{equation} 
Then conditioning on the event $\{\de_{m,n} \| \vec{\ov \eta} (d)\|_2  \le \ep/2 \}$ we see that the left-hand side of (\ref{5.14mmz}),
 \begin{equation} 
  \le P\(  \sup_{t'\in T'_{n,1}}{ \sum_{i=1}^{k}( \ov\eta^{2}_{i }(t')-\ov\eta^{2}_{i }(d)) \over   2 \psi(t)\| \vec{\ov \eta} (d)\|_2 }\ge  1-  \ep/2 \)(1-\ep(n)) +\ep(n).
\end{equation}
Therefore,
 \bea \label{4.30nn}
 \lefteqn{P\(  \sup_{t'\in T'_{n,1}}{ \sum_{i=1}^{k}( \ov\eta^{2}_{i }(t')-\ov\eta^{2}_{i }(d)) \over   2 \psi(t) \| \vec{\ov \eta} (d)\|_2 }\ge  1-  \ep/2 \)}\\
  &&\quad \ge   P\(  \sup_{t'\in T'_{n,1}}{ \(\ov \eta (t')- \ov \eta (d)\)   \over  \psi (t)} \ge 1-  \ep- \de_{m,n}  |\ov\eta(d)| \)\nn -\frac{ \ep(n)}{1-\ep(n)}.
\eea
    By (\ref{d3.42}) for all $t'_{j}\in T'_{n,1}$,
  \begin{equation}\label{4.30mm}
 \frac{\ov\eta (t'_{j})-\ov\eta   (d)} {  \psi(t_{j})}  =    \frac{\eta( t'_{j})-\eta (d)} {\psi(t_{j})}+\frac{a(m,n)_{j}-a(m,n)_{0 } }{ \psi(t_{j})}\xi.
 \end{equation}
Consequently, 
 \bea
  \lefteqn{P\(  \sup_{t'\in T'_{n,1}}{ \(\ov \eta (t')- \ov \eta (d)\)   \over  \psi (t)} \ge   1-  \ep -   \de_{m,n} |\ov\eta(d)|\)}\\
   &&   \ge  P\( \sup_{t'\in T'_{n,1}}{ \(  \eta (t')-   \eta (d)\)   \over  \psi (t)}\ge   1-  \frac\ep2-  \de_{m,n} |\ov\eta(d)|  ,\right. \nn\\
   &&\hspace{1in}  \left.  \sup_ {1\le j\le m}\frac{|a(m,n)_{ {j}}-a(m,n)_{0 }| }{ \psi(t_{j})}|\xi|\le \frac\ep2 \)\nn 
  \eea
\vspace{-.2in}  \bea && \quad =     P\( \sup_{t'\in T'_{n,1}}{ \(  \eta (t')-   \eta (d)\)   \over  \psi (t)}\ge   1-  \frac\ep2-  \de_{m,n} |\ov\eta(d)|  \)\nn\\
   &&\hspace{.9in}P\(    \sup_ {1\le j\le m}\frac{|a(m,n)_{ {j}}-a(m,n)_{0 }| }{ \psi(t_{j})}|\xi|\le \frac\ep2  \).\nn 
    \eea 
By (\ref{3.47z}) and (\ref{link1}), for all $n$ sufficiently large,   
 \begin{equation}
  \sup_ {1\le j\le m(n)} \frac{|a(m,n)_{ {j}}-a(m,n)_{0 }| }{ \psi(t_{j})}|\xi|
\le     \de'_{m,n} |\xi|,
    \end{equation}
 Set,
\begin{equation} 
  P\(\de'_{m,n}  |\xi| \le \ep/2\)= 1-\ep'(n), 
\end{equation}  
so that,
 \begin{equation}\label{4.33}
    P\( \sup_ {1\le j\le m(n)} \frac{|a(m,n)_{ {j}}-a(m,n)_{0 }| }{ \psi(t_{j})}|\xi|
 \le  \frac{\ep}{2 }  \)= 1-   \ep'(n).
    \end{equation}
  Combining (\ref{4.30mm})--(\ref{4.33}) we have, 
 \bea\label{4.37qq}
  \lefteqn{P\(  \sup_{t'\in T'_{n,1}}{ \(\ov \eta (t')- \ov \eta (d)\)   \over  \psi (t)} \ge   1-  \ep -   \de_{m,n} |\ov\eta(d)|\)}\\
   &&   \ge \( 1-   \ep'(n)\) P\( \sup_{t'\in T'_{n,1}}{ \(  \eta (t')-   \eta (d)\)   \over  \psi (t)}\ge   1-  \frac\ep2-  \de_{m,n} |\ov\eta(d)| \). \nn   \eea
Now set,
\begin{equation} 
  P\(\de_{m,n}  | \ov \eta(d)  |  \le \ep/2\)= 1-\ep''(n). 
\end{equation}
This gives, 
\bea 
  &&P\( \sup_{t'\in T'_{n,1}}{ \(  \eta (t')-   \eta (d)\)   \over  \psi (t)}\ge   1-  \frac\ep2-  \de_{m,n} |\ov\eta(d)| \)\\&&\qquad\ge \nn P\( \sup_{t'\in T'_{n,1}}{ \(  \eta (t')-   \eta (d)\)   \over  \psi (t)}\ge   1-   \ep     \)\( 1-   \ep''(n)\).
\eea
Using this in   (\ref{4.37qq}) we get,
\bea\label{4.37qqk}
  \lefteqn{P\(  \sup_{t'\in T'_{n,1}}{ \(\ov \eta (t')- \ov \eta (d)\)   \over  \psi (t)} \ge   1-  \ep -   \de_{m,n} |\ov\eta(d)|\)}\\
   &&   \ge \( 1-   (\ep'(n)+\ep''(n)\) P\( \sup_{t'\in T'_{n, 1}}{ \(  \eta (t')-   \eta (d)\)   \over  \psi (t)}\ge   1-   \ep   \), \nn   \eea 
   Using this and (\ref{4.30nn}) we have, 
\bea 
 \lefteqn{P\(  \sup_{t'\in T'_{n,1}}{ \sum_{i=1}^{k}( \ov\eta^{2}_{i }(t')-\ov\eta^{2}_{i }(d)) \over   2 \psi(t) \| \vec{\ov \eta} (d)\|_2 }\ge  1-  \ep/2 \)}\\
  &&\quad \ge   \( 1-   (\ep'(n)+\ep''(n)\) P\( \sup_{t'\in T'_{n,1}}{ \(  \eta (t')-   \eta (d)\)   \over  \psi (t)}\ge   1-   \ep   \)\nn -\frac{ \ep(n)}{1-\ep(n)}.
\eea

  It follows from this that   for any $\ep'>0$, for all $n$ sufficiently large, 
    \bea \label{4.33aa}
  &&P\(  \sup_{t'\in T'_{n,1}}{ \sum_{i=1}^{k}( \ov\eta^{2}_{i }(t')-\ov\eta^{2}_{i }(d)) \over   2 \psi(t) \| \vec{\ov \eta} (0)\|_2 } \ge  1-  \ep'\)\\
  &&\qquad\ge (1-\ep')P\(  \sup_{t'\in T'_{n,1}}{ \(  \eta (t')-   \eta (d)\)   \over  \psi (t)} \ge   1-  \ep'  \)-\ep'.\nn
\eea
Using  this and  (\ref{3.5mm}) we   get what we need, that for any $\ep''>0$,  for all $n$ sufficiently large, 
 \begin{equation} \label{4.34}
  P\(  \sup_{t'\in T'_{n,1}}{ \sum_{i=1}^{k}( \ov\eta^{2}_{i }(t')-\ov\eta^{2}_{i }(d)) \over   2 \psi(t)} \ge \| \vec{\ov \eta} (d)\|_2    ( 1-  \ep'' )\)\ge 1-\ep ''.
\end{equation}

  We can now apply   Lemma \ref{theo-4.1}. Note that in the notation of  Lemma \ref{theo-4.1},
\begin{equation} 
  { \sum_{i=1}^{k}( \ov\eta^{2}_{i }(t')-\ov\eta^{2}_{i }(d)) \over   2   }=  Y_{k/2, (m,n)}( t')- Y_{k/2, (m,n)}(d), \end{equation}
      and
  \begin{equation}
\| \vec{\ov \eta} (d)\|_2 =\(2\sum_{i=1}^{k}\frac{\ov \eta^{2}_{i}  (d)}{2}\)^{1/2}=\sqrt2 Y_{k/2, (m,n)}^{1/2}(d).\label{2.20mmy}
\end{equation} 
  Therefore, 
by (\ref{4.34}) we have, 
   \be     P\(  \sup_{t'\in T'_{n, 1}}{ Y_{k/2, (m,n)}(t')- Y_{k/2, (m,n)}(d) \over  \psi(t)} \label{6.16qq} \ge  \sqrt 2 Y^{1/2}_{k/2, (m,n)}(d)\(1- \ep\)    \) \ge    1-\ep,
   %\label{4.39mm}
     \ee 
  and  finally,
  using  Lemma \ref{theo-4.1}   with $\wh X_{k/2,(m,n)} $ replaced by $  X_{k/2,(m,n)} $  we see that for all $n$ sufficiently large,  
 \bea
   &&         P\(  \sup_{t'\in T'_{n,1}}{ X_{k/2, (m,n)}(t')- X_{k/2, (m,n)}(d) \over  \psi(t)}   \ge  \sqrt 2 X^{1/2}_{k/2, (m,n)}(d)\(1- \ep''\)  \)\nn\\
   &&\hspace{3 in}  \ge    \nu_{m,n}^{- k/2} (1-\ep'' ). \label{2.2app}
     \eea

Now suppose that condition $(i)$ holds. Then it follows from Theorem  \ref{lem-3.3} that $\lim_{n\to\ff}\nu_{n,m} =1$.
Using   this in      (\ref{2.2app}) we have that for all    $\de_j>0$, for all $\ep_{j}>0$ sufficiently small,   
   \be         P\(  \sup_{t\in (0,\de_{j}]}{ X_{k/2 }(d+t)- X_{k/2 }(d) \over  \psi(t)}   \ge  \sqrt 2 X^{1/2}_{k/2 }(d)\(1- \ep_j\)  \)  \ge  (1-\ep_j ).   
     \ee 
 Choose  subsequences   $\{\de_{j_{i}},\ep_{j_{i}}\}$ of $\{\de_j,\ep_j\}$ such that $\lim_{i\to\ff}\de_{j_{i}}=0 $ and $\sum_{i=1}^{\ff}\ep_{j_{i}}<\ff$. It then follows from the Borel-Cantelli Lemma    that  (\ref{1.34sz}) holds. 
 
 Using Theorem \ref{lem-3.3q} the same argument gives (\ref{1.34sz}) when condition $(ii)$ holds.  
 
\medskip The last statement in this theorem 
 has the same proof as the proof of (\ref{1.34sz}) if we replace $T_n$ in (\ref{4.62ms}) by $\{-t_0,-t_1,-t_2,\ldots, -t_{m(n)}\}$.
  
 \qed

    \begin{remark}  {\rm   As we point out in   (\ref{1.29nn}), when the L\'evy processes in Theorem \ref{theo-1.2zz} have a Gaussian component,  $(\si ^{  \bb  })^2(x)\sim {|x|}/{C}$,  as $x\to 0$, for some constant $C$.  Consequently,  by (\ref{3.7q}),  condition  $(i)$ in Theorem \ref{theo-1.2z} does not hold. Nevertheless we can still use Theorem \ref{theo-1.2z} when condition $(ii)$ holds. 
 We can give insight into why   we require the condition 
$f'(d)=0$  in  Theorem  \ref{theo-1.2z}, $(ii)$. 

Consider a permanental vector with kernel,  
\begin{equation}\label{3.130q}
 G(m, n)_{j,k}+f(t_k) , \qquad j,k=0, 1,\ldots,m.
\end{equation}
For this vector,
\begin{equation} 
\nu_{n,m} =1+ {\bf  r}G(m, n){\bf  v}^{T}- {\bf  h}G(m, n){\bf  h}^{T}; \end{equation}
 see (\ref{15.17}).   In  our   proofs of lower bounds for LILs we require that $\lim_{n\to \ff}  \nu_{n,m} \newline  =1$, so we can use Theorem \ref{cor-7.1jq}.  
 Consider (\ref{u0}) when $\vf (x)=x$ so that,  
\begin{equation} 
 u^{(0)} (x,y)=x\wedge y ,\qquad x,y\in ( 0,\ff),  
\end{equation}
 which is the potential of Brownian motion killed the first time it hits $0$, (so that $u(t_j,t_k)= G(m, n)_{j,k}=t_j\wedge t_k$).
  The next two lemmas show that for this potential density  $\lim_{n\to \ff}  \nu_{n,m}\neq 1$ unless     take  $f'(d)=0$.  

\begin{lemma} \label{lem-3.1} Let $ G(m, n)$ be an  $(m+1)\times (m+1)$ matrix with elements, 
\begin{equation}\label{3.130}
 G(m, n)_{j,k}=  \( t_{j}\wedge t_{k}\) , \qquad j,k=0, 1,\ldots,m,
\end{equation}
in which  $t_{j}$ is a strictly increasing sequence. Then
\be G(m, n)^{-1}=  \left (
\begin{array}{ cccccc cc}  
{ a_{0}}+ { a_{1}}&- { a_{1}}&0&\dots &0 &0  \\
- { a_{1}}& { a_{1}}+ { a_{2}}&- { a_{2}}&\dots &0&0  \\
\vdots&\vdots&\vdots&\ddots&\vdots &\vdots  \\
0&0&0&\dots & { a_{m-1}}+ { a_{m}}&- { a_{m}}   \\
0&0&0&\dots &- { a_{m}}& { a_{m}} \end{array}\right ),\label{harr2.211w}
  \ee
  where 
    \begin{equation}
a_{0}=\frac{1}{t_{0}},\quad\mbox{and}\quad a_{j}=\frac{1}{t_{j}-t_{j-1}},\quad j\ge 1.\label{3.8bb}
\end{equation}
  \el

\Proof    It is easy to verify that this is the inverse of $G( m,m)$. See \cite[Lemma 2.4]{MRasym}. \qed

\begin{lemma}  For $G(m,n)$ in (\ref{3.130q}) with $g\equiv 1$, 
 \begin{equation} \label{slick}
 \nu_{n,m}=1+\frac{f(d+t_1)-f(d)}{2t_1}.\end{equation}
In particular,
\begin{equation} \label{slicka}
  \lim_{n\to \ff}  \nu_{n,m}=1+\frac{f'(d)}{2}.
\end{equation}
\end{lemma}

Obviously, to have  $\lim_{n\to \ff}  \nu_{n,m}=1$ we must take $f'(d)=0$.

\medskip\Proof \begin{equation} 
 {\bf  h}G(m, n){\bf  h}^{T}=\sum_{j,k=0}^{m}(r_jv_j)^{1/2}G(m, n) _{j,k}(r_kv_k)^{1/2}=r_0v_0G(m, n) _{0,0}=v_0,
\end{equation}
since $r_0=Ct'_0$ and $G(m, n)_{0,0}=t'_0/C$. In addition,
\bea 
 {\bf  r}G(m, n){\bf  v}^{T}&=&\sum_{j,k=0}^{m}r_jG(m, n)_{j,k} v_k =r_0 \sum_{ k=0}^{m}G(m, n)_{0,k} v_k\\
 &=& r_0G(m, n)_{0,0}\(v_0+\sum_{ k=1}^{m}  v_k\)=v_0+\sum_{ k=1}^{m}  v_k
 \nn ,
\eea
since $G(m, n)_{0,k}=G(m, n)_{0,0}$.
Consequently,
\begin{equation} 
  {\bf  r}G(m, n){\bf  v}^{T}- {\bf  h}G(m, n){\bf  h}^{T}=\sum_{ k=1}^{m}  v_k.
\end{equation}
It follows from this that,\begin{equation} 
 \sum_{ k=1}^{m}  v_k=\sum_{ k=1}^{m} \sum_{ l=0}^{m} G(m, n)^{k,l}f(t'_l)= \sum_{ l=0}^{m}\(\sum_{ k=1}^{m} G(m, n)^{k,l}\)f(t'_l).
\end{equation}
For each $l$, $S_{l}:=\sum_{ k=1}^{m} G(m, n)^{k,l}$ is the sum over all elements in the $l$--th column of $G(m, n)^{ -1}$,   except  the first element,  $G(m, n)^{0,l}$ . Therefore, by Lemma \ref{lem-3.1}, we have 
 $S_{0}=- a_{1}/2$, $S_{1}=-a_{1}/2$ and $S_{l}=0$ for all $l\geq 2$,  where ,\be a_{1}=\frac{1}{t'_{1}-t'_{0}}=\frac{1}{t_{1}}.\ee
This shows that 
\begin{equation} 
 \sum_{ k=1}^{m}  v_k=-\frac{ a_{1}f(t'_0)}2+\frac{ a_{1}f(t'_1)}2=\frac{f(t'_1)-  f(t'_0)}{2(t'_{1}-t'_{0})},
 \ee
  which gives (\ref{slick}).\qed
}\end{remark}

\section{Proof of Theorems \ref{theo-1.2zz}--\ref{theo-1.5}}\label{sec-4}

\subsection{Proof of lower bounds in Theorems \ref{theo-1.2zz}--\ref{theo-1.2aa}}

The increments variance $(\si^\bb)^{2}(x)$, $\bb\ge 0$ plays a critical role in all these theorems. The next lemma gives some properties of this function.

 \begin{lemma} \label{lem-3.5mm}\begin{itemize} 
  \item[(i)]For L\'evy processes, the L\'evy exponent,
     \be \label{3.79nn}
  \psi(\la)=C\la^2+\psi_1(\la),
  \ee
  for some $C\ge0$, where  $\psi_1(\la)=o(\la^2)$ as $\la\to\ff$;
\item[(ii)] For  $(\si^\bb)^{2}(x)$   in (\ref{lv.21}),  when $C>0$,
 \begin{equation} \label{3.82nn} 
  \frac{(\si^\bb)^2(x)}{|x|}\le\frac{1}{C},\qquad  \forall\, \bb\ge 0 ,\ee
and, \begin{equation} \label{3.94nn}
  \limsup_{x\to0}\frac{|x|}{(\si^\bb)^{2}(x)}=C,\qquad  \forall\, \bb\ge 0 ; 
\end{equation}
\item[(iii)]   When $C=0$, in which case  in  the L\'evy process does not have a Gaussian component,  
\begin{equation} \label{3.94nnq}
  \lim_{x\to0}\frac{|x|}{(\si^\bb)^{2}(x)}=0,\qquad  \forall\, \bb\ge 0 . 
\end{equation}
\end{itemize}
\end{lemma}

\Proof The statement in $(i)$  is well known. (See e.g., \cite[page 138]{book}.)  

\medskip The statement in (\ref{3.94nn})  is given in (\ref{1.29nn}).  The statement in (\ref{3.82nn})  is trivial since   
    \begin{equation} \label{3.103}
 (\si^\bb)^2 (x)\le  \frac{4}{ \pi}  \int_{0}^{\ff}\frac{ \sin^2\la x/2}{ C\la^2  } \,d\la,
\end{equation}
and  by  \cite[3.741.3]{GR},
\begin{equation} \label{3.95mm} 
  \frac{4}{ \pi}  \int_{0}^{\ff}\frac{ \sin^2\la x/2}{ C\la^2  } \,d\la=\frac{4|x|}{ C\pi}  \int_{0}^{\ff}\frac{ \sin^2s/2}{ s^2  } \,ds= \frac{|x|}{C} .
\end{equation}

We now obtain (\ref{3.94nnq}).
Consider  
\be  \label{3.96mm}
(\si^\bb)^2(x)   = \frac{4}{ \pi}  \int_{0}^{\ff}\frac{ \sin^2\la x/2}{\bb+\psi_1(\la)   } \,d\la\\ \ee 
Let,  
 \begin{equation} 
  M_\ep =\sup \{\la: \bb+\psi_1(\la) \ge  \ep  \la^2\} ,
\end{equation} 
 and  write (\ref{3.96mm}) as,
 \bea &&     \frac{4}{ \pi}  \int_{0}^{M_{\ep}}\frac{ \sin^2\la x/2}{\bb+\psi_1(\la)   } \,d\la +\frac{4}{ \pi}  \int_{ M_{\ep}}^{\ff}\frac{ \sin^2\la x/2}{\bb+\psi_1(\la)    } \,d\la\qquad\\
 &&\qquad \nn :=  (\si_1^\bb)^2(x) +(\si_2^\bb)^2(x).
 \eea
 We have,  
\be 
(\si_2^\bb)^2(x)\ge  \frac{4}{\ep \pi}  \int_{ M_{\ep}}^{\ff}\frac{ \sin^2\la x/2}{   \la^2  } \,d\la= \frac{4|x|}{\ep \pi}  \int_{ xM_{\ep}}^{\ff} \frac{ \sin^2s /2}{ s^2  } \,ds.
\ee
Consequently,  
  \begin{equation} 
 \frac{|x|}{(\si_2^\bb)^2 (x)}\le \ep\(\frac{4 }{  \pi}  \int_{ xM_{\ep}}^{\ff} \frac{ \sin^2s /2}{ s^2  } \,ds\)^{-1}.
\end{equation}
Using this,   (\ref{3.95mm}) and the fact that $(\si_1^\bb)^2 (x)>0$ we see that,
 \begin{equation} 
 \limsup_{x\to 0}\frac{|x|}{(\si ^\bb)^2 (x)}\le   \limsup_{x\to 0}\frac{|x|}{(\si_2^\bb)^2 (x)}\le \ep.
\end{equation}
   Since this holds for all $\ep>0$, we get (\ref{3.94nnq}).   \qed
 
 \begin{lemma} \label{lem-2.9}   When $\vf(x)$
 satisfies  the first line of  (\ref{3.7q}) it satisfies (\ref{2.149}), i.e. 
 \begin{equation} \label{2.149ww}
 \limsup_{x\to 0} \frac{|x|}{\vf(x)}=O(1) .
\end{equation}
\end{lemma}
  
 \Proof This follows immediately from  (\ref{3.94nn}) and  (\ref{3.94nnq}). \qed
  
   We  show that the potentials considered in Theorems \ref{theo-1.2zz}--\ref{theo-1.2aa} satisfy (\ref{3.7q}) and show   what   $\vf_d(x)$ is. In the case of Theorem \ref{theo-1.2zz} this is simple because, 
 \begin{equation} \label{3.11aw}
    u^\bb (x,x)-u^\bb (x,y)= \frac{(\si^\bb)^2(x-y)}{2} ,\qquad\forall x,y\in  \De_d ,\,\forall d\ge 0 .
\end{equation}
By hypothesis, $(\si^\bb)^2(x)$ is a regularly varying function at 0 with index $0<r\le 1$.    It follows from  \cite[Lemma 7.2.4]{book},  in particular (7.125),   that  a regularly varying function with index not equal to 0, is asymptotic  to a function with a positive derivative in some neighborhood of $(0,\de]$,  for some $\de>0$. Label this function  $(\rho^\bb)^2(x)$,
 For each $\bb>0$ we take   
\be\label{3.69nn}
 \vf(x)=(\rho^\bb)^2(x).
\ee
 Therefore in the cases we consider, in some neighborhood of 0,    $\vf$ is a continuous strictly increasing  regularly varying function with index $0<r\le 1$.  
(Actually as we point out in (\ref{1.26mmw2}),
under mild regularity on $\psi$, all the   functions $(\si^\bb)^2(x)$, $\bb\ge 0$, are asymptotically equivalent at 0.)

  \medskip In Theorem \ref{theo-1.2zzk} we consider potential densities,
 \begin{equation}
 v^\bb(x,y)=u^\bb(x,y)-\frac{u^\bb(x )u^\bb( y)}{u^\bb(0,0)},\qquad x,y\in\De_d(\de),\label{80.4aa}
 \end{equation}
 for some $\de>0$, which may depend on $x$ and $y$ and where   $\{u^\bb(x,y)  \}$ is translation invariant.  The next lemma shows the relationship between  $v^\bb(x,y)$ and $u^\bb(x,y)$ when $\psi\(\la\)$  does not have a Gaussian component. This is the case in Part I, since  (\ref{1.38nnb}) implies that  \begin{equation} 
  \limsup_{x\to0}\frac{|x|}{(\si^\bb)^{2}(x)}=0,\qquad  \forall\, \bb\ge 0 ; 
\end{equation}
which contradicts (\ref{3.94nn}).
 
\bl \label{lem-80} Let $v^\bb(x,y)$, $u^\bb(x,y)$ and $(\si^\bb)^2(x)$ be as in Theorem \ref{theo-1.2zzk},  in which $\bb>0$,   and assume that $\psi\(\la\)$  does not have a Gaussian component. Then for  any  $d\neq 0$,  
\begin{equation} \label{80.7}
\frac{v^\bb(x,x)-v^\bb(x,y)}{(\si^\bb)^2(x-y)}  \sim \frac{1}{2} ,   \qquad\mbox{as } \,\,x,y\to d, \,\,x\neq y,
\end{equation}
and consequently,
\begin{equation} \label{80.7t}
\frac{v^\bb(x,x)+v^\bb(y,y)-2v^\bb(x,y)}{(\si^\bb)^2(x-y)}  \sim 1,   \qquad\mbox{as } \,\,x,y\to d, \,\,x\neq y,
\end{equation}
\el  \Proof  \begin{eqnarray}
  v^\bb(x,x)-v^\bb(x,y)&=&\frac{  (\si^\bb)^2 (x-y)}{2 }-
 \frac{u^\bb(x )\(u^\bb( x)-u^\bb( y)\)}{u^\bb(0 )} 
 \label{80.8qx}\\
 &=&\frac{  (\si^\bb)^2(x-y)}{2 }-
 \frac{u^\bb(x )\( (\si^\bb)^2(y)- (\si^\bb)^2(x)\)}{2u^\bb(0 )}.\nn  
 \end{eqnarray}

  By hypothesis,   $ (\si^\bb)^2\(x\)\in C^{ 1}\(T\)$. Therefore there exists a constant $C$ such that,
 \begin{equation} 
| (\si^\bb)^2(y)-(\si^\bb)^2(x)|\leq C|x-y|, \hspace{.2in} \forall x,y\in [d,d+\de].  \label{80.9}
 \end{equation}
  Using (\ref{3.94nnq}) we see that, 
\begin{equation} \label{80.9u}
 \lim_{|x-y|\to 0} \frac{| x-y|}{(\si^\bb)^2(x-y)}=0.
\end{equation} 
This gives  (\ref{80.7}).\qed

When $ \vf(x)=(\rho^\bb)^2(x)$ as in (\ref{3.69nn}),   (\ref{3.7q}) holds.

\medskip In  Theorem \ref{theo-1.2zzk}, II  we add the hypothesis that the L\'evy exponent $\psi(\la)$  of $\,Y'$ is a stable mixture. In this case (\ref{3.7q})   follows from the next   lemma:

  \begin{lemma} \label{lem-3.6m}   When  the L\'evy exponent $\psi(\la)$  of $\,Y'$ is a stable mixture, for all $\bb>0$     and all $d> 0$  there exists a $\de>0$ such that,
   \begin{equation} \label{3.11aa}
    v^\bb(x,x)-v^\bb (x,y)\asymp  (\si^\bb)^2(x-y) ,\qquad\forall x,y\in \De_{d}(\de) ,
\end{equation}
where $\de $ may depend on   $ d$. 
\end{lemma}

 \Proof We have, \bea 
   v^\bb(x,x)-v^\bb(x,y)&=&\frac{  (\si^\bb)^2 (x-y)}{2 }-
 \frac{u^\bb(x )\(u^\bb( x)-u^\bb( y)\)}{u^\bb(0 )} \label{80U.8qx}\\
 &=&\frac{  (\si^\bb)^2 (x-y)}{2 }-
 \frac{u^\bb(x )\((\si^\bb)^2 (y)-(\si^\bb)^2 (x)\)}{2u^\bb(0 )} . \nn
\eea
  so that,
\be  
   \frac{v^\bb(x,x)-v^\bb(x,y)}{(\si^\bb)^2 (x-y)}=\frac{ 1 }{2 }-
  \frac{u^\bb(x )\((\si^\bb)^2 (y)-(\si^\bb)^2 (x)\)}{2u^\bb(0 )(\si^\bb)^2 (x-y)} . \label{80U.8qs} 
 \ee  
Note that for all $x\in R^1-\{0 \}$ there exists an $\ep>0$ such that,
\begin{equation} 
  \frac{u^\bb(x ) }{ u^\bb(0 ) }\le 1-\ep,
\end{equation}
this follows from the definition of $u^\bb$, see (\ref{1.21nn}), and the fact that $\psi(\la)$ is a continuous function. Therefore to get the lower bound in (\ref{3.11aa}) it suffices to show that   for any $d> 0$ we can find  a $\de>0$,  which may depend on $d$, with
\begin{equation} \label{3.111mm}
    \frac{| (\si^\bb)^2 (y)-(\si^\bb)^2 (x) |}{ (\si^\bb)^2 (x-y)}\le 1+\frac{\ep}{2},\qquad\forall x,y\in \De_d(\de).
\end{equation}
 Since $ (\si^{\bb})^2(x)\in C^{1}\(T\) $ there exists at least one  $d^*\in \De_d(\de)$ for which, 
\begin{equation} 
  |((\si^\bb)^2 (d^*))'|=\max_{0\le u \le |x-y|}\big|((\si^\bb)^2 (u))'\big|,\end{equation}
which implies that,
\begin{equation} \label{3.111}
  \big|(\si^\bb)^2 (y)-(\si^\bb)^2 (x)\big|\le  |((\si^\bb)^2 (d^*))'||x-y|.
\end{equation}
  If $\psi\(\la\)$ in (\ref{3.79nn}) is such that $C=0$, it follows from (\ref{3.94nnq}) that 
\begin{equation} \label{3.14nn}
    \frac{ |x-y| }{ (\si^\bb)^2 (x-y)}\to 0,  \qquad\text{as $|x-y|\to0$ },
\end{equation}
which   gives (\ref{3.111mm}).

If $\psi\(\la\)$ in (\ref{3.79nn}) is such that $C>0$, it follows from (\ref{4.58aqe}) and (\ref{3.82nn}) that,
\be \label{3.111ax}
  \big|(\si^\bb)^2 (y)-(\si^\bb)^2 (x)\big|\le  \frac{ (\si^\bb)^2 (d^*) }{d^*}|x -y |\\ 
   \le   \frac{|x-y|   }{C }. 
\ee 
  Set $x=d+x'$ and $y=d+y'$.  Then, using  (\ref{3.94nn}) we see that, 
  \begin{equation} \label{3.111fgh}
 \limsup_{|x'-y'|\to 0} \frac{\big|(\si^\bb)^2 (y)-(\si^\bb)^2 (x)\big|}{(\si^\bb)^2 (x-y)}= 1,
\end{equation}
which gives (\ref{3.111mm}).\qed

%  We have,
%\begin{equation} \label{3.11aa}
%    v^\bb(x,x)-v^\bb (x,y)\asymp   \vf^\bb  (x-y) ,\qquad\forall x,y\in [d,d+\de],\,\forall d> 0,
%\end{equation}
%for $\vf^\bb$ as given in (\ref{3.69nn}).

  In  Theorem \ref{theo-1.2aa} the Gaussian process   $\eta=\{\eta(x),x\in R^1 \}$  defined in (\ref{lv.22})--(\ref{1.44mm})  has   covariance,   
\begin{equation} \label{u0q}
  u^{(0)} (x,y)=\phi(x)+\phi(y)-\phi(x-y) ,\qquad x,y\in R^1.  
\end{equation} 
In this case,
\bea \label{sig0q}
  (\si^{(0)})^2(x,y)&=&E\(\eta(x)-\eta(y)\)^2 \\
  &=&\nn u^{(0)} (x,x)+u^{(0)} (y,y) -2u^{(0)} (x,y)=2\phi(x-y).
\eea
Since $ u^{(0)}(0,0)=0$  it follows from (\ref{poscond.8}) that we restrict $d$ to $R^1-\{0 \}$. Then,  
\be u^{(0)}(d+x,d+x)-  u^{(0)}(d+x,d+y)=\phi( d+x) -\phi( d+y)  + \phi( x-y).
\ee
If  $ \phi (x)\in C^{1}(\De_d)$   for some   $d\neq 0$, then
\be 
\phi( d+x)-\phi( d+y) \asymp\phi'(d) |x-y|   ,\qquad\text{as}\quad |x-y|\to 0,
\ee 
and
\begin{equation} 
   u^{(0)}(d+x,d+x)-  u^{0}(d+x,d+y)\asymp  \phi'(d)|x-y|+\phi( x-y).
\end{equation}

In  Theorem \ref{theo-1.2aa}, I,  it follows from  (\ref{1.38nnqq}) that,  
\begin{equation} 
  \limsup_{|x-y|\to 0}\frac{|x-y|}{\phi( x-y)}=0
\end{equation}
and we have,
\bea \label{2.31}
   u^{(0)}( x, x)-  u^{(0)}( x, y)&\asymp &  \phi( x-y),\,\, \text{and,}\\  (\si^{(0)})^2(x,y)&= &  2\phi( x-y),\qquad \forall x,y\in   \De_d(\de),\nn
\eea
 for some $\de>0$  sufficiently small.

  By hypothesis, $(\si^0)^2(x)=2\phi( x)$ is a regularly varying function at 0 with index $0<r\le 1$. Such a function is asymptotic to a strictly increasing  function at zero, say $(\rho^0)^2(x)$. (See e,g,  \cite[Lemma 7.2.4]{book}.)  Therefore, it follows from (\ref{2.31}) that   (\ref{3.7q}) holds with
\be\label{3.69nnaa}
 \vf (x)=(\rho^0)^2(x).
\ee
 It   follows from (\ref{u0}) that,
\begin{equation} \label{u02x}
u^{(0)} (x,x)-u^{(0)} (x,y)=\phi(x)-\phi(y)+\phi(x-y) ,\qquad x,y\in R^1-\{0\}.  
\end{equation}
 In  Theorem \ref{theo-1.2aa}, II  we add the condition that   $ (\si^{(0)})^2(x) =2\phi( x)$ is increasing in some neighborhood of   $d$. Therefore  for $x,y$
in this neighborhood,
\begin{equation} \label{u02}
u^{(0)} (x,x)-u^{(0)} (x,y)\ge \phi(x-y) ,\qquad y\le x.  
\end{equation}
Since   $2\phi(x-y)=(\si^{(0)})^2(x,y) $, defining  $\vf(x)$ as in (\ref{3.69nnaa}), we see that   (\ref{3.7q}) is satisfied.
This completes the proof that the potentials considered in Theorems \ref{theo-1.2zz}--\ref{theo-1.2aa} satisfy (\ref{3.7q}). 
  
 \medskip  We now show that for the Gaussian processes with covariances given by the potentials in Theorems \ref{theo-1.2zz}--\ref{theo-1.2aa},  (\ref{3.5mm}) holds. The inequalities in the next two lemmas are critical because they allow us to use Slepian's Lemma.

   \begin{lemma} \label{lem-4.1mm}   Let  $\{\eta (x),x\in R^{1}\}$ be a stationary Gaussian process with covariance    $ u^\bb (x-y)=u^\bb(x,y)$, $ \bb> 0$ as in Theorem \ref{theo-1.2zz}.  Set  
\begin{equation} \label{sidef.2}
 (\si^\bb)^{2}\(x\)=  E \(\eta(x)-\eta(0)\)^2   =2\(u^\bb(0)-u^\bb(x)\). 
 \end{equation} 
 Then for     $0< x\le y\le \De'$,   
  \begin{equation} \label{1.21wq}
\frac{  E\(\(\eta(d+x)-\eta(d)\)\(\eta(d+y)-\eta(d)\)\) }{ \si^\bb  (x)\si^\bb  ( y)}\le  C_{\De'}\frac{ \si^\bb  (x)}{ \si^\bb   (y)},  \end{equation} 
for some   constant $C_{\De'}$ depending on $\De'$. 
\end{lemma}

\Proof 
  Since $\eta$ is a stationary Gaussian process it suffices to show  (\ref{1.21wq})  for    $d=0$, i.e., that,  
  \begin{equation} \label{3.71mm}
\frac{  E\(\(\eta(x)-\eta(0)\)\(\eta(y)-\eta(0)\)\) }{\si^\bb   (x)\si^\bb   ( y)}\le C_{\De'}\frac{ \si^\bb (x)}{ \si^\bb (y) }.
\end{equation}  
Note that,
\be  \label{3.72mm}
  E\(\(\eta(x)-\eta(0)\)\(\eta(y)-\eta(0)\)\) = \frac12\((\si^\bb)^2(x)+ (\si^\bb)^2(y)-(\si^\bb)^2(y-x)\). 
\ee 
It follows from  \cite[Lemma 7.4.2]{book} and    (\ref{poscond.8})   with $T=R^{1}$ that when $\bb>0$,  
 \bea \label{3.73mm}
(\si^\bb)^2(y)- (\si^\bb)^2(y-x) & \le &\frac{2u^\bb(0)}{u^\bb(x)+u^{\bb}(y)} (\si^\bb)^2(x).\\ 
& \le &C'_{\De'}\, (\si^\bb)^2(x),\nn
\eea
for some   constant $C'_ {\De'}$ depending on $\De'$.
 Using   (\ref{3.72mm}) and (\ref{3.73mm}) we get  (\ref{3.71mm}).   

\qed

     \begin{lemma} \label{lem-4.3}       Let 
  $\{\eta (x),x\in R^{1}\}$ be a Gaussian process with covariance \newline  $ \{u^{(0)}(x,y),x,y\in R^1\}$ as in Theorem \ref{theo-1.2aa}. Let
  \begin{equation} \label{1.14nn}
  (\si^{(0)})^2(x,y)= E (\eta( x)-\eta( y))^2,
\end{equation}
 and note that $(\si^{(0)})^2(x,y)=(\si^{0})^2(x-y)$.
Then for all   $x,y,d \in R^1-\{0 \}$,
  \begin{equation} \label{1.21wz}
\frac{  E\(\(\eta(d+x)-\eta(d)\)\(\eta(d+y)-\eta(d)\)\) }{\si^{0}  (d+x,d ) \si ^{0} (d+y,d  )}\le \frac{ \si^{0}(x )}{ \si^{0}(y )}. \end{equation}  
\end{lemma}

\Proof Since $\eta $ is a Gaussian process with stationary increments and $\eta(0)=0$,  the   left--hand side of (\ref{1.21wz}) is equal to,
  \begin{equation} \label{1.21wza}
\frac{  E\( \eta(x) \eta(y) \) }{ \si^{0}(x )  \si^{0}(y )}= \frac{ u^{(0)} (x,y)}{\si^{0}(x )   \si^{0}(y )} .
\end{equation} 
By hypothesis  $u^{(0)}=\{ u^{(0)}(x,y),x,y\in T\}$ is the potential density of a  symmetric  Borel right  process with state space   $R^1-\{0 \}$.  Therefore, by (\ref{smp.1}) for all   $x,y\in T$,
\begin{equation} 
   \frac{ u^{(0)} (x,y)}{\si^{0}(x) \si^{0}    (y  )}\le \frac{ u^{(0)} (x,x)}{\si^{0}(x )   \si^{0}(y )} =  \frac{ \si^{0}(x )}{ \si^{0}(y )},
\end{equation}
  where we use    (\ref{u0}) and (\ref{1.44mm}).     \qed

  We now show that condition (\ref{3.5mm}) in   Theorem \ref{theo-1.2z} is satisfied by the Gaussian processes considered in Lemma \ref{lem-4.1mm}.     We write it with the following change of notation.
   For all $\ep>0$,   
 \be  P\(\sup_{1\le j\le m(n)} \frac{\eta(d+t_j)-\eta(d)}{\si^\bb (t_j) (2\log \log  1/t_j)^{1/2}}\ge  1-\ep\) \ge      1-\ep   \label{4.51}  ,
\ee
for all $n$ sufficiently large.

 \medskip It follows from Lemma \ref{lem-4.1mm} that for $j<k$, which implies that $t_j<t_k$,  
 \begin{equation} \label{1.21w}
\frac{  E\(\(\eta (d+t_j)-\eta (d)\)\(\eta (d+t_k)-\eta (d)\)\) }{\si^\bb  ( t_j)\si^\bb   ( t_k)}\le C\frac{ \si^\bb  (t_j)}{ \si^\bb  (t_k)} . \end{equation}  
Since $(\si^\bb)^2(x)$ is regularly varying at $0$ with index  greater than 0,
\begin{equation} \label{4.44}
  \frac{ \si^\bb (t_j)}{ \si^\bb   (t_k)}\le \sup_{1<q\le  t_{m(n)} } \frac{ \si^\bb   (\th^{q})}{ \si^\bb   (\th^{q-1}) }\le \ep(\th).
\end{equation} 
Therefore,   the left--hand side of  
(\ref{1.21w}) is bounded by $C\ep(\th)$.
This term goes to 0 as $\th\to 0$ and $n\to\ff.$  

 This allows us to use Slepian's Lemma. Let,  
\begin{equation} 
  \ze_j=\frac{  \eta (d+t_j)-\eta (d) }{ \si^\bb  ( t_j) },\qquad j\in\{1,\ldots,m(n)\}.
  \end{equation}
  It follows from (\ref{1.21w}) and   (\ref{4.44}) that for any $0<\th<1$ and $\ep'(\th)$, 
 \begin{equation} \label{3.121z}
  E\ze_j\ze_k\le   \ep'(\th),\qquad \forall\,j,k\in \{1,\ldots,t_{m(n)} \},\end{equation} 
  for all $n$ sufficiently large.
  Let   $\{ \xi_j\}_{j=0}^{m(n)}$ be independent standard normal random variables and consider the sequence, 
\begin{equation} 
  \wt\xi_j=(1- \ep'(\th))^{1/2} \xi_j+\ep'(\th)^{1/2} \xi_{0},\qquad j \in \{1,\ldots,t_{m(n)}\}.
  \end{equation}
Note that,
\begin{equation} \label{4.62}
  E\ze_j^2= E\wt\xi_j^2,\qquad \text{and}\qquad  E\ze_j\ze_k\le  E \wt\xi_j \wt\xi_k,\qquad  \forall\,j,k\in \{1,\ldots,t_{m(n)} \}. 
\end{equation}
Set,
\begin{equation} 
  b_j=\(2\log\log 1/t_j\)^{1/2}. \label{4.62mr}
\end{equation}  
Using (\ref{4.62}), it follows from Slepian's Lemma, (see e.g.  \cite[Lemma 5.5.1]{book}
), that,     
   \be 
  P\(\sup_{1\le j\le m(n)}\frac{\ze_{j}}{b_j}\ge  (1-\ep'(\th))\)  \ge    P\(\sup_{1\le j\le m(n)}  \frac{\wt\xi_{j}}{b_j}\ge   (1-\ep)  \)  
  \ee
 \bea &&\qquad \ge    P\(\sup_{1\le j\le m(n)} \frac{\(1- \ep'(\th)\)^{1/2}\xi_{j}}{b_j}\ge   (1-\ep/2), \ep'(\th)^{1/2} \xi_{0}<\ep/2  \).\nn
  \eea

Therefore, by taking $\th>0$ sufficiently small we have,    
\bea \label{3.124}
 && P\(\sup_{1\le j\le m(n)}\frac{\ze_{j}}{b_j}\ge  (1-\ep) \)\\
 &&\qquad \ge  (1-\ep)  P\(\sup_{1\le j\le m(n)} \frac{\xi_{j}}{b_j}\ge   (1-\ep/2)  \).\nn
\eea

  \begin{lemma} \label{3.148z} Let $\{\xi _j\}_{j=1}^{m(n)}$ be independent standard normal random variables, Then  for all $\ep>0$, for all $ n$ sufficiently large,  
 \begin{equation} \label{3.125}
 P\(\sup_{1\le j\le m(n)}  \frac{\xi_{j}}{b_j}\ge  (1-\ep) \)\ge  1-e^{-n^{\ep  }  }.\end{equation}
\end{lemma}

\Proof   Refer to (\ref{4.62mr}) and (\ref{4.62ms})
to see that,
\bea 
 b_j&=&\(2\log\log \(\frac{1}{  \th^{n+1-j}}\)\)^{1/2}\\
 &=&\(2\log (n+1-j)+ 2\log\log \(\frac{1}{  \th }\)\)^{1/2}\nn\\
 &:=&\(2\log (n+1-j)+   h(\th)\)^{1/2} \nn.
\eea 
Note that for all $\th$ sufficiently small,
\begin{equation} \label{3.77x}
 \inf_{1\le j\le m(n)}b_j\ge \(2\log (n+1-[n^q])\)^{1/2}. 
\end{equation}

Let $\wt b_j=(1-\ep)b_j$.  It follows from   \cite[(5.19)]{book} and the fact that $P( \xi_{j}\ge \wt  b_{j})=(1/2)P( |\xi_{j}|\ge \wt  b_{j})$ that for all $\th$ and $\ep$ sufficiently small,  for all $1\le j\le m(n)$,  for all $n$ sufficiently large, 
\bea \label{3.79nnk}
 P\( \xi_{j}\ge \wt  b_{j}\) &\ge &\frac{e^{-\wt b_{j}^{2}/2}}{2\sqrt{2\pi}\,\,\wt b_{j}}=     \frac{e^{-(1-\ep)^{2}    (\log (n+1-j)+   h(\th)/2)   }}{2\sqrt{2\pi}\,\,(1-\ep)\(2\log (n+1-j)+   h(\th)\)^{1/2}}\nn\\
  &\ge &\frac{e^{-(1-\ep)^{2}        h(\th)/2    }\,e^{-(1-\ep)^{2}     \log (n+1-j) }}{2\sqrt\pi(1-\ep)[(2\log (n+1-j)^{1/2}+h^{1/2}(\th)]}\\
   &\ge &\frac{e^{-(1-\ep)^{2}        h(\th)/2    }\,e^{-(1-\ep)^{2}     \log (n+1-j) }}{4\sqrt\pi(1-\ep) ( \log  n )^{1/2} }\nn\\
   &:=& \frac{ C(\th,\ep)   }{  ( \log  n )^{1/2} (n+1-j)^{(1-\ep)^2}}.\nn
   \eea
 We now note that,  
 \bea \label{3.75nn}
  P\(\sup_{1\le j\le m(n)}  \frac{\xi_{j}}{\wt b_{j}}\ge 1\) &=&1-\prod_{1\le j\le  m(n)} P\( \xi_{j}<\wt b_{j}\)\\
  &=& 1-\prod_{1\le j\le m(n)}  \(1-P\( \xi_{j}\ge \wt  b_{j}\)\)\nn. 
\eea
Let $x_j=P( \xi_{j}\ge \wt  b_{j})$. By the well known estimate, (see  e.g., \cite[(5.18)]{book}), 
\begin{equation} \label{3.77}
 \sup_{1\le j\le m(n)}x_j\le  \frac{1}{(n+1 -[n^{q}])^{1-\ep}}. 
\end{equation}
Since this implies that $\lim _{n\to\ff} \sup_{1\le j\le m(n)}x_j=0$
we see that, for all $\vep>0$, for all  $n$ sufficiently large  
\begin{equation} \label{3.81}
  \prod_{1\le j\le m(n)} (1-x_j)=\exp\(\sum_{j=1}^{m(n)}\log(1-x_j)\)\le \exp\(-(1-\vep)\sum_{j=1}^{m(n)} x_j \).
\end{equation}
Using (\ref{3.79nn}) we see that,  
\bea 
\sum_{j=1}^{ m(n)}x_j &\ge&  \frac{C (\th,\ep) }{(\log n)^{1/2}}    \sum_{j=1}^{ m(n)}\frac{1}{ (n+1-j)^{ (1-\ep )^{2}}   }\\
&= &\frac{C (\th,\ep) }{  ( \log n)^{1/2}}  \sum_{k=   [n]^q}^n\frac{1}{       k^{(1-\ep )^{2}}}\nn  \\
 &\ge&\nn  \frac{C (\th,\ep) }{2\ep(\log n)^{1/2}}  n^{ 2\ep-\ep^{2}},
\eea
 for all $n$ sufficiently large. Using this in (\ref{3.81}) and then in (\ref{3.75nn}) we get (\ref{3.125}).\qed

  Combining Lemma \ref{3.148z} and (\ref{3.124}) we get that for  all $\ep>0,$  
\be  P\(\sup_{1\le j\le m(n)}\frac{\ze_{j}}{(2\log\log 1/t_j)^{1/2}}\ge  (1-\ep)\) \ge (1-\ep)   \(1-e^{-n^{\ep }  }\),
\ee
 or equivalently,  
 \bea &&  P\(\sup_{1\le j\le m(n)} \frac{\eta(d+t_j)-\eta(d)}{\si^\bb (t_j)(2\log\log 1/t_j)^{1/2}}\ge  (1-\ep) \) \\
 &&\qquad\qquad \ge     (1-\ep)   \(1-e^{-n^{\ep  }  }\),\nn 
\eea
for all $n$ sufficiently large.   This gives    (\ref{4.51}) or , equivalently, (\ref{3.5mm}).

\medskip    It should be obvious that, using Lemma \ref{lem-4.3} exactly the same proof  shows that condition (\ref{3.5mm}) in   Theorem \ref{theo-1.2z} is satisfied by the Gaussian processes considered in Lemma \ref{lem-4.3},    and that it shows that 
   for all $\ep>0$,    
 \be  P\(\sup_{1\le j\le m(n)} \frac{\eta(d+t_j)-\eta(d)}{\si^0 (t_j)}\ge  (1-\ep)(2\log \log  1/t_j)^{1/2}\) \ge     (1-\ep)  \label{4.51q}  ,
\ee
for all $n$ sufficiently large.

\medskip

  Now let  $\{\wh\eta (x);x\in R^{1}\}$ be a Gaussian process with covariance   $ v^{\bb}(x,y)$, $\bb>0, $ as in Theorem \ref{theo-1.2zzk}. Note that,
  \begin{equation} 
  \wh\eta(x)=\eta(x)-\frac{u^\bb(x)}{u^\bb(0)}\eta(0),\qquad x\in R^1-\{0 \},
\end{equation}
  so that,
\bea \label{4.74}
   \wh \eta(d+t_j)-\wh\eta(d)&=&\eta(d+t_j)-\eta(d)-\frac{u^\bb(d+t_j)-u^\bb(d)}{u^\bb(0)}\eta(0)\\
   &=&\eta(d+t_j)-\eta(d)-\frac{(\si^\bb(t_j))^2}{2u^\bb(0)}\eta(0).\nn
\eea

\medskip
Set,  \be 
 (\wh \si^\bb)^{2}\( d+t_j,d\)=  E( \wh\eta(   d+t_j)- \wh \eta  (d))^2.
 \label{80.10rf}   
\ee
It follows from Lemma \ref{lem-80} that. \begin{equation} \label{80.7rr}
\frac{(\wh \si^\bb)^{2}\( d+t_j,d\)}{(\si^\bb)^{2}\( t_j\)} \sim 1, \quad\mbox{as } \,\,  t_j\to 0.
\end{equation}
Also, since obviously, 
 \begin{equation} 
  \frac{ \si^\bb(t_j)  }{( \log\log 1/t_{j})^{1/2}}\to 0,\qquad \forall \,1\le j\le m(n),\quad\text{as }\,n\to\ff,
\end{equation} 
 it follows from  (\ref{4.51})    that  
 for all $\ep>0$,   
 \be  P\(\sup_{1\le j\le m(n)} \frac{\wh\eta(d+t_j)-\wh\eta(d)}{\wh{\si}^\bb(d+t_j, d)(2\log\log 1/t_j)^{1/2}}\ge  1-\ep\) \ge     1-\ep,\label{80.10re}
\ee
for all $n$ sufficiently large.   
Using Theorem \ref{theo-1.2z} we see that (\ref{1.34sz}) holds for the permanental processes considered in Theorem \ref{theo-1.2zzk}. Using (\ref {80.7rr}) we see that this is the same as (\ref{80.6}).

\medskip To complete the proof of Theorems \ref{theo-1.2zz}--\ref{theo-1.2aa} it remains to show that either conditions $(i)$ or $(ii)$ of Theorem \ref{theo-1.2z} is satisfied.
 In Part I, condition $(ii)$ of these theorems we assume that for some $p>1$ and $\de>0$,
  \begin{equation}
  (\si^{\bb})^2(x) \ge x\( \log 1/x\)^{ p},\qquad x\in [0,\de], 
  \end{equation}
  for some $\bb\ge 0$.
  It follows from (\ref{3.7q}), which we. have shown holds for $(\si^{\bb})^2$ that
\begin{equation} 
   (\si^{\bb})^2(x)\le C\vf(x),\qquad\text{as $x\to 0.$} 
\end{equation}
This implies that, 
 \begin{equation} \label{5.10bbq}
\lim_{n\to\ff} m(n) \sup_{1\le j\neq k\le  m(n) }   |t_k-t_j| \(\vf\( |t_k-t_j|\)\)^{-1}  = 0,
\ee
 since, for $n+1-[n^{q}]\ge n/2$,   
 \bea 
 && m(n)    |t_k-t_j| \(\vf\( |t_k-t_j|\)\)^{-1}\\&&\qquad\le   n\sup_{|x|\le \th^{n/2}}\frac{x}{ \vf (x)}\le Cn \sup_{|x|\le \th^{n/2}}\frac{x}{ (\si^{\bb})^2(x)}\nn\\
  &&\qquad\le Cn\sup_{|x|\le \th^{n/2}}\frac{1}{  (\log 1/x)^p}=\frac{C2^{p}}{n^{p-1}\log^{p} (1/\th)}. \nn
\eea 
 Therefore,   condition  $(i)$  of Theorem \ref{theo-1.2z} is satisfied by the processes considered in Part I of Theorems \ref{theo-1.2zz}--\ref{theo-1.2aa}. 
 
\medskip Finally, it is clear that condition  $(ii)$  of Theorem \ref{theo-1.2z} is satisfied by the processes considered in Part II of Theorems \ref{theo-1.2zz}--\ref{theo-1.2aa} because the hypotheses in Part II of Theorems \ref{theo-1.2zz}--\ref{theo-1.2aa} simply repeat condition  $(ii)$  of Theorem \ref{theo-1.2z}. (The additional conditions in Part II of Theorems \ref{theo-1.2zz} and \ref{theo-1.2zzk} are used to show that (\ref{3.7q}) holds in these cases.)\qed

In Section \ref{sec-ex} we show that there are   many examples of the functions $(\si^\bb)^2$
that satisfy the hypotheses in Theorems \ref{theo-1.2zz}--\ref{theo-1.2aa}. In Section \ref{sec-ex1} we show the same thing for the excessive functions considered in these theorems.

    \subsection{Upper bounds   in Theorems \ref{theo-1.2zz}--\ref{theo-1.2aa}}    \label{sec-4.2}

  \medskip  The following theorem follows    from   \cite[Theorems 1.1 and  1.2]{MRejp}.

 \begin{theorem}  \label{theo-4.1aa} Let $\wt X_{\al}=  \{\wt X_{\al}(t),t\in  \De_d(\de)\}$, (see (\ref{deld})), be an $\al-$permanental process with kernel $\{v(x,y);x,y\in \De_d(\de)\} $
 and assume that,  
 \begin{equation} \label{upperbound.1}
  v(x,x)+ v(y,y) -2(v(x,y) v(y,x))^{1/2}\le \vf^{2}(|x-y|),\qquad x,y\in \De_d(\de),
\end{equation}
 for some $\de>0$ sufficiently small, where $\vf^{2}(|x|)$ is a regularly varying function at $0$  with   index $0<p\le 1$. Then     \begin{equation}
   \limsup_{  x\to 0}\frac{| \wt X_{\al}(d+x)- \wt  X_{\al}  (d)|}{  (2\vf^{2}(x)   \log\log 1/|x|)^{1/2}}\le \sqrt {2 }  X ^{1/2}_{\al} (d) \qquad a.s.\label{1.34sw}
   \end{equation} 
   \end{theorem}

   \Proof In \cite{MRejp} Theorem 1.1 is used to prove Theorem 1.2 which is the same as Theorem \ref{theo-4.1aa} except that $\De_d(\de)$ is replaced by $[0,1]$.
  The same proof gives the   statement in Theorem \ref{theo-4.1aa}.\qed  

   We   now use this theorem to get upper bounds 
 for the permanental processes in Theorems \ref{theo-1.2zz}--\ref{theo-1.2aa}. To begin we must estimate (\ref{upperbound.1}) for the kernels of the permanental processes in these theorems.
Let $X_{\al}=\{X_{\al}(t),t\in[d,d+\de]\}$,  $d\in T$, be an $\al-$permanental process with kernel 
\be
u _{f,g}(x,y)= u (x,y)+g(x)f(y),\qquad x,y \in[d,d+\de].\label{4.1nnq}
  \ee   
 Set,
 \begin{equation} \label{3.107a}
   \si  ^2(x,y) =u(x,x)+u(y,y)-2u(x,y)
\end{equation}
 and 
 \be \label{3.107}
  \wh\si_{ g,f} ^2(x,y) = u_{g,f}(x,x) +u_{g,f}(y,y) -2\(u_{g,f}(x,y) u_{g,f}(y,x)\)^{1/2}.
 \ee
 
  Note that by (\ref{poscond.8}) we can choose $\de$ so that,
\begin{equation} \label{3.111a}
  \sup_{x,y\in [d,d+\de]}u\(x,y\)>0.
\end{equation}
The next lemma generalizes part of  \cite[Corollary 4.6]{MRejp}.

 \begin{lemma} Let $d,\de$ be such that (\ref{3.111a}) holds. Then,
  if  $g,f\in C^1([d,d+\de])$,       \be \label{3.109}
  \wh\si_{g,f} ^2(x,y)\le     \si ^2(x,y) +C_{d,\de} |y-x|^2,\qquad x,y\in [d,d+\de]. 
 \ee
 for some constant $C_{d,\de}$ depending on $d$ and $\de$.

\end{lemma}

\Proof Add and subtract $u_{g,f}(x,y) +u_{g,f}(y,x)$ from the right--hand side of (\ref{3.107}) to see that,  
 \bea
\wh\si  _{g,f}^2(x,y)  &=&u_{g,f}  (x,x) +u_{g,f}  (y,y)\\
 &&   - u_{g,f}  (x,y) -u_{g,f}(y,x) +\( u _{g,f} ^{1/2}(x,y)- u _{g,f} ^{1/2}(y,x)\)^2\nn.
\eea
Note that,
\bea 
 && u_{g,f}  (x,x) +u_{g,f}  (y,y) - u_{g,f}  (x,y) -u_{g,f}(y,x)  \\
  &&\qquad =\si^2(x,y)+(f(x)-f(y))(g(x)-g(y))\nn\\
  &&\qquad \le \si^2(x,y)+ O|x-y|^2 ,\nn
\eea
and 
\bea 
  &&\( u _{g,f} ^{1/2}(x,y)- u _{g,f} ^{1/2}(y,x)\)^2 \\
  &&\qquad=\(\(u (x,y)+g(x)f(y)\)^{1/2}-\(u (y,x)+g(y)f(x)\)^{1/2}\)^2\nn .
\eea

 Using the simple fact that   for $a,b>0$, $|b^{1/2}-a^{1/2}|=|b-a|/(a^{1/2}+b^{1/2})$,
 we see that,  
 \begin{equation} 
  \(u _{g,f} ^{1/2}(x,y)- u _{g,f} ^{1/2}(y,x)\)  ^2\le \frac{ |g(x)f(y)-g(y)f(x)|^2}{2(a^*(d,\de))^2}=O\(|x-y|^2\),
\end{equation}
where,  
\begin{equation} 
 ( a^*(d,\de))^2=\inf_{x,y\in[d,d+\de]} u (x,y)+\inf_{x \in[d,d+\de]}(g(x)\wedge f(x)).
\end{equation}
Since $u^\bb$ is continuous and $u^\bb(d,d)>0$, (see (\ref{poscond.8})),  and $g,f\ge 0$, $a^*(d,\de)>0$, for all $\de$ sufficiently small.
 Putting all this together gives (\ref{3.109}).\qed

\medskip\noindent  \textbf{Proof of upper bounds in Theorem \ref{theo-1.2zz}   }  
In   this theorem $u^{\bb}$ is translation invariant. Therefore, $(\si^{\bb})^2 (x,y)=(\si^{\bb})^2 (y-x,0)$, which we write simply as $ (\si^{\bb})^2 (y-x)$. Therefore,  by (\ref{3.109}),   for any $\ep>0$,
 \begin{equation} \label{3.114}
 (\wh \si^{\bb})^2(x,y) \leq (1+\ep) (\si^{\bb})^2(x-y),\qquad\text{as $x,y\to d.$}
\end{equation} 
 Here we use the fact that in   this theorem,  
   \begin{equation} \label{equiv.5ee}
  \frac{|x-y|^2}{ (\si^{\bb})^2(x-y)}\to 0,\qquad\text{as $x,y\to d,\,\, x\neq y.$}
\end{equation} 
 It follows from (\ref{3.114}) that for any $\ep>0$, (\ref{upperbound.1})  holds for all $\de$ sufficiently small with    $\vf^2(|x-y|)=(1+\ep) (\si^{\bb})^2(x-y)$.  Using this in Theorem \ref{theo-4.1aa}  we get the upper bounds in   Theorem \ref{theo-1.2zz}  and see that they   hold for $\al-$permanental processes with kernel 
 $u^{\bb}_{g,f}(x,y)$ for all $\al>0$.\qed
  
  \noindent  \textbf{Proof of upper bounds in Theorems \ref{theo-1.2aa}  }  In this case $ (\si^{\bb})^2(x,y)$ in  (\ref{3.107a})
  is equal to $2\phi(x-y)$. Consequently the proof proceeds exactly as in the proof above  with $ (\si^{\bb})^2(x-y)$ replaced by $2\phi(x-y)$.\qed
  
  \noindent \textbf{Proof of upper bounds in Theorem \ref{theo-1.2zzk} }    In this case $ (\si^{\bb})^2(x,y)$ in  (\ref{3.107a})
  is equal to,
     \be 
  v^\bb(x,x)+v^\bb(x,x)-2v^\bb(x,y)  .\label{80.8q}
 \ee  
  By Lemma \ref{lem-80},   
  \begin{equation} \label{80.7ems}
\frac{  v^\bb(x,x)+v^\bb(x,x)-2v^\bb(x,y)}{ (\si^{\bb})^2(x-y)}  \to 1    , \qquad\mbox{as } \,\,x,y\downarrow d,\,\, x\neq y.
\end{equation}
  Therefore the proof follows as in the proof of the upper bounds in Theorem \ref{theo-1.2zz}. \qed

    \begin{remark} \label{rem-4.1}{\rm  Let $\eta$ be a Gaussian process with stationary increments and $\eta(0)=0.$  For a very large class of Gaussian processes;   see e.g.,  \cite[Theorem 7.2.15]{book},
 \begin{equation}
  \limsup_{  t \downarrow 0}\frac{|  \eta(  t)    |}{ \( 2\si^{2}\(t\)\log\log 1/t\)^{1/2}}=  1\qquad a.s.,\label{1.3aa}
  \end{equation}
where $
\si^{2}\(t\)=  E(\eta^2(  t) ).
$
  Consequently,
\bea
  \limsup_{  t \downarrow 0}\frac{   (\eta^2(  t) - \eta^2(  0))/2}{ \(  \si^{2}\(t\)\log\log 1/t\) }&=&  1\qquad a.s.,\label{1.3bbq} 
  \eea
This is the result we expect to get for permanental process with kernels of the form of (\ref{4.1nnq}) when $u(x,y)$ is the covariance of $\eta.$

}\end{remark}

  The following remark relates to statements made on pages \pageref{ref-3.2a} and \pageref{rem3.2}.

\begin{remark} \label{rem-3.2}{\rm  Let $\xi =\{\xi(x);x\in R^1 \}$ be a Gaussian process with    $\xi(0)=0$ and the property that,
\begin{equation} \label{3.118}
  E \( \xi (x)- \xi (y)    \)^2=\rho(x-y),\qquad x,y \in R^1, 
\end{equation} 
  for some positive function $\rho$.  It is well known and easy to show that $\xi$ has stationary increments.
  Using the equality,
  \begin{equation} 
\frac{1}{2}  (a-b)(c-d)=-(a-c)^2+(a-d)^2+(b-d)^2-(b-c)^2,
\end{equation}
 we see that since for any $x,y,z,w\in R^1$,
 \begin{eqnarray}
 &&\frac12 E\(\( \xi (x)- \xi (y)    \) \( \xi (z)- \xi (w)    \)  \)
 \label{stat}\\
 &&\qquad  =-\rho(x-z)+\rho(x-w)+\rho(y-z)-\rho(y-w),  \nonumber
 \end{eqnarray}
 which shows that the covariance of the increments of $\xi$ are translation invariant.  Since the joint distributions Gaussian process are determined by it's covariance we see that $\xi$ has stationary increments.
 We might also note that (\ref{3.118}) also shows that,
 \begin{equation} 
  E(\xi(x)\xi(y))=\frac{1}{2}\(\rho(x)+\rho(y)-\rho(x-y)\).
\end{equation}
}\end{remark}

\subsection{Proof of Theorem  \ref{theo-1.5} } 

Note that the conditions imposed on the excessive functions $f$
and $g$ in Theorem  \ref{theo-1.5} are the same as those in condition $(ii)$ in Theorem \ref{theo-1.2z}. Therefore to prove this theorem it remains to show that
 the potentials in (\ref{diff.1}), (\ref{diff.22}) and (\ref{diff.5}) satisfy the conditions in (\ref{3.7q}) and (\ref{2.149}) and that (\ref{1.34sz}) holds for the Gaussian processes with covariances that are these potentials. This follows immediately from the following computations:
 
 \medskip   For $x,y\in R^{1}$ define, 
 \begin{equation} \label{diff.1q}
\wt u^\bb (x,y)= \left\{
 \begin{array} {cc}
 p(x)q(y),& \quad x\leq y  
 \\
 q(x)p(y),& \quad y\leq x   
\end{array}  \right. ,
\end{equation}
where   $p,q\in C^{2}(R^{1})$    are positive convex functions such that $p$ is  strictly increasing and  
$q$ is  strictly  decreasing.    Note that
 \begin{equation} \label{diff.2}
 \wt u^\bb (x,x)-\wt u^\bb (x,y)= \left\{
 \begin{array}  {cc}
 p(x)\( q(x)-q(y)\) ,& \quad x\leq y  
 \\
 q(x)\(p(x)-p(y)\),& \quad y\leq x   
\end{array}  \right. .
\end{equation}
   It follows from \cite[Lemma 5.1.8]{book} that $\{\wt u^\bb (x,y),x ,y\in R^1 \}$ is a strictly positive definite function on $R^1\times R^1$. Let $\wt \eta=\{\wt\eta(x),x\in R^{1} \}$  be a Gaussian process with covariance $\wt u^\bb (x,y)$ and set,
\begin{equation} 
  (\wt \si^\bb ) ^{2}(x,y)=E\(\wt \eta(x)-\wt \eta(y)\)^2, 
\end{equation}
so that
\begin{equation} \label{diff.3}
   (\wt \si^\bb )  ^{2}(x,y)= \left\{
 \begin{array}  {cc} 
 p(x)\( q(x)-q(y)\)+\( p(y)-p(x)\)q(y),& \quad x\leq y  
 \\
 q(x)\(p(x)-p(y)\)+p(y)\( q(y)-q(x)\),& \quad y\leq x.   
\end{array}  \right. 
\end{equation}

\begin{lemma} \label{lem-4.9} For   $x,y\in R^{1}$,
 \be  \label{diff.14q}
  (\wt \si^\bb )  ( x, y)\sim
 \tau(x)|x-y|, \qquad \text{as $|x-y|\to 0$},
\ee
where, 
 \begin{equation} \label{8.51}
 \tau(x)=\wt u^\bb(x,x)\frac{d}{dx}\(\log \frac{ p(x)}{q(x)}\)  >0 .
\end{equation} 
In addition for $y\le x$,
\begin{equation} \label{4.111}
  \wt u^\bb(x,x)-\wt u^\bb(x,y)\sim q(x)p'(x)|x-y|, \qquad \text{as $|x-y|\to 0$}. 
\end{equation}
\end{lemma}

\Proof Note that,   
 \bea \label{4.112}
 q( x)-q( y)&\sim & q'(x ) (x-y),\qquad \text{as  $|x-y|\to 0$} \\
  p( y)-p( x)&\sim & p'(x ) (y-x),\qquad \text{as  $|x-y|\to 0$}\nn.
\eea
Therefore,  by (\ref{diff.3}), when  $  x\leq y$,   as $|y-x|\to 0$, 
\begin{eqnarray}
\label{diff.14} 
 (\wt \si^\bb )^{2} ( x, y)&  =& p( x)\( q( x)-q( y)\)+q( y)\( p( y)-p( x)\) \nonumber\\
& \sim &p(x )   q'(x ) (x-y)  + (q(x)-q'(x ) (x-y))p'(x)(y-x)\nn\\
&\sim& \(q(x)p'(x)-p(x )   q'(x )\)(y-x)=\tau(x)(y-x). \nonumber
\end{eqnarray}
Since $(\wt \si^\bb )^{2}$ is symmetric we get  (\ref{diff.14q}).    Furthermore, since $p(x)/q(x)$ is strictly increasing we see that (\ref{8.51}) holds.

  When  $y\le x$,
\begin{equation} 
  \wt u^\bb(x,x)-\wt u^\bb(x,y)= q(x)(p (x)-p(y))\sim q(x)p'(x)|x-y|, \quad \text{as $|x-y|\to 0$}, 
\end{equation}
which is (\ref{4.111}).\qed

\begin{lemma} \label{lem-4.10}
The increments of  $\wt\eta $ are negatively correlated, i.e., for  $ x<y<z<w$,
\begin{equation} 
 E(\wt\eta(y)-\wt\eta(x))(\wt\eta(w)-\wt\eta(z))\le 0.
\end{equation}
Consequently, for $d,x,y\in T$ where  $  x\le y $,
 \begin{equation} \label{diff.10}
\frac{  E\(\(\wt\eta(d+x)-\wt\eta(d)\)\(\wt\eta(d+y)-\wt\eta(d)\)\) }{ (\wt \si^\bb)    (d,d+x)(\wt \si^\bb)   ( d,d+y)}\le \frac{ (\wt \si^\bb)   (d,d+x)}{ (\wt \si^\bb)     (d,d+y))}.  \end{equation}
\end{lemma}

\Proof If follows from (\ref{diff.1q}) that
\begin{equation} 
 E(\wt\eta(y)-\wt\eta(x))(\wt\eta(w)-\wt\eta(z))=(p(y)-p(x))(q(w)-q(z))\le 0.
\end{equation}

\qed

  Let
\begin{equation}
\wt v^\bb  (x,y)=\wt u^\bb (x,y)-\frac{\wt u^\bb (x,0)\wt u^\bb (0,y)}{\wt u^\bb (0,0)} . \label{diff.22w}  
\end{equation}

This is the covariance of  the Gaussian process, $\wt\xi=\{\wt\xi(x),x\in R^1 \}$ where,
\begin{equation} \label{4.119}
\wt \xi(x)=\wt\eta(x)-\frac{\wt u^\bb(x,0) }{\wt u^\bb(0,0)}\wt \eta(0).
\end{equation}
Note that $\wt\xi(0)=0$. 

 Set,
\begin{equation} \label{9.23}
 (\ov\si^\bb)  ^{2}(x,y):=E\(\wt \xi(x)-\wt \xi(y)\)^2 =  (\wt\si^\bb)  ^{2}(x,y)-\frac{\(\wt u^\bb (x,0)-\wt u^\bb(0,y)\)^2}{\wt u^\bb (0,0)}.
\end{equation}

 \begin{lemma} \label{lem-4.11} For $x,y>0$,
 \be  \label{9.7aa}
 (\ov\si^\bb)  ^{2}(x,y)\sim  (\wt\si^\bb)  ^{2}(x,y)\sim
 \tau(x)|x-y|
 , \qquad \text{as $\,|x-y|\to 0$}.
\ee
 When $0< y\leq x$,
\be  \label{4.122}
  \wt v^\bb(x,x)- \wt v^\bb(x,y)    \sim  C_{x}|x-y|, \qquad \text{as $\,|x-y|\to 0$}, 
\ee  
where $C_{x}>0.$ 
\end{lemma}

\Proof  Consider (\ref{9.23}) and note that,
\begin{equation} 
 \frac{\(\wt u^\bb(x,0)-\wt u^\bb(0,y)\)^2}{\wt  u^\bb(0,0)}=\frac{p(0)}{q(0)}(q(x)-q(y))^2.
\end{equation}
When $0<y\le x$,
\begin{equation} 
\frac{p(0)}{q(0)}(q'(x))^2(x-y)^2
\le \frac{\(\wt u^\bb(x,0)-\wt u^\bb(0,y)\)^2}{u^\bb(0,0)}\le\frac{p(0)}{q(0)}(q'(y))^2(x-y)^2.
\end{equation}
Using this and (\ref{diff.14q}) we get (\ref{9.7aa}). 

  When $0<y\le x$,  
\bea  
  \wt v^\bb(x,x)- \wt v^\bb(x,y) &=&\wt u^\bb(x,x)- \wt u^\bb(x,y) - \frac{\wt u^\bb (x,0)(\wt u^\bb (x,0)-\wt u^\bb (0,y))}{\wt u^\bb (0,0)} \nn\\
  &=&\wt u^\bb(x,x)- \wt u^\bb(x,y) - \frac{p(0)}{q(0)}q (x)(q (x)- q(y)).  \eea
 As   $ |x-y|\to 0$  this is,
  \bea 
 &  \sim& q(x)p'(x)|x-y|- \frac{p(0)}{q(0)}q (x) q'(x)|x-y|   \nonumber\\
&  = &\( p'(x) - \frac{p(0)}{q(0)} q'(x)\) q (x)   |x-y|,\nonumber
\end{eqnarray}  
where we use (\ref{4.112}). Since $p'>0$ and 
$q'>0$ we get (\ref{4.122}).\qed

  \begin{lemma} \label{lem-9.2a}
The increments of  $\wt\xi $ are negatively correlated, i.e., for $0<x<y<z<w$,
\begin{equation} \label{9.29}
 E(\wt\xi(y)-\wt\xi(x))(\wt\xi(w)-\wt\xi(z))\le 0.
 \end{equation}
 Consequently, for $d,x,y\in T$ where  $0< x\le y $,
 \begin{equation} \label{diff.10q}
\frac{  E\(\(\wt\xi(d+x)-\wt\xi(d)\)\(\wt\xi(d+y)-\wt\xi(d)\)\) }{ (\ov\si)^{\bb}    (d,d+x)(\ov\si)^{\bb}  ( d,d+y)}\le \frac{(\ov\si)^{\bb}  (d,d+x)}{(\ov\si)^{\bb}     (d,d+y))}.  \end{equation}

\end{lemma}

\Proof
\bea 
\wt \xi(y)-\wt\xi(x)&=& \wt\eta(x)-\wt\eta(y)-\frac{\wt u^\bb(0,y)-\wt u^\bb(0,x)}{\wt u^\bb(0,0)}\wt\eta(0)\\
  \wt \xi(w)-\wt\xi(z)&=& \wt\eta(w)-\wt\eta(z)-\frac{\wt u^\bb(0,w)-\wt u^\bb(0,z)}{\wt u^\bb(0,0)}\wt\eta(0)\nn 
\eea
Therefore, using Lemma \ref{lem-4.10}, the left--hand side 
of (\ref{9.29}) is equal to, 
\bea \label{9.31}
 &&E(\wt\eta(y)-\wt\eta(x))(\wt\eta(w)-\wt\eta(z))-\frac{(\wt u^\bb(0,y)-\wt u^\bb(0,x))(\wt u^\bb(0,w)-\wt u^\bb(0,z))}{\wt u^\bb(0,0)} \nn\\
 &&\qquad\le 0-\frac{p (0)}{q (0)}\(q(y)-q(x)\)\(q(w)-q(z)\) \le 0.
\eea
\qed

Let
 \begin{equation} \label{diff.71}
w(x,y)=s(x)\wedge s(y),\qquad  x,y>0.
\end{equation}
where $s\in C^{2}(R^+-\{0 \})$ and  is positive  and  strictly increasing. 
Therefore, for $x\le y$  
\begin{equation} 
  w(x,x)-w(x,y)=0
\end{equation}
and for $y\le x$
\begin{equation} \label{4.134}
  w(x,x)-w(x,y)=s(x)-s(y)\ge 0.
  \ee

Obviously, $w(x,y)$ is the covariance of a time changed Brownian motion   $B=\{B(s(x)),x\in R^+-\{0 \}) \}$.
Therefore,
\begin{equation} 
 \wh\si^2(x,y):=E\(B(s(y))=B(s(x))\)^2=|s(x)-s(y)|.
\end{equation}  
and
\begin{equation}
\wh \si ^{2}(x,y)\sim s'(x) |x-y|, \qquad \text{as $|x-y|\to 0.$}\hspace{.2 in}\label{diff.4aa}
\end{equation}

Since $B$ has independent increments we have: 
\begin{lemma} \label{lem-4.13} For $d,x,y\in T$ where  $0< x\le y $,
  \begin{equation} \label{diff.20}
\frac{  E\(\(B(d+x)-B(d)\)\(B(d+y)-B(d)\)\) }{\wh \si   (d,d+x)\wh\si   ( d,d+y)}=\frac{ \wt\si   (d,d+x)}{ \wh\si    (d,d+y))}.  \end{equation} 
\end{lemma}

  We now note that it follows from Lemmas \ref{lem-4.9} and \ref{lem-4.11} and (\ref{diff.4aa}) that we can take,
\begin{equation} 
  \vf_d(|x-y|)=|x-y|
\end{equation}
for $(\wt \si^\bb)^2$, $(\ov \si^\bb)^2$ and $\wh \si^2 $ in the first line of (\ref{3.7q}) and, obviously (\ref{2.149}) is satisfied. The second line of (\ref{3.7q}) for each of the different potentials follows from   (\ref{4.111}),   (\ref{4.122}), and (\ref{4.134}). \qed

\medskip The proof that (\ref{1.34sz}) holds for the Gaussian processes with covariances that are the potentials that we are studying follows from the inequalities in Lemmas \ref{lem-4.10} and \ref{lem-4.13} and (\ref{4.119}) exactly as in the proof of Theorems \ref{theo-1.2zz}--\ref{theo-1.2aa}.
This completes the proof of the lower bounds in Theorem \ref{theo-1.5}.

\medskip The proof of the upper bounds in Theorem \ref{theo-1.5} is exactly the same as the proof of the upper bounds in Theorem \ref{theo-1.2aa}. For the processes we are studying $\phi(x)=C|x|$, for some constant $C$.

\section{  Exact moduli of continuity for Markov local times}\label{sec-emodMP}

Let   $ Z=
(\Om,  \FF_{t}, Z_t,\th_{t}, P^x)$ be a  Borel right process  on a locally compact metric state space $(S,d)$ 
strictly positive $0$--potential densities
$\{v(x,y),x,y\in S\}$ with respect to a $\si$--finite measure $m$.  Let 
$  L=\{  L^y_{t}\,,\,(y,t) \in S\times R_+\}$
denote the local   times for
$  Z$  normalized so that,
\begin{equation}
  E^{ x}\(  L^{y}_\ff\)=v(x,y). \label{az.1}
\end{equation}
  Let $X_{1/2}=\{X_{1/2}(y),y\in S\}$ denote a $1/2$--permanental process with kernel $\{v(x,y),x,y\in S\}$. 

We assume that $  L$ is jointly continuous on $S\times [0,\ze)$, where $\ze$ is the lifetime of $Z$
and that $X_{1/2}$  is continuous almost surely. 

\medskip  The next    theorem uses   the Eisenbaum--Kaspi isomorphism theorem \cite{EK} and the local   modulus of continuity of $X_{1/2}$ to obtain the local  modulus of continuity of $L$.

\bt\label{theo-LILMP} Let $\phi $ be a positive continuous real valued function. Suppose that for some $v\in (S,d)$,
\begin{equation}
 \limsup_{  u \to v}\frac{|  X_{1/2}( u )-  X_{1/2}  (v )|}{\phi(d(u,v))}=    \(      X _{1/2} (v )\)^{1/2}, \qquad  a.s.\label{az.2}
 \end{equation} 
Then
\begin{equation}
 \limsup_{  u \to v}
\frac{ | L^{u}_t-  L^{v}_t|}{\phi(d(u,v))} =
   \(  L^{v}_t\)^{1/2},\qquad \text{for a.e. $t<\ze$}\label{az.3},
\end{equation}
  $ P^{x}$ almost surely, for all $x\in S$. 
 
\et

We state the Eisenbaum--Kaspi theorem for the convenience of the reader:

\bt [Eisenbaum and Kaspi,   \cite{EK}]\label{theo-ke}  Let    $ Z=
(\Om,  \FF_{t}, Z_t,\th_{t}, P^x)$ be a  Borel right process in $S$ with 
$0$--potential density
$v(x,y)$. Let $h_{x}(z)=v(z,x)$ and assume that   $h_{x}(z)>0$ for all $x,z\in S$. Let 
$L=\{L^y_{t}\,,\,(y,t) \in S\times R_+\}$
denote the local   times for
$  Z$, normalized so that
\begin{equation}
 E^{ x}\( L^{y}_\ff\)=v(x,y), \label{mp4.8.5 }
\end{equation}
as defined above. Let $X_{1/2}=\{X_{1/2}(y),y\in S\}$ denote a $1/2$--permanental process with kernel $v(x,y)$.  Then, for any
  countable subset $D\subseteq S$,  
\be
\Big\{ L^u_{\ff}+ X_{1/2}(u),u\in D,P^{y/h_{y}}\times
P_{X_{1/2}}\Big\}\stackrel{law}{=}
\Big\{X_{1/2}(u),u\in D\,,\,\frac{2X_{1/2}(y)}{v(
y,y)}P_{X_{1/2}}\Big\},\label{it1.8}
\ee

Equivalently, for all $x_{ 1},\ldots,x_{ n}$ in $S$ and bounded 
measurable functions
$F$ on $R^n_+$, for all
$n$,
\be
     E^{y/h_{y}}E_{X_{1/2}} \(F\(L^{x_{ i}}_{\ff}+ X_{1/2}(x_{ i})\)\)
=
  E_{ X_{1/2}}\(\frac{2X_{1/2}(y)}{v(
y,y)}F\(X_{1/2}(x_i)\)\),\label{it1.9x} 
   \ee 
 where we use the notation $F( f(x_{ i})) =F( f(x_{ 1}),\ldots, 
f(x_{ n}))$. 
\et

The   probability measures $P^{x/h_y}$ are the probability measures for the $h_{y}$--transform of $Z$,  see \cite[(62.24)]{S} or \cite[(2.22)]{FR}. In this proof we express them using Azema's formula \cite{Azema}, as  \be
 P^{x/h_y}(A) =\frac{1}{h_y(x)}  P^x\left(\int_0^\infty 1_A\circ k_s\,d_s  L_s^y\right),\qquad A\in \Om,
\label{az.4}
\ee
 where $k_s$ is the killing operator for $Z$, i.e., $k_s\om (r)=\om (r)$ when $r<s $, and $k_s\om (r)=\De$ when $r\geq s $.

\medskip\noindent {\bf  Proof of Theorem \ref{theo-LILMP}  }   
 Let $F( \mathcal{C})$  be a set of  functions on a countable dense   subset $\mathcal{C}$ of $S$ and let, 
\be
A  =\Big\{g\in F( \mathcal{C})\,\Big |  \,\limsup_{u\in \mathcal{C},   u\to v}\frac{| g(u)-  g(v)|}{\phi(d(u,v))}=  \({g(v)}\)^{1/2} \Big\}.
\ee
Since $X_{1/2}$ is continuous   and satisfies (\ref{az.2}) it follows that,   
\be
  P_{X_{1/2}}\(X_{1/2} \in A\)=1.\label{az.2p}
\ee 
Note that it follows from   (\ref{int.1}) that   $E_{X_{1/2}}\(X_{1/2}(y)\)= v(y,y)/2$.  
Therefore  by Theorem \ref{theo-ke},  for almost all $\om'\in
\Om_{X_{1/2}}$, the probability space of $X_{1/2}$,  
\be  \limsup_{ u\in \mathcal{C},   u\to v}
\Bigg|\frac{  L^{u}_\ff- L^{ v}_\ff}{\phi(d(u,v))}\!+  \frac{   X_{1/2}( u,\om')- X_{1/2}  (v,\om') }{\phi(d(u,v))} \Bigg| =
  \,  \({L^{ v}_\ff+X_{1/2}(v,\om')}\)^{1/2}, \label{az.7}  
\ee
 $  P^{y/h_{y}}$ almost surely, for all $y\in S$. 
 Using (\ref{az.2p}) we see that  for almost all
$\om'\in
\Om_{X}$  with respect to $P_{X}$,
\begin{equation}
 \limsup_{u\in \mathcal{C}, u\to v}\frac{|  X_{1/2}( u,\om')-  X_{1/2}  (v,\om')|}{\phi(d(u,v))}=    \(      X _{1/2} (v,\om')\)^{1/2} .  \label{az.8}
 \end{equation} 
 As we point out immediately after the statement of Theorem \ref{theo-1.2zz},     $X_{1 /2} (d)$ is the square of a mean zero Gaussian random variable. Consequently,  it takes values arbitrarily close to 0 with probability greater than 0.  Therefore it follows from (\ref{az.7}) and the triangle inequality that, 
\begin{equation}
\limsup_{ u\in \mathcal{C},   u\to v}
\Bigg|\frac{  L^{u}_\ff-  L^{ v}_\ff}{\phi(d(u,v))} \Bigg|=     \(   L^{ v}_\ff\)^{1/2},\label{az.9}
\end{equation}
 $ P^{y/h_{y}}$ almost surely, for all $y\in S$.
 (A similar argument is used in     \cite[Theorem 9.5.1]{book}.)   

   \medskip    For $t\in R^1_+$, including   $t=\ff$, set
\begin{equation}
D_{t} =\Bigg\{  \limsup_{ u\in \mathcal{C},   u\to v}
\Bigg|\frac{   L^{u}_t-   L^{ v}_t}{\phi(d(u,v))} \Bigg|=    \,  \(   L^{ v}_t\)^{1/2}  \Bigg \}. \label{az.10}
 \end{equation}
It follows from  (\ref{az.9})  that
\be
P^{y/h_y}\left(D_{\ff} \right)=1.
\label{az.11}
\ee
We now show that (\ref{az.11}) implies that  
  \begin{equation} \label{az.17}
  \limsup_{ u\in \mathcal{C},  \, u\to v}
 \frac{|  L^{u}_s-  L^{ v}_s|}{\phi(d(u,v))} =\(    L_{s}^{v}\)^{1/2}, 
\end{equation}
 for  almost every $s<\ze$, $ P^{x}$ almost surely, for all $x\in S$.   Then, since  by hypothesis   $  L$ is jointly continuous on $S\times [0,\ze)$,  
we   get (\ref{az.3}).

 To obtain (\ref{az.17}),   we first note that by (\ref{az.4}) and  (\ref{az.11}),
\be
\frac{1}{h_y(y)}  P^y\left(\int_0^\infty 1_{\{ D^{c}_{\ff}\}}\circ k_s\,d _s L_s^y\right)=0.\label{az.4qq}
\ee
Since $1_{\{D^{c}_{\ff}\}}\circ k_s=1_{\{D^{c}_{s}\}}$ and  $h_y(y)>0$ this implies  that,
\be
  P^y\(\int_0^\infty 1_{\{D^{c}_{s}\}}\,d_{s} L^y_s\)=0.\label{az.12}
\ee

 Let $T_y$ denote the first hitting time of $y$ and $\th_{T_{y}}$ denote the shift operator. Then for any $x\in S$,
  \begin{eqnarray}
    P^{x}\(\int _{0}^{\ff} 1_{\{D^{c}_{s}\}} \,d_{s}  L^{y}_{s} \) 
  \label{az.13}  &=&  P^{x}\(\int _{T_{y}}^{\ff} 1_{\{D^{c}_{s}\}} \,d_{s}  L^{y}_{s} \)   \\
  &=&    P^{x}\( \(\int _{0}^{\ff} 1_{\{D^{c}_{s}\}} \,d_{s}  L^{y}_{s} \)\circ\th_{T_{y}}  \)  \nonumber\\
  &=&    P^{x}\(T_{y}<\ff \)   P^{y}\(\int _{0}^{\ff} 1_{\{D^{c}_{s}\}} \,d_{s} L^{y}_{s} \)=0, \nonumber
  \end{eqnarray}
by (\ref{az.12}). %This shows that
%\begin{equation}
%P^{x}\(\int _{0}^{\ff}e^{-ps}1_{\{B^{c}_{s}\}} \,ds \) =\int _{R^{1}}P^{x}\(\int _{0}^{\ff}e^{-ps}1_{\{B^{c}_{s}\}} \,d_{s}L^{y}_{s} \)\,dy =0.\label{dt.14}
%\end{equation}
 
    We now examine  (\ref{az.13}). 
We first note that by the occupation density formula \cite[(2.28)]{FR},
\begin{equation}
\int_{S}    L^{y}_{s}   \,dm (y)=s, \hspace{.2 in}\forall s<\ze \quad a.s. \label{az.14}
\end{equation}
Therefore, using (\ref{az.13}) we have,\begin{equation}
 P^{x}\(\int _{0}^{\ze} 1_{\{D^{c}_{s}\}} \,ds \)=\int_S   P^{x}\(\int _{0}^{\ff} 1_{\{D^{c}_{s}\}} \,d_{s}  L^{y}_{s} \) \,dm (y)=0.\label{az.15}
\end{equation}
  Here we have used the fact that for any bounded  measurable
$f(\om, s)$
\begin{equation}
 P^{x}\(\int _{0}^{\ze} f(\om, s) \,ds \)=\int_S   P^{x}\(\int _{0}^{\ff} f(\om, s) \,d_{s}  L^{y}_{s} \) \,dm (y).\label{az.15a}
\end{equation}
 To see this, as in the proof of Fubini's theorem,  it suffices to prove it for $f(\om, s)$ of the form $g(\om)h(s)$, in which case, using Fubini's theorem, (\ref{az.15a}) becomes
\begin{equation}
 P^{x}\(  g(\om)   \int _{0}^{\ze} h( s) \,ds \)=  P^{x}\(    g(\om) \int_S  \int _{0}^{\ff} h( s) \,d_{s}  L^{y}_{s} \,dm (y)\).\label{az.15am}
\end{equation}
Therefore, to obtain (\ref{az.15a})
it only remains to show that,  
\begin{equation}
 \int _{0}^{\ze} h( s) \,ds =\int_S   \int _{0}^{\ff} h( s) \,d_{s}  L^{y}_{s}  \,dm (y), \qquad a.s.\label{az.15b}
\end{equation}
  Moreover, to show this it suffices to verify it for functions of the form $h( s)=1_{[0,a]} (s)$ with $a<\ze$. This is obvious, since by (\ref{az.14}),
 \be  
  a= \int_S       L^{y}_{a}  \,dm (y)=\int_S   \int _{0}^{a}   \,d_{s}  L^{y}_{s}  \,dm (y) 
  =\int_S   \int _{0}^{\ff} 1_{[0,a]} (s) \,d_{s}  L^{y}_{s}  \,dm (y).  
\ee 
   Using  (\ref{az.15})  we now see  that,
\begin{equation} 
  \int_0^\zeta 1_{\{D_s^c \}}\, ds = 0 ,\qquad    P^x \quad a.s. 
\end{equation} 
Therefore, for $ P^x$   almost every $\omega$ we have $\omega\in D_s$ for  almost every $s\in[0,\zeta(\omega))$, 
which is (\ref{az.17}). \qed

    \begin{remark} {\rm 
 Considering the importance of  Az\'ema's formula in the proof of Theorem \ref{theo-LILMP}, we show that Az\'ema's formula,  
\be
 P^{x/h_y}(A) =\frac{1}{h_y(x)}  P^x\left(\int_0^\infty 1_A\circ k_s\,d_{s}  L_s^y\right),
\label{dt.400}
\ee
satisfies  the more common definition of   $P^{y/h_y}(A)$;   i.e.,  for $A_{t}\in \mathcal{F}_{t}$,
  
\begin{equation}
 P^{x/h_y}(A_{t}1_{\{t<\ze\}}) =\frac{1}{h_y(x)}  P^x\left(A_{t} h_y(Z_{t})\right).\label{dt.401}
\end{equation}
This follows since 
\begin{equation}
\( A_{t}1_{\{t<\ze\}}\)\circ k_s= A_{t}1_{\{t<s\wedge\ze\}}\label{dt.402},
\end{equation}
and $ \int_t^\infty \,d_{s}  L_s^y=L_\ff^y-L_t^y=L_\ff^y\circ \th_{t}$,
so that
\begin{eqnarray}
&&P^x\left(\int_0^\infty 1_{A_{t}}\circ k_s\,d_{s}  L_s^y\right)=P^x\left(1_{A_{t}} \int_t^\infty \,d_{s}  L_s^y\right) \label{dt.403}
\nn\\
&&\qquad =P^x\left(1_{A_{t}} L_\ff^y\circ \th_{t} \right)=P^x\left(1_{A_{t}} E^{Z_{t}}(L_\ff^y)\right)  =P^x\left(A_{t} h_y(Z_{t})\right).\nn
\end{eqnarray}
}\end{remark} 

  The assumption in Theorem \ref{theo-ke} that  $Z$  has a finite $0$--potential density   implies that $Z$ is transient. The next theorem gives the local modulus of continuity for  recurrent processes.

   \begin{theorem}\label{cor-LILMP} Let   $ \wt Z\!=\!
 ( \wt \Om, \wt \FF_{t}, \wt Z_t , \wt P^x
 )$ be a recurrent  Borel right process with locally compact metric state space $(S,d)$, and 
strictly positive $p$--potential density
$v(x,y)$ with respect to a $\si$--finite measure $m$ on $S$,   where $p>0$.   Let 
$ \wt  L=\{ \wt  L^y_{t}\,,\,(y,t) \in S\times R_+\}$
denote the local   times for
$\wt   Z$, normalized so that
\begin{equation}
  \wt    E^{ x}\(\int_{0}^{\ff} e^{-pt} d_{t}\wt L^{y}_t\)=v(x,y), \label{az.1c}
\end{equation}
and assume that $ \wt  L$ is jointly continuous on $S\times R_+$.  

 Let $X_{1/2}=\{X_{1/2}(y),y\in S\}$ denote a $1/2$--permanental process with kernel $v(x,y)$ and assume that $X_{1/2}$  is continuous almost surely. Then,
 if (\ref{az.2}) holds for $X_{1/2}$ then (\ref{az.3}) holds for $\wt L$ for almost every  $t\in R_+$,   $ \wt  P^{x}$ almost surely, for all $x\in S$. 
\end{theorem}

   \Proof  
 Let  $ Z=
(\Om,  \FF_{t}, Z_t,\th_{t}, P^x)$ be  a Borel right process obtained by killing $ \wt Z$ at an independent
exponential time $\la$ with mean $1/p$.  Since the 0--potential densities for $Z$ are the $p$--potential densities
for $\wt Z$, $ Z$ is transient
Borel right process with 
  continuous strictly positive 0--potential densities $v(x,y)$.

The process $ Z$ is defined on a product space $\Om=\wt \Om\times [0,\infty)$ with probability $ P^x=\wt  P^x\times \ov P$ where $\ov P$  has probability density $pe^{-ps}$
  and 
  \be\label{5.34}
   \begin{array}{ll }
   Z_t (\omega,s)=\wt Z_t (\omega),&\qquad t<s   \\
 Z_t (\omega,s)=\De,& \qquad t\ge s  
\end{array}   .
\ee
Let $\ze (\omega,s)$ denote the lifetime and $  L^y_t$ the local time of $  Z$. We see from (\ref{5.34}) that $\ze (\omega,s)=s$ which implies that $ L^y_t(\omega,s)=\wt L_{t\wedge s}^y(\omega)$. (This is \cite[Remark 3.6.4, (3)]{book}). 
  Therefore, by Theorem \ref{theo-LILMP},
\begin{equation}
 \limsup_{  u \to v}
\frac{ | L^{u}_t(\omega,s )-   L^{v}_t(\omega,s)|}{\phi(d(u,v))} =
   \(   L^{v}_t(\omega,s)\)^{1/2},\qquad \text{for a.e. $t<s$}\label{az.3mar},
\end{equation}
  $ \wt P^{x}\times pe^{-ps}\,ds$ almost surely, for all $x\in S$. Using    Fubini's Theorem we then see that   for a.e. $s$,
\begin{equation}
 \limsup_{  u \to v}
\frac{ | L^{u}_t(\omega,s )-   L^{v}_t(\omega,s)|}{\phi(d(u,v))} =
   \(   L^{v}_t(\omega,s)\)^{1/2},\qquad \text{for a.e. $t<s$}\label{az.3marc},
\end{equation}
  $ \wt P^{x}$ almost surely, for all $x\in S$.   Using this for  a sequence $s_{n}\to \ff$ shows that (\ref{az.3}) holds for $\wt L$ for almost every $t$,   $ \wt  P^{x}$ almost surely, for all $x\in S$.  
  \qed
 
    An analogue of  Theorem \ref{theo-LILMP} also holds for the uniform modulus of continuity:
  \begin{theorem} \label{theo-5.2newa}
  Let $\xi $ be a positive continuous real valued function. Suppose that there exists a compact  set $K\subset S$ such that,  
 \begin{equation}
 \lim_{h\to 0}\sup_{\stackrel{u,v\in K}{d(u,v)\leq h}}\frac{|  X_{1/2}( u )-  X_{1/2} (v )|}{\xi(d(u,v))}= \sup_{v\in K}   \(     X _{1/2} (v )\)^{1/2}, \qquad  a.s.\label{az.2ka}
 \end{equation} 
and  for all $\ep>0$, 
\begin{equation} \label{5.5newa}
    P\( \sup_{v\in K} X_{1/2}  (v)  <\ep \)>0.
\end{equation} 
Then,
\begin{equation}
 \lim_{h\to 0}\sup_{\stackrel{u,v\in K}{d(u,v)\leq h}}
\frac{ | L^{u}_t-  L^{v}_t|}{\xi(d(u,v))} =
   \sup_{v\in K}   \( L^{v}_t\)^{1/2},\qquad \text{for a.e. $t<\ze$}\label{az.3k},\end{equation}
 $ P^{x}$ almost surely, for all $x\in S$.
\end{theorem}

\Proof   The proof,   with obvious modifications, is exactly the same except for one point. In the proof of Theorem \ref{theo-LILMP} we use the fact that $X_{1 /2} (d)$ is the square of a mean zero   Gaussian random variable to see that for all $\ep>0$, $P\(X_{1/2}(d)<\ep\)>0$. Here, in (\ref{5.5newa}), we add the analogue of this for $\sup_{v\in K}X_{1/2} (v)$ as an hypothesis.\qed

Similarly, under the conditions of  Corollary \ref{cor-LILMP} we can show that 
if  (\ref{az.2ka}) and (\ref{5.5newa}) hold for $X_{1/2}$ then (\ref{az.3k}) holds for $\wt L$ for   almost every $t$,   $  \wt P^{x}$ almost surely, for all $x\in S$.

\medskip At this point we   can not show that there are any 1/2--permanental processes   that are not symmetric  that satisfy  (\ref{az.2ka}).

    \section{Partial rebirthing of  transient Borel right processes } \label{sec-LILpr}

  In this section we apply Theorem \ref{theo-LILMP} and the LILs for the permanental processes in Theorems \ref{theo-1.2zz}--\ref{theo-1.5} to obtain LILs for the local times of   related Borel right processes. 

 \medskip    Let $T$ a be locally compact set with a countable base. 
  Let      $ Y = 
( \Om,   \FF_{t},   Y_t,  \th_{t},\newline   P^x
)$ be a transient Borel right process with state space $T$,  and  continuous strictly positive  potential   densities $u(x,y)$ with respect to some $\si$--finite measure $m$ on $T$.    Let $\mu$ be a finite  positive measure on $T$   with mass $|\mu|=\mu(T)$. 
We call the   function   
\begin{equation}
 f(x)= \int_{T} u(y,x)\,d\mu(y)\label{rp.1}
\end{equation}
 a left potential for $  Y$. %Since $u(x,y) $ is continuous in $y$ uniformly in $x$ and $\mu$ is a finite measure
 %we see that $f(y)  $ is continuous. See \cite[Section 2]{FR}.
 (It follows from \cite[Chapter 1, Proposition 14]{Bertoin} that for the   Markov processes $Y$ considered in Theorems   \ref{theo-1.2zz}, all integrable excessive functions are of this form.)
The following is \cite[Theorem 6.1]{MRejp}.

\bt\label{theo-borelpr}
Let   $ Y = 
( \Om,   \FF_{t},   Y_t,  \th_{t},    P^x
)$ be a transient Borel right process with state space $T$,   as above. Then for any left potential $  f$  for $   Y$, there exists a transient Borel right process  
 $\wt Z \!=\!
(\wt \Om, \wt \FF_{t}, \wt Z_t, \wt\th_{t}, \wt P^x)$ with state space $ S=T\cup \{\ast\}$, where $\ast$ is  an isolated point, such that  $\wt Z $ has potential densities, 
\begin{eqnarray}
&&  \wt    u(x,y)= u(x,y)+ f(y), \hspace{.2 in}x,y\in T
\label{rp.2}\\
&&   \wt    u(\ast,y)=   f(y), \hspace{.2 in}\mbox{ and }\hspace{.2 in} \wt    u(x, \ast)=\wt    u(\ast, \ast) =1, 
\nonumber
\end{eqnarray}
with respect to the measure $\wt m$ on $S$ which is equal to $m$ on $T$ and assigns 
a unit mass to $\ast$.  \et

 Note that in the proof of Theorem \ref{theo-borelN} we actually show that there exists a $1/2$--permanental process  $ X_{1/2}=\{X_{1/2}(y)\,,\,y\in S\}$ with kernel $ \wt    u(x,y)$.  
All the permanental processes considered in Theorems \ref{theo-1.2zz}--\ref{theo-1.5} are actually   defined on $S=T\cup  \{ *\}$ and    we can extend  the metric on $T$ to $S$ so that $*$ is an isolated point.

  \begin{remark} \label{rem-6.1}{\rm The process  $\wt Z $ is a partial rebirthing of  $   Y$.  When   $\wt Z $ starts in $T$, it proceeds just like $ Y$ until the lifetime of $  Y$. It then goes to the isolated point $*$. It remains in $*$ for an independent exponential time with mean,  
  \begin{equation} 
   \frac{1}{1+|\mu|}.  \end{equation} 
   Then it returns to $T$ with law $\mu/(1+|\mu|)$.     (The  probability that $\wt Z $ returns to $T$ is $|\mu|/(1+|\mu|)$.  With probability $1/(1+|\mu|)$  it goes to the cemetery state.)   We continue to repeat this procedure. At each stage the  process $\wt Z $ is ``reborn'', but   not with probability $1$. We call this a  partial rebirthing.  (This is explained in greater detail in the   proof of 
  \cite[Theorem 6.1]{MRejp}.) 
}
  \end{remark} 
  
  The analogues of $u(x,y)$, in (\ref{rp.2}), in the kernels of the permanental processes in Theorems \ref{theo-1.2zz}--\ref{theo-1.5} are symmetric. The next lemma uses this property.

\begin{lemma} \label{lem-6.1a} Let  $  Y$  be a strongly symmetric  transient Borel right process with transition densities $p_{s}(x,y)$ with respect to a $\si$--finite measure  $m$.  Then $  f=\{  f(x), x\in T \}$, in (\ref{rp.1}), is an excessive function for $ Y$. 
 %\textbf{  Do we need to say this. It is rather obvious and we don't use it now.}  In %particular, for any fixed $y$, $ u(x,y)$ is an excessive function for  $\wt Y$.
\end{lemma}

 \Proof We show that  
 \begin{equation} \label{6.7mm}
E^x\( f(  Y_{t})\)\uparrow   f(x),\quad \text{as $t\downarrow 0$,}\quad   \forall\, x\in T.  
\end{equation} 
The condition that $  Y$  is strongly symmetric  implies that $u(x,y)=u(y,x)$, for all $x,y\in T$. We also note that,
\begin{equation}
u(x,y)=\int_{0}^{\ff}p_{s}(x,y)\,ds.\label{td.1}
\end{equation}
Therefore, for any finite measure $\mu$, 
\begin{equation}
  f(x)= \int_{T} u(y,x)\,d\mu(y)= \int_{T} u(x,y)\,d\mu(y)=\int_{T} \int_{0}^{\ff}p_{s}(x,y)\,ds\,d\mu(y).\label{td.2}
\end{equation}
We have, 
\begin{eqnarray}
E^x\(  f(  Y_{t})\)&=& \int_{T} p_{t}(x,z)  f(z)\,dm(z)\label{td.3}\\
&=&\int_{T} \int_{0}^{\ff}\(  \int_{T} p_{t}(x,z)p_{s}(z,y)\,dm(z)\)\,ds \,\,d\mu(y)  \nonumber\\
&=&\int_{T} \int_{0}^{\ff}   p_{s+t}(x,y) \,ds \,\,d\mu(y) =\int_{T} \int_{t}^{\ff}   p_{s}(x,y) \,ds \,\,d\mu(y).\nonumber
\eea
Since,
\be
 \int_{T} \int_{t}^{\ff}   p_{s}(x,y) \,ds \,\,d\mu(y)\uparrow \int_{T} \int_{0}^{\ff}p_{s}(x,y)\,ds\,d\mu(y)=  f(x),\quad \mbox{as }t\downarrow 0,     
\ee
we get (\ref{6.7mm}).\qed

   The next theorem  is an  immediate application of Theorem \ref{theo-LILMP}. It gives LILs for the local times of the partially rebirthed Markov processes $\wt Z$ in Theorem \ref{theo-borelpr}. The spaces $S$ and $T$ are as defined in Theorem \ref{theo-borelpr}.

  \begin{theorem} \label{theo-6.2a} 
  Let  $\wt L=  \{ \wt L_{t}^{y}, (y,t)\in  S \times R_{+}^1\}$ be the local times of  $\wt  Z  $, normalized so that
\begin{equation}
   \wt E^{ x}\( \wt L^{y}_\ff\)=\wt u(x,y). \label{az.1u}
\end{equation}
%Define, 
%\begin{equation} \label{6.5nnq}
  %    \si^2(x,y)=   u(x,x)+   u(y,y)-2    u(x,y).
%\end{equation}\nc \textbf{ S not T !! }
If  a permanental process with kernel,
\begin{equation}
\wt    u_{ T}(x,y)= u(x,y)+ f(y), \hspace{.2 in}x,y\in T,\label{}
\end{equation}
satisfies  the conditions in any of Theorems \ref{theo-1.2zz}--\ref{theo-1.5}, 
then $ \wt L_{t}^{y}$ is jointly continuous on $S\times [0,\ze)$ and  for any  $d\in T$  for which the LIL for the permanental process holds,
\begin{equation}
   \limsup_{  x\to 0}\frac{| \wt L_{t}^{d+x}-\wt  L_{t}^{d}|}{  (  2       \si^2(d+x,d)\log\log 1/|x|)^{1/2}}=  \({2\wt L_{t}^{d}}\)^{1/2} \qquad a.s.  \mbox{ for a.e. } t<\ze, \label{pr.al}
   \end{equation}
     $\wt P^{v}$ almost surely,   for all $v\in T$, where 
\begin{equation} \label{6.5nn}
      \si^2(x,y)=   u(x,x)+   u(y,y)-2    u(x,y).
\end{equation} 
\end{theorem}
 (The relationship between $u$ and $\wt u$ is given in (\ref{rp.2}).)

\medskip\Proof 
This follows from Theorem \ref{theo-LILMP}.
 Since $\ast$ is an isolated point, we only need to take the $\limsup$ over  $x\in T$. 
 That $\wt L$ is jointly continuous on $S\times [0,\ze)$,  follows by \cite[Theorems 1.1 and 1.2]{MRsuf}. 
 The fact that the hypotheses on $u$ and $f$  imply  the conditions in these theorems is well known;   see  e.g., \cite[Chapter 6]{book}.

  (Note that the terms that take the place of   $   \si^2(d+x,d)$ in (\ref{1.34szqa}), (\ref{80.6}), (\ref{1.34szw}), (\ref{1.34szq}) and (\ref{diffa.3q}) are asymptotically  equivalent to $    \si^2(d+x,d)$ for the relevant potential densities $ u(x,y)$. This is shown in the proofs of Theorems  \ref{theo-1.2zz}--\ref{theo-1.5}.)
\qed
 
\begin{remark}{\rm  The measure $\mu$ that determines the partially rebirthed Markov process $\wt Z$ does not influence the denominator on left--hand side of (\ref{pr.al}). It is the same for all $\mu$ and consequently for all potentials $  f$ in (\ref{rp.2}). The measure $\mu$, through $  f$, influences  the local times $\wt L_{t}^{y}$,   as one can 
see in (\ref{az.1u}).
}\end{remark}

  There are six different functions   $u(x,y)$
in Theorems  \ref{theo-1.2zz}--\ref{theo-1.5} and nine sets of conditions on $T$ and $  f$. (Take $g\equiv 1$.) To clarify   Theorem \ref{theo-6.2a} we give  three specific examples.

 \begin{example} {\rm \label{ex-6.1}  Consider the permanental processes in Theorem \ref{theo-1.2zz}, I. and the corresponding partially rebirthed process  $\wt Z $ in Theorem \ref{theo-borelpr} for some finite positive measure $\mu$. The process $\wt Z $ is a rebirthing of a symmetric L\'evy process killed at the end of an exponential time with mean $1/\bb$ as described prior the the statement of Theorem \ref{theo-1.2zz}.    In this case     $T=R^1$,
$u(x,y)$
 in (\ref{rp.2}) is given by $u^\bb(x,y)=u^\bb(x-y)$ in (\ref{1.21nn}) and $(\si^\bb)^2(x)$ is given by (\ref{sigti.1}) and (\ref{lv.21}) for any $\bb>0$. The function $  f$ in 
 (\ref{rp.2}) is   a $C^{1}(R^1)$ excessive function  for   $Y$ in Theorem \ref{theo-1.2zz}.
 Note that the   requirement that $    f\in C^1(R^1)$ imposes restrictions on the measures $\mu$ that can be used for the rebirthing.   Similar restrictions  apply to all the applications of   Theorems \ref{theo-1.2zz}--\ref{theo-1.5} and Theorem \ref{theo-borelpr} in Theorem \ref{theo-6.2a}.}
 \end{example}
 
\medskip The following corollary is an application of Theorem \ref{theo-6.2a}. As in   Theorem \ref{theo-6.2a} the spaces $S$ and $T$ are as defined in Theorem \ref{theo-borelpr}. In this case $T=R^1$

\begin{corollary} \label{cor-6.1q} 
  Let  $\wt L=  \{ \wt L_{t}^{y}, (y,t)\in  S \times R_{+}^1\}$ be the local times of  $\wt Z $, normalized so that
\begin{equation}
   \wt E^{ x}\( \wt L^{y}_\ff\)= \wt u(x,y). \label{az.1uqqq}
\end{equation}
When $ f\in C^1(R^1)$ and $(\si^\bb)^2$    satisfies conditions $(i)$ and $(ii)$ in Theorem \ref{theo-1.2zz},  $\wt L$ is jointly continuous on $S\times [0,\ze)$, 
 and  for  all $d\in R^1$, 
\begin{equation}
   \limsup_{  x\to 0}\frac{|\wt  L_{t}^{d+x}- \wt L_{t}^{d}|}{  (  2    (\si^\bb)^2(x)\log\log 1/|x|)^{1/2}}=  \({2\wt L_{t}^{d}}\)^{1/2} \qquad a.s.  \mbox{ for a.e. } t<\ze, \label{pr.al3}
   \end{equation}
     $\wt P^{v}$ almost surely,  for all  $v\in R^1$.   
\end{corollary}

   \begin{example} {\rm 
Consider the permanental processes in Theorem \ref{theo-1.2aa}, I. and the corresponding partially rebirthed processes $\wt Z$ in Theorem \ref{theo-6.2a} with  respect to a finite   positive measure $\mu$ on  $T=R^{1}-\{0\}$.    The process $ \wt  Z $ is a rebirthing of a symmetric process  $\ov Y$ which is obtained by killing a    symmetric L\'evy  process    $ Y$   the first time it hits 0,  as described prior to the statement of Theorem \ref{theo-1.2aa}.      In this case  $u(x,y)$
 in (\ref{rp.2}) is given by $u^0(x,y)$ in (\ref{u0}) and $(\si^0)^2(x)$ is given by (\ref{1.44mm}).   The function $  f$ in 
 (\ref{rp.2}) is any $C^{1}(R^1)$ excessive function  for   $\ov Y$ in Theorem \ref{theo-1.2aa}.

 \medskip The following corollary is an application of Theorem \ref{theo-6.2a} in this case:

    \begin{corollary} \label{cor-6.2}
  Let  $\wt L=  \{\wt  L_{t}^{y}, (y,t)\in  S \times R_{+}^1\}$ be the local times of  $\wt  Z $, normalized so that
\begin{equation}
  \wt E^{ x}\( \wt L^{y}_\ff\)= \wt u(x,y). \label{az.1uqq}
\end{equation}
Then the local times of $\wt   Z  $ are jointly continuous on $S\times [0,\ze)$.
When $  f\in C^1(R^1-\{0\})$ and $(\si^0)^2$    satisfies conditions $(i), (ii)$ and $ (iii)$ in Theorem \ref{theo-1.2aa},
then for any   any $d\neq 0$,
\begin{equation}
   \limsup_{  x\to 0}\frac{| \wt L_{t}^{d+x}- \wt L_{t}^{d}|}{  (  2 (\si^0)^2(x)\log\log 1/|x|)^{1/2}}=  \({2 \wt L_{t}^{d}}\)^{1/2} \qquad a.s.  \mbox{ for a.e. } t<\ze, \label{pr.alsq}
   \end{equation}
     $\wt P^{v}$ almost surely, for all $v\in T$. 
Furthermore,\begin{equation}
   \limsup_{  x\to 0}\frac{| \wt L_{\ff}^{d+x}- \wt L_{\ff}^{d}|}{  (  2 (\si^0)^2(x)\log\log 1/|x|)^{1/2}}=  \({2\wt L_{\ff}^{d}}\)^{1/2} \qquad a.s.,    \label{pr.alsda}
   \end{equation}
 $\wt P^{v}$ almost surely, for all   $v\in   R^1-\{0 \}$.
 \end{corollary}
}\end{example}

 \Proof We only need to prove (\ref{pr.alsda}).   This follows  because in this case the Markov process  is killed   when it hits $0$. Since the paths are right continuous,   for any  $d\neq 0$, there is a  random interval  $I=[s,\ze)$ during which the path is not in some neighborhood, say B, of $d$. Consequently, for any $t\in I$   and any $d+x\in B$, we have that $\wt L_{t}^{d+x}=\wt L_{\ff}^{d+x}$,  so (\ref{pr.alsda}) follows from (\ref{pr.alsq}). \qed

\begin{example} {\rm 

   Consider the permanental processes   $Z'$ in Theorem \ref{theo-1.5}, {\bf b}. and the corresponding partially rebirthed process  $\wt   Z $ in Theorem \ref{theo-6.2a} with  respect to a finite   positive measure $\mu$   on     $T=(0,\ff)$.    The process   $Z' $ is a diffusion  killed at the end of an exponential time with mean $1/\bb$, killed the first time it hits 0, as described prior to the statement of Theorem \ref{theo-1.5}. In this case the generic potential density $    u(x,y)$ in (\ref{rp.2}) is  $\wt v^\bb(x,y)$ in (\ref{diff.22}) and $f$ is an excessive function for $Z'$ in (\ref{diff.22}).
  
   \medskip The following corollary is an application of Theorem \ref{theo-6.2a} in this case:

   \begin{corollary} \label{cor-6.1}
  Let  $\wt  L=  \{ \wt  L_{t}^{y}, (y,t)\in  S \times R_{+}^1\}$ be the local times of   $ \wt  Z $, normalized so that
\begin{equation}
  \wt  E^{ x}\( \wt  L^{y}_\ff\)= \wt u(x,y). \label{az.1uq}
\end{equation}
Then the local times of   $\wt   Z $ are jointly continuous on $S\times [0,\ze)$.
When    $f\in C^2(\De_{d}(\de))$ for some $\de>0$, then 
   for  any $d>0$  for which $ f'(d)=0$,  \begin{equation}
   \limsup_{  x\to 0}\frac{|\wt   L_{t}^{d+x}- \wt  L_{t}^{d}|}{  (  2 \tau(d)|x|\log\log 1/x)^{1/2}}=  \({2\wt  L_{t}^{d}}\)^{1/2} \qquad a.s.  \mbox{ for a.e. } t<\ze, \label{pr.alsm}
   \end{equation}
     $\wt  P^{v}$ almost surely, for all $v\in (0,\ff)$.
    
Furthermore,\begin{equation}
   \limsup_{  x\to 0}\frac{| \wt  L_{\ff}^{d+x}-  \wt L_{\ff}^{d}|}{  (  2 \tau(d)|x|\log\log 1/x)^{1/2}}=  \({2\wt  L_{\ff}^{d}}\)^{1/2} \qquad a.s.,    \label{pr.alsdm}
   \end{equation}
\end{corollary}
 $\wt  P^{v}$ almost surely, for all $v\in (0,\ff)$. 
}\end{example}

 \Proof The proof is the same as the proof of Corollary \ref{cor-6.2}.\qed
 
  \begin{example}\label{ex-6.4}  {\rm 

   Consider the permanental processes   $\ov Z$ in Theorem \ref{theo-1.5}, {\bf c}. and the corresponding partially rebirthed process  $\wt   Z $ in Theorem \ref{theo-6.2a} with  respect to a finite   positive measure $\mu$ on     $T=(0,\ff)$.    The process   $\ov Z  $ is a diffusion   killed the first time it hits 0, as described prior to the statement of Theorem \ref{theo-1.5}. In this case the generic potential density $    u(x,y)$ in (\ref{rp.2}) is  $\wt u_{T_0} (x,y)$ in (\ref{diff.5}) and $f$ is an excessive function for $\ov Z$ in (\ref{diffa.1}). (
   We take $g\equiv 1$.)
  
   \medskip The following corollary is an application of Theorem \ref{theo-6.2a} in this case:

   \begin{corollary} \label{cor-6.4}
  Let  $\wt  L=  \{ \wt  L_{t}^{y}, (y,t)\in  S \times R_{+}^1\}$ be the local times of   $ \wt  Z $, normalized so that
\begin{equation}
  \wt  E^{ x}\( \wt  L^{y}_\ff\)= \wt u_{T_0} (x,y)+f(y). \label{az.1uq}
\end{equation}
Then the local times of   $\wt   Z $ are jointly continuous on $S\times [0,\ze)$.
When   $f\in C^2(\De_{d}(\de))$ for some $\de>0$, then 
   for  any $d>0$  for which $ f'(d)=0$, 
\begin{equation}
   \limsup_{  x\to 0}\frac{|\wt   L_{t}^{d+x}- \wt  L_{t}^{d}|}{  (  2 s'(d)|x|\log\log 1/x)^{1/2}}=  \({2\wt  L_{t}^{d}}\)^{1/2} \qquad a.s.  \mbox{ for a.e. } t<\ze, \label{pr.alsmz}
   \end{equation}
     $\wt  P^{v}$ almost surely, for all $v\in (0,\ff)$.
    
Furthermore,\begin{equation}
   \limsup_{  x\to 0}\frac{| \wt  L_{\ff}^{d+x}-  \wt L_{\ff}^{d}|}{  (  2  s'(d)|x|\log\log 1/x)^{1/2}}=  \({2\wt  L_{\ff}^{d}}\)^{1/2} \qquad a.s.,    \label{pr.alsdmz}
   \end{equation}
\end{corollary}
 $\wt  P^{v}$ almost surely, for all $v\in (0,\ff)$. 
}\end{example}

 \Proof The proof is the same as the proof of Corollary \ref{cor-6.2}.\qed

 \begin{remark}{\rm \label{rem-6.3} The LILs in (\ref{1.66mm}) and (\ref{1.67mm}) are given in (\ref{pr.alsmz}) and (\ref{pr.alsdmz}) with, obviously, $d$ replaced by $x_0$. We just have to check that the conditions on $f$ are satisfied. We do this in Remark \ref{rem-9.1}. }\end{remark}

  \medskip There are   five more   sets 
 of conditions on  $u(x,y)$, $T$ and $  f$ in Theorems \ref{theo-1.2zz}--\ref{theo-1.5} for which specific realizations of (\ref{pr.al})
 can be written out. We leave this to the interested reader.

 \section{Rebirthing   transient Borel right processes }\label{sec-Borel}

  In this Section we apply Corollary   \ref{cor-LILMP} to recurrent Markov processes and  the LILs for permanental processes in Theorems \ref{theo-1.2zz}, \ref{theo-1.2zzk} and \ref{theo-1.5} to obtain LILs for the local times of  certain fully rebirthed    Borel right processes.

In Theorem \ref{theo-borel} we obtain a general result about the potential densities of   fully rebirthed Borel right processes: 
    Let $S$ a be locally compact space with a countable base.  
Let   $\YY=(\Om,  \FF_{t},  \YY_t,\th_{t}, P^x)$ be a transient Borel right process with state space $S$  and  continuous strictly positive  $p-$potential   densities  
\be
\ov u^{p}=\{\ov u^{p}(x,y), x,y\in S \},
\ee
 with respect to some $\si$--finite  positive measure $m$ on $S$. Let $\ze=\inf\{t\,|\,\YY_t=\De \}$,  where $\De$ is the cemetery state for   $\YY$  and assume       that $\ze<\ff$   almost surely.    We do not require that   $\YY$
is symmetric.

  Let $\mu$ be a probability measure on $S$. We       modify $\YY$ so that instead of going to $\De$ at the end of its lifetime it is immediately ``reborn'' with  measure $\mu $.  (I.e., the process goes to the set $B\subset S$ with probability $\mu(B) $, after which  it continues to evolve the way $\YY$ did, being reborn with probability   $\mu $ each time it dies.)   We denote this rebirthed process by $\wt Z\!=\!
(\wt \Om, \wt  \FF_{t}, \wt  Z_t, \wt \th_{t},\wt P^x)$.   To emphasize  the parameters that determine $\wt Z$ we sometimes write $\wt   Z=\wt Z(\ov u^{p}, m ,\mu)$.

 \bt\label{theo-borel} 
The process  
 $\wt   Z=\wt Z(\ov u^{p}, m ,\mu)$ is a recurrent Borel right process with state space $ S$ and $p$--potential densities, 
\begin{equation}
w^{p}(x,y)= \ov u^{p}(x,y)+\(\frac{1}{p}-\int_S \ov u^{p}(x,y)\,dm(y)\)\frac{f(y)}{\|f\|_1}\label{193.3a},
\end{equation}
with respect to   $m$, where  
\begin{equation} \label{7.9a}
  f(y)=\int_S \ov u^{p}(x,y)\,d\mu(x),
\end{equation}
and the $L_1$ norm is taken with respect to $m$.
\et

  \noindent \textbf{Proof of Theorem \ref{theo-borel}}       It is a classical result that   $\wt Z$  is a recurrent Borel right process. (See, e.g., \cite[Theorem 1]{Meyer} or \cite[(14.17)]{S}.) We continue to show that it has the potential density in (\ref{193.3a}).  
 
    We use the following notation.   For measurable functions  $h$ on $S$ we  define  the $p$--th potential operator,
\begin{equation} \label{6.1mm}
  \ov U^{p}h(x)=\int_S\ov   u^{p}(x,y)h(y)\,dm(y). 
\end{equation}  
As a probabilistic statement this is 
\begin{equation} \label{6.2mm}
  \ov U^{p}h(x)=  E^{x}\( \int_{0}^{\ff}e^{-pt}h\(\YY_{t}\) \,dt \)=  E^{x}\( \int_{0}^{\ze}e^{-pt}h\( \YY_{t}\) \,dt \),
\end{equation}
as in (\ref{potdef7}).
 For positive measures $\mu$  on $S$ we define,  
  \begin{equation} 
  \mu h =\int_S   h(x)\,d\mu(x). 
\end{equation}
%and note that,  
%\begin{equation} 
%  \mu \ov u^{p}(\cdot,y)= \int \ov  u^{p}(x,y)\,d\mu(x):=\mu \ov U^p(y).
%\end{equation} 
  Note that by (\ref{6.1mm}),
\begin{equation} \label{6.5mm}
  \mu \ov U^ph=\int_S \int_S \ov u^{p}(x,y)h(y)\,dm(y)\,d\mu(x).
\end{equation}
In particular,     when $h\equiv 1 $,  (\ref{6.1mm}) is
\begin{equation}
\ov U^p1(x)= \int_S \ov u^{p}(x,y)\,dm(y),\label{6.5mml}
\end{equation}
and (\ref{6.5mm}) is
\begin{equation} \label{6.5mmma}
  \mu \ov U^p1=\int_S \int_S \ov u^{p}(x,y)\,dm(y)\,d\mu(x).
\end{equation}

   We now take $h$ be bounded and consider  $W^p$,   the $p$--th potential operator for $\wt Z$.
 For any $x\in S$,  \bea \label{6.6mk}
 && W^{p}h(x)= \wt E^{x}\( \int_{0}^{\ff} e^{-pt}h\(\wt Z_{t}\) \,dt \) \\
  &&\quad=\nn  \wt E^{x}\(\int_{0}^{\ze_1}e^{-pt}h\(\wt Z_{t}\)\,dt\)+ \sum_{n=1}^{\ff}\wt E^{x}\(\int_{\wt \ze_n}^{\wt \ze_{n+1}}e^{-pt}h\(\wt Z_{t}\)\,dt\)\\
  &&\quad=\nn
   \wt  E^{x}\(\int_{0}^{\ze }e^{-pt}h\(   \wt Z_{t}\)\,dt\)+ \sum_{n=1}^{\ff}\wt E^{x}\(e^{-p\wt \ze_n} \int_{0}^{ \ze_{n+1}}e^{-pt}h\(\wt Z_{t}\circ \wt\th_{\wt \ze_n}\)\,dt\).
\eea
 where $\wt \ze_n=\sum_{j=1}^{n}\ze_j$, and $\{\ze_j \}_{j=1}^{\ff}$ are the consecutive lifetimes of the rebirthed process, so that $ \ze_{n+1}= \ze\circ \wt\th_{\wt \ze_n}$. 
 
 Using the fact that  $Z$ is a Markov process we have,
\bea \label{7.12}
&&\wt E^{x}\(e^{-p\wt \ze_n} \int_{0}^{ \ze_{n+1}}e^{-pt}h\(Z_{t}\circ\wt\th_{\wt \ze_n}\)\,dt\)\\
&&\qquad =\wt E^{x}\(e^{-p\wt \ze_n}\( \int_{0}^{ \ze }e^{-pt}h\(\YY_{t}\)\,dt\)\circ\wt\th_{\wt \ze_n}  \)\nn\\
&&
\qquad=\wt E^{x}\(e^{-p\wt \ze_n}\wt E^{\wt Z_{\wt \ze_n}}\( \int_{0}^{ \ze }e^{-pt}h\(\YY_{t}\)\,dt\) \)\nn\\
&&\qquad
=\wt E^{x}\(e^{-p\wt \ze_n} \)\wt E^{\mu}\( \int_{0}^{ \ze }e^{-pt}h\(\YY_{t}\)\,dt\),\nn
\eea
where we use  the fact that $\wt Z_{\wt \ze_n}$ has   probability distribution $\mu  $ and is independent of $\wt\FF_{\wt \ze_n}$.    We also use the notation that for any random variable $\cal Z$ on $\wt \Om$,  
\begin{equation} 
 \wt  E^{\mu}\(\cal Z\)=\int \wt  E^{y}\( \cal Z\)  \,d\mu(y).
\end{equation}
  Note that,
\begin{equation} \label{7.11mm}  \wt  E^{x}\(\int_{0}^{\ze }e^{-pt}h\(   \wt Z_{t}\)\,dt\)= E^{x}\(\int_{0}^{\ze }e^{-pt}h\(  {\cal Y}_{t}\)\,dt\), 
\end{equation}
 and
 \begin{equation} \label{7.14mm}  \wt  E^{\mu}\(\int_{0}^{\ze }e^{-pt}h\(   \wt Z_{t}\)\,dt\)= E^{\mu}\(\int_{0}^{\ze }e^{-pt}h\(  {\cal Y}_{t}\)\,dt\). 
\end{equation}

    Using the Markov property and  induction, in a similar manner  we see  that for $n\ge 1$,  \be  \label{7.14qq}
 \wt   E^x(e^{-p\wt \ze_n}) = E^{x}\( e^{-p\ze }\)\wt   E^{\mu}\( e^{-p\wt \ze_{n-1}}\)
  = E^{x}\( e^{-p\ze }\) \(  E^{\mu}\( e^{-p \ze }\)  \)^{n-1}.  
\ee 
  Note that  by simple integration and (\ref{6.2mm}),\begin{equation}
E^{x}\( e^{-p\ze}\)=1-pE^{x}\( \int_{0}^{\ze}e^{-pt} \,dt \)=1-p\ov U^{p}1(x),  \label{rp.4}
\end{equation}
and   
\begin{equation} 
  E^{\mu}\( e^{-p\ze}\)= 1-p\mu \ov U^{p}1,  \label{rp.4a}
\end{equation}
which implies, in particular, that $p\mu \ov U^{p}1<1$.  
 Using this  (\ref{7.12}),   (\ref{7.11mm}) and (\ref{7.14mm})  in (\ref{6.6mk})
we have, 
\begin{eqnarray}
W^{p}h(x)&=&E^{x}\(\int_{0}^{\ze}e^{-pt}h\(\YY_{t}\)\,dt\) \\
&&\qquad+E^{x}\( e^{-p\ze}\)   \sum_{n=1}^{\ff}\(    E^{\mu}\(e^{-p\ze}\)  \)^{n-1}E^{\mu}\(\int_{0}^{\ze}e^{-pt}h\(\YY_{t}\)\,dt\)
\nn \\
&=& \ov  U^{p}h(x)+\frac{E^{x}\( e^{-p\ze}\)}{1-E^{\mu}\(e^{-p\ze}\)}\mu\ov  U^{p}h,\nn
\eea
where we use the notation in (\ref{6.1mm})--(\ref{6.5mm}).
Consequently,
\begin{equation}
W^{p}h(x)=\ov U^{p}h(x)+\(1-p\ov U^{p}1(x) \)\frac{\mu \ov U^{p}h}{p\mu \ov U^{p}1}. \label{rp.5}
\end{equation}
Since this is true for all $h$  we see that (\ref{193.3a}) follows by  using (\ref{6.1mm}), (\ref{6.5mml}) and   (\ref{6.5mmma}).
 \qed

        \begin{lemma} \label{lem-7.1}  The function 
        \begin{equation} 
  g(x)=1-p\int_S \ov u^{p}(x,y)\,dm(y)
\end{equation}   is an excessive function   for  the   process  
 obtained by killing  $\YY$  at the end of an independent exponential time with mean $1/p$. 
 \end{lemma} 
 
  \Proof    It follows from  (\ref{6.5mml}) and (\ref{rp.4})   that  $g\geq 0$.  To complete the proof we show that,
\begin{equation}
e^{-pt}E^{x}\(g\(\YY_{t} \)  \)\uparrow g(x), \qquad \mbox{as }t\downarrow 0.\label{pexcc.1}
\end{equation}
  It follows from (\ref{6.5mml}) and  (\ref{rp.4}) that to show this  it suffices to show that $E^{x}\( e^{-p\ze}\)$ is   a $p$--excessive function for $\YY$. I.e., that,
\begin{equation}
e^{-pt}E^{x}\(E^{\YY_{t}}\( e^{-p\ze}\)   \)\uparrow E^{x}\( e^{-p\ze}\), \quad \mbox{as }t\downarrow 0.\label{pexcc.2}
\end{equation}
  For $t<\ze$ we have $\ze\circ \th_{t}=\ze-t$; see \cite[(2.137)]{book}. Consequently, 
\begin{eqnarray}
E^{x}\( e^{-p\ze},t<\ze \)&=&e^{-pt}    E^{x}\( e^{-p\ze\circ \th_{t}},t<\ze \)
\nn\\
&=&  e^{-pt}E^{x}\(E^{\YY_{t}}\( e^{-p\ze}\)\) \nonumber
\end{eqnarray}
 by the strong Markov property. This gives(\ref{pexcc.2}). \qed

 \medskip  We   now use Theorem \ref{theo-borel} along with    Theorem \ref{cor-LILMP} and Theorems \ref{theo-1.2zz}, \ref{theo-1.2zzk} and  \ref{theo-1.5} to obtain LILs for the local times $\wt  L=  \{ \wt  L_{t}^{y}, (y,t)\in  S \times (0,\ff)\}$ of several fully rebirthed Markov processes. In particular whenever $\{\ov u^p(x,y),x,y \in T \}$ is any of the symmetric potentials in Theorems \ref{theo-1.2zz}, \ref{theo-1.2zzk} or  \ref{theo-1.5} and it and $f$ and $g$ satisfy the conditions in the relevant theorem we have,
    \begin{equation}
  \limsup_{  x \to 0}\frac{|  \wt   L_t^{ d+x}-  \wt   L_t^{ d }|}{ \(2   (\ov \si^p)^{2}\(d+x,d\)\log\log 1/|x|\)^{1/2}}=    \( \wt   L_t^{ d }\)^{1/2}, \qquad \text{ for a.e. t},\label{7.22mm}
  \end{equation} 
$ \wt  P^v$ almost surely, for all $d$ for which the relevant result in Theorems \ref{theo-1.2zz}, \ref{theo-1.2zzk} or  \ref{theo-1.5} hold, where,
\begin{equation} 
  (\ov \si^p)^{2}\(d+x,d\)=\ov u^p(d+x,d+x)+\ov u^p(d ,d )-2\ov u^p(d+x,d ).
  \end{equation}

  We give   two broad examples:

 \begin{example} {\rm\label{ex-7.1} Let $Y=\{Y_t,t\in R_{+}^1 \}$ to be a strongly symmetric Borel right process with locally compact state space $S$ and  $\ga$--potential densities $u^{\ga}(x,y)$ with respect to some $\si$--finite measure $m$ on $S$. Assume that $Y$ has a  Markov transition semigroup  $P_{t}$, so that in particular,  $P_{t}(x, S)=1$ for all $t$ and all $x\in S$. Note that for all   $\ga>0$,
\begin{equation}
\int_S u^{\ga}(x,y)\,dm(y)=\int_{0}^{\ff}e^{-\ga t}P_{t}(x, S)\,dt=\frac{1}{\ga}.\label{7.100q}
\end{equation}
 }\end{example}

 Let $\ov Y=\{\ov Y_t;t\in R^1_{+} \}$ be the    Markov  process  
 obtained by killing  $Y$ at the end of an independent exponential time with mean $1/\al$. ($\ov Y$ is an example of $\cal Y$ introduced in the beginning of this section.) Let $\ov u^{p}=\{\ov u^{p},x,y\in S \}$ denote  the  $p$--th potential density of   $\ov Y$, for  $p>0$. Note that $\ov u^{p}$ is also the potential density of $Y$ killed at the end of time $\xi\wedge \rho$, where $\xi$ and   $\rho$ are independent exponential random variables with means $1/\al$ and $1/p$.  It is well known that $\xi\wedge \rho$ is an exponential random variable with mean $1/(\al+p)$;  see e.g., \cite[Remark 5.2.8]{book}. Consequently,  
 \begin{equation}
\ov u^{p}(x,y)  =  u^{\al+p}(x,y) .\label{7.101a}
\end{equation}
We have
\begin{equation} 
   \int_{S }\ov u^{p}(x,y)\,dm(y)= \int_{S }u^{\al+p}(x,y)\,dm(y)=\frac{1}{\al+p}
\end{equation}
and  since $\mu$ is a probability measure
\begin{equation} 
  \|f\|_1=\int_{S }u^{\al+p}(x,y)\,dm(y)\,d\mu(x)=\frac{1}{\al+p}.
\end{equation}
 Using  these we can write  (\ref{193.3a}) as, \begin{equation} 
  w^{p}(x,y)=   u^{\al+p}(x,y)+\frac{\al}{p }f(y). 
\end{equation}
 
 It follows from Lemma \ref{lem-6.1a} that  $\{ w^{p}(x,y),  x,y \in S\}$ has the form of the kernels of the permanental processes considered in Theorems \ref{theo-1.2zz}
and \ref{theo-1.5} a.      Therefore, by Theorem \ref{theo-LILMP} whenever these permanental processes satisfy an LIL as in (\ref{1.34szqa}) or (\ref{1.34szq}),
(\ref{az.3}) holds for the local times of the fully rebirthed Markov processes $\wt Z(u ^{\al+p},m,\mu)$, where the local times $\wt  L=  \{ \wt  L_{t}^{y}, (y,t)\in  S \times (0,\ff)\}$ are normalized so that
\begin{equation}
  \wt E^{ x}\(\int_{0}^{\ff}  e^{-ps}\, d_{s}\wt L^{y}_s\)=   w^p(x,y). \label{az.1vw}
\end{equation}

 Theorems \ref{theo-1.2zz}
and \ref{theo-1.5} a. give conditions for 
when these permanental processes satisfy the LILs  in (\ref{1.34szqa}) or (\ref{1.34szq}). These conditions also imply that  the 1/2 permanental processes with  kernels $w^p(x,y)$ are continuous almost surely  and that $\{\wt L_s^y, s,y\in R^1\times R_+^1\}$ are jointly continuous almost surely.
These properties follow from \cite[Theorems 1.1 and 1.2]{MRsuf}. (The fact that the hypotheses on $u ^{\al+p}$ and $f$  imply  the conditions in these Theorems is well known;   see  e.g., \cite[Chapter 6]{book}.) Consequently   when Theorems \ref{theo-1.2zz}
and \ref{theo-1.5} a. satisfy the LILs  in (\ref{1.34szqa}) or (\ref{1.34szq}) it follows from Corollary \ref{cor-LILMP}.
that the local times of the fully rebirthed Markov processes $\wt Z(u^{\al+p},m,\mu)$,  normalized so that (\ref{az.1vw}) holds,
satisfy the LIL in  (\ref{az.3})   for almost every  $t\in R_+$,   $ \wt  P^{x}$ almost surely, for all $x\in S$.

    \begin{example}{\rm   \label{ex-7.2}  Let $Y$ and $\ov Y$
  be as in Example \ref{ex-7.1} so that, as in (\ref{7.101a}), the  $p$--th potential density of $\ov Y$ is
$u^{\al+p}(x,y)$. However, in   that  example we  take $Y$ and $\ov Y$ to have to  state space $R^{1}$. Now, let $\ov Y'=\{\ov Y'_t, t\in R^1-\{0\} \}$ be a  symmetric  Borel right  process with state space   $S=R^1-\{0\}$ that is obtained by killing $\ov Y$  the first time it hits $0$. It   follows from  \cite[(4.165)]{book} that the  $p$--th potential density of   $\ov Y'$ is  
\begin{equation}
\ov v^{p}(x,y)=v^{\al+p}(x,y)    =   u^{\al+p}(x,y)-\frac{u^{\al+p}(x,0)u^{\al+p}(0,y)}{u^{\al+p}(0,0)}.\label{7.107a}
\end{equation} 
We now use  (\ref{7.100q})  and (\ref{7.101a}) to see that, 
\be\label{7.102mq}
\int_{R^{1}} \ov v^{p}(x,y)\,dm(y)= \frac{1}{\al+p}\(1-\frac{u^{\al+p}(x,0)}{u^{\al+p}(0,0)} \) .\ee
Therefore, by (\ref{193.3a}), we have,  
\begin{equation}
  w^{p}(x,y)=   v^{\al+p}(x,y)+\frac{1}{\al+p}\(\frac{\al}{p}+\frac{u^{\al+p}(x,0)}{u^{\al+p}(0,0)}    \) \frac{ f(y)}{\|f\|_1},\label{7.24x}
\end{equation} 
where
\begin{equation} \label{7.104mb}
  f(y)=\int_S  v^{\al+p}(x,y)\,d\mu(x).
\end{equation} 
(It is interesting to note that,
\begin{equation} 
  \frac{1}{ p(1+p/\al)} <\frac{1}{\al+p}\(\frac{\al}{p}+\frac{u^{\al+p}(x,0)}{u^{\al+p}(0,0)}    \)\le \frac{1}{\al+p}\(\frac{\al}{p}+1  \)=\frac{1}{p}.
\end{equation}
This, and the fact that $f/\|f\|_1$ has $L^1$ norm 1, shows that the excessive functions added to $v^{\al+p}$ in (\ref{7.24x})   are not be too large.)

\medskip 
 Let $\ov Y''=\{\ov Y''_t,t\in R^1_{+} \}$ be the   process  
 obtained by killing  $\ov Y'$ at the end of an independent exponential time with mean $1/p$. Consequently, $\ov Y''$ has $0$--potential densities $v^{\al+p}(x,y)$.   Since $\{v^{\al+p}(x,y),x,y\in S \}$ is symmetric it follows from Lemma \ref{lem-6.1a} that  $f$ is an excessive function for $\ov Y''$   and from Lemma \ref{lem-7.1} that $g$  is an excessive function for $\ov Y''$.    \ Consequently, $\{\ov w^{p}(x,y),x,y \in S\}$ has the form of the kernels of the permanental processes considered in  Theorem \ref{theo-1.2zzk}.  
($\ov Y''$ is an example of $\cal Y$ introduced in the beginning of this section.)

 We are now at the same place as in the paragraph containing (\ref{az.1vw}) but with respect to the kernels of the permanental processes in   Theorem \ref{theo-1.2zzk}. 
  Therefore, by Theorem \ref{theo-LILMP}, whenever these permanental processes satisfy an LIL as in (\ref{80.6}) or (\ref{1.34szq}),
(\ref{az.3}) holds for the local times $\wt  L=  \{ \wt  L_{t}^{y}, (y,t)\in  S \times (0,\ff)\}$ of the fully rebirthed Markov processes $\wt  Z(v^{\al+p},m,\mu)$, where the local times are normalized so that
\begin{equation}
  \wt   E^{ x}\(\int_{0}^{\ff}  e^{-ps}\, d_{s}\wt  L^{y}_s\)=  \ov w^p(x,y). \label{az.1vwx}
\end{equation}

 Now we repeat the paragraph following (\ref{az.1vw}) with Theorems \ref{theo-1.2zz} replaced by   Theorem  \ref{theo-1.2zzk}   to see that when they satisfy the LILs in (\ref{80.6}) or (\ref{1.34szq})   the local times of the fully rebirthed Markov processes $\wt Z(v^{} ,m,\mu)$,  normalized so that (\ref{az.1vwx}) holds,
satisfy the LIL in  (\ref{az.3})   for almost every  $t\in R_+$,   $ \wt  P^{x}$ almost surely, for all $x\in S$. 
 }\end{example}

  \section{L\'evy processes   that  satisfy the hypotheses of Theorems  \ref{theo-1.2zz}--\ref{theo-1.2aa}}\label{sec-ex}

 In    Theorems \ref{theo-1.2zz}--\ref{theo-1.2aa} we require that
\begin{equation}
 (\si^{\bb})^2\(x\)=\frac{2}{\pi}\int_{0}^{\ff}\frac{1- \cos\la x}{\bb+ \psi(\la)}\,d\la,\qquad   \bb\ge 0,\label{lv.21q}
\end{equation}
   satisfies various conditions. With some regularity conditions on $\psi$ it does satisfy   them.   
We show in  \cite[Lemma 7.3.1]{book}
that  when $\psi$ is regularly varying at infinity with index $1<p\le  2$,
 \begin{equation} \label{1.26ee0}
  (\si ^{ \bb })^2(h)\sim C_p \frac{1}{|h|\psi(1/h)},\qquad  \text{as $h\to 0$ for all   $  \bb\ge 0$}.
\end{equation}
  where, 
  \be \label{1.26z}
  C_p=\frac{4}{\pi}\int_{0}^{\ff}  \frac{\sin^2 s/2}{s^p}\,ds. 
\end{equation}
 
In (\ref{1.32nn}) we point out that
   \begin{equation} \label{7.3n}
  \psi(\la)=C\la^2+\psi_1(\la),\qquad\text{where $\psi_1(\la)=o(\la^2)$ as $\la\to\ff$}.
\end{equation}
Clearly, when $C>0$, $\psi(\la)$ is regularly varying at infinity with index 2 and
  \begin{equation} \label{1.26ee}
  (\si ^{ \bb })^2(h)\sim   \frac{|h|}{C },\qquad  \text{as $h\to 0$ for all   $  \bb\ge 0$},
\end{equation}   
as we point out in (\ref{1.29nn}). 

  \medskip Suppose $C=0$.
 We now show that we can obtain functions $\psi_1(\la)$ with nice regularity conditions. We can represent $\psi_1(\la)$ by the equation,
\begin{equation} \label{4.1z}
  \psi_{1}(\la)=2\int_{0}^{\ff}(1-\cos\la x)\nu(dx),
\end{equation}
where $\nu$ is a   L\'evy measure on $(0,\ff)$,   i.e.,
\begin{equation} \label{4.2z}
  \int_{0}^{\ff}(1\wedge x^2)\nu(dx)<\ff.
\end{equation}  
We consider L\'evy measures $\nu$ that can be represented by the equation,  
  \begin{equation} \label{4.3z}
  \nu([u,\ff))=\int_{u}^{\ff}\frac{1}{\tau(x)}\,dx,
\end{equation}
  where $\tau$ is a regularly varying function at 0 with index $2<p\le 3$.  Note that the integral in (\ref{4.2z}) is finite at 0 when $\tau$ is a regularly varying function at 0 with index $2<p< 3$  and in certain cases when $p=3$.    
  
   We now write  (\ref{4.1z}) as 
\begin{equation} \label{5.6zz}     \psi_1(\la)=4\int_0^{\ff} \frac{\sin^2\la  u/2}{ \tau(u)}\,du.
\end{equation} 
It follows from \cite[Lemma 7.4.10]{book}
that when $\tau$ is a regularly varying function at 0 with index $2<p< 3$,
\begin{equation} \label{4.6zz}
   \psi_1(\la)\sim \frac{\pi C_p}{\la\tau(1/\la)},\qquad\text{as $\la\to\ff$},
\end{equation}
where $C_p$ is given in (\ref{1.26mmpa}).     This   shows   that $\psi$   can be taken to be   asymptotic to any  regularly varying function at infinity with index $1<p<2$, and   consequently,  that for every regularly varying function at zero  of index $0<q<1$,   say $g_q(x)$ it follows from    (\ref{1.26ee0}) that  there is a L\'evy process for which the functions 
\be
(\si^{ \bb } )^2(x)\sim   g_q(x),\quad  \text{as $x\to 0$ for all   $  \bb\ge 0$ }. \label{5.1mm}
\ee

When $\tau$ is a regularly varying function at 0 with index $3$, 
  \begin{equation} \label{4.6zzz}
   \psi_{1}(\la)\sim  2\la^2\int_0^{1/\la}\frac{u^2}{\tau(u)}\,du,\qquad\text{as $\la\to\ff$}.
\end{equation}
  This shows that  $\psi_1$  can also be   asymptotic to certain regularly varying functions at infinity with index 2.   See Remark \ref{rem-8.3nn}.  
(The integral in (\ref{4.6zzz}) is  \cite[7.284]{book} which contains the typo that  the term $\la^2$ is omitted.) 

  When the integral in     (\ref{4.6zzz}) is a slowly varying function     as $\la\to\ff$
we can write,
  \begin{equation} \label{z}
   \psi_1(1/h)\sim   \frac{2}{h^2} \int_0^{h}\frac{u^2}{\tau(u)}\,du:=\frac{1}{h^2L(1/h)},\qquad\text{as $h\to0$ for all   $  \bb\ge 0$}.
\end{equation}
where 
  $L(1/h)$ is a slowly varying function at 0  such that  $\lim_{h\to 0}L(1/h)=\ff$. (See e.g.  \cite[Theorem 14.7.4]{book}.) Therefore, by  (\ref{1.26ee0}),   \begin{equation} \label{1.26mmws}
  (\si ^{ \bb })^2(h)\sim C_2  {h}{L(1/h) },\quad  \text{as $h\to 0$ for all   $  \bb\ge 0$}.
\end{equation} 
In particular   if $\tau(u)=u^3(\log1/u)^{q+1}$, for any $q>0$,  
 \begin{equation} \label{1.26mmpb}
  (\si ^{ \bb })^2(h)\sim C_q  h(\log 1/h)^{q},\qquad  \text{as $h\to 0$ for all   $  \bb\ge 0$},
\end{equation}
 for constants $C_q$ depending on $q$.

\medskip We now consider specifically, the hypotheses on the functions $ (\si ^{ \bb })^2$ in each of the Theorems \ref{theo-1.2zz}--\ref{theo-1.2aa}. We consider the conditions on the excessive functions   in these theorems in Section \ref{sec-ex1}.

\subsection{ Theorem  \ref{theo-1.2zz} }The material above shows that there are many examples of the functions $ (\si ^{ \bb })^2$, $0<\bb\le 1$ that are regularly varying at zero and  satisfy the hypotheses of this  theorem. Furthermore, we see from  (\ref{5.1mm}) and (\ref{1.26mmpb}) that  (\ref{1.38nna})  holds for all  such functions with $0<\bb<1$  and  for all $p>1$ when $\bb=1$. In addition (\ref{1.26mmpb}) 
  provides examples in which Theorem  \ref{theo-1.2zz}, II   holds for functions $ (\si ^{ \bb })^2$, $\bb>0$, for L\'evy processes   without a Gaussian component, for which 
(\ref{1.38nna}) does not hold. (For example, as pointed out   following the statement of Theorem  \ref{theo-1.2zz},  part II applies to all L\'evy processes   with a Gaussian component. For such processes (\ref{1.38nna}) does not hold.)

  \subsection{Theorem  \ref{theo-1.2zzk}  }\label{sec-ex2}
 
  In this theorem   there is a condition that must be satisfied that is not required in Theorem  \ref{theo-1.2zz}, which is that for all $\bb\ge 0$, 
$ (\si ^{ \bb })^2\in C^{1}([d,d+\de])$   for   $d>0.$ With an additional regularity condition on $\psi$ we can show that this is the case.
  We consider L\'evy processes $Y$ with L\'evy exponents $\psi$ such that,\begin{equation} \label{4.13}
 \psi\(\la\)=\int_{1}^{2}  \la ^{s}\,d\mu(s),
\end{equation}
where $\mu$ is a finite positive measure on $(1,2]$.  We refer to these L\'evy processes as stable mixtures. (See Remark \ref{rem-sm}.)

Stable mixtures are a good  class of L\'evy processes to work with because we can get good  estimates   for the  derivatives of  $\psi(\la)$.

\medskip  The next  two lemmas follow  immediately from the representation of $\psi$ in (\ref{4.13}).  
\begin{lemma} \label{lem-4.1}Suppose that $\mu$ is supported on   $[\ga_0,\ga_1]$ where $1<\ga_0\le \ga_{1}\le 2$ and set  $\mu([\ga_0,\ga_{1}])=|\mu|$. Then   $\psi(\la)$ in (\ref{4.13}) satisfies,
  \begin{equation} \label{4.22qq}
  \psi(\la)\le \left\{
  \begin{matrix} 
  |\mu||\la|^{\ga_0},& \quad |\la|\le 1 \\
  \\
  |\mu||\la|^{\ga_{1}} ,& \quad |\la|\ge 1   
\end{matrix}  \right.
\end{equation}
and,
\begin{equation} \label{4.23qq}
  \psi(\la)\ge \left\{
  \begin{matrix} 
   \mu( [\ga_0,\ga_0+\ep])|\la|^{\ga_0+\ep},& \quad |\la|\le 1 \\
  \\
   \mu( [\ga_{1}-\ep,\ga_{1}] )|\la|^{\ga_{1} -\ep},& \quad |\la|\ge 1,   
\end{matrix}  \right.
\end{equation}
for all $\ep>0$ sufficiently small.  (Without loss of generality we can assume that 
$\ga_0$ and $\ga_1$   are such that for all $\ep>0$   both $ \mu( [\ga_0,\ga_0+\ep])> 0$ and  $\mu( [\ga_{1}-\ep,\ga_{1}] )>0$.)
\end{lemma}

 \begin{lemma}  Suppose that $\mu$ is supported on $[\ga_0,\ga_{1}]$ where $1<\ga_0\le \ga_{1}\le 2$.   Then $\psi$  in (\ref{4.13}) is in  $C^{2}\( (0,\ff)\) $ and for all  $\la>0$,
  \bea \label{4.14qq} 
\ga_0  \,\psi(\la)&\le  & \la\psi'(\la)\le \ga_{1}  \,\psi(\la) \\
\ga_{0} (\ga_{0}-1) \,\psi(\la)&\le  & \la^2\psi''(\la)\le \ga_{1} (\ga_{1}-1) \,\psi(\la).\nn
\eea

\end{lemma}

\begin{lemma} \label{lem-7.3}
  When $\psi$ is the   characteristic exponent of a L\'evy   process  that is a stable mixture, for all   $\bb\ge 0$,    $(\si^{\bb})^{2}(x)\in C^{2}\(R^1- \{0 \}\)$ and,  
\be  \label{4.58aqe}
 |(\si^{\bb}) ^2 (x)'|\le \frac{  (\si^{\bb}) ^2 (x)}{|x|} ,
 \ee 
   where $(\si^{\bb}) ^2 (x)':=\displaystyle \frac{d}{dx}(\si^{\bb}) ^2 (x)$.
\end{lemma}

 \medskip   Differentiating under the integral sign is critical in  the proof of this lemma. Because we can not find a suitable reference  we include the following lemma: 

\bl\label{lem-8.4}
Let $f(t,x)$ be a function on $(0,\ff)\times (a,b)$  such that for each $t\in (0,\ff)$,  $f_{x}(t,x):= {d} f(t,x)) /{dx}$
 is continuous on $(a,b)$.  Assume furthermore that $\sup_{x\in (a,b)}|f_x(\cdot,x)|\in L^{1}\((0,\ff), \,d\mu \)  $    for some  measure $\mu$ on $(0,\ff)$. Then 
\be  \label{9.7w}
  \int_{0}^{\ff} f_{x} (t,x)\,d\mu (t) \text{ is continuous on $ (a,b)$, } 
 \ee
% \be\label{9.7q} \int_{0}^{\ff} f (t,x)\,d\mu (t)    \text{ is differentiable on $  (a,b),$ }
%  \ee 
and  for all $x\in (a,b)$,
\be \frac{d}{dx}\int_{0}^{\ff} f (t,x)\,d\mu (t)=\int_{0}^{\ff} f_{x} (t,x)\,d\mu (t).\label{9.8w}
\end{equation}
\end{lemma} 

  \Proof   The statement  in (\ref{9.7w}) follows from the Dominated Convergence Theorem and the assumption that    $f_{x}(t,x) $  is continuous for $x\in(a,b)$. For (\ref{9.8w}) we note that for any $a<u<v<b$ and $t\in (0,\ff)$,
\begin{equation}
f(t,v)-f(t,u)=\int_{u}^{v}f_{x}(t,x)\,dx.\label{dui.2}
\end{equation}
Since $\|f_x(\cdot,x)\|_1$ is uniformly bounded on $(a,b)$ it follows by 
 Fubini's Theorem that,   
\begin{equation}
\int_{0}^{\ff}f(t,v)\,d\mu (t)-\int_{0}^{\ff}f(t,u)\,d\mu (t)=\int_{u}^{v}\int_{0}^{\ff} f_{x}(t,x)\,d\mu (t)\,dx.\label{dui.3}
\end{equation}
Therefore,
\begin{equation} \label{duia}
  \frac{1}{\de}\(\int_{0}^{\ff}\(f(t,x+\de) - f(t,x)\)\,d\mu (t)\)=\frac{1}{\de}\int_{x}^{x+\de}\int_{0}^{\ff} f_{x}(t,x)\,d\mu (t)\,dx.
\end{equation}
Taking the limit as $\de\to 0$ and using (\ref{9.7w}) gives (\ref{9.8w}).\qed

\noindent {\bf  Proof of Lemma \ref{lem-7.3} }    Since $(\si^{\bb}) ^2 (x)$ is symmetric we can restrict our attention to  $x>0$.    Write $1-\cos  \la x =2\sin^2 \la x/2$ in (\ref{lv.21q}). Then, by a change of variables,   when $x\neq 0$, for all $\bb\ge 0$,
\begin{equation} \label{4.59}
 (\si^{\bb}) ^2 (x) =\frac{4}{\pi x}\int_{0}^{\ff}\frac{ \sin^2t/2}{ \bb+ \psi(t/x)}\,dt,  
\end{equation}
 and,  
\bea \label{4.58a}
 (\si^{\bb}) ^2 (x)'&=&-\frac{  (\si^{\bb}) ^2 (x)}{x}+\frac{4}{\pi x}\frac{d}{dx}\(\int_{0}^{\ff}\sin^2 t/2 \frac{ 1}{\bb+  \psi(t/x)}\,d t\)\nn\\
   &:=& -\frac{(\si^{\bb}) ^2 (x) }{x}+\nu_{\bb}(x). 
\eea
 
Note   that using  (\ref{4.14qq}) we have,
\be\label{4.24ww}
 \frac d{dx}\frac{1}{ \bb+\psi(t/x)  }   
 =\frac{1}{x} \frac{(t/x)\psi'(t/x)}{(\bb+  \psi (t/x))^2 } \le   \frac{\ga_{1}}{x  (\bb+  \psi (t/x)) }.
 \ee
It  follows from  (\ref{4.24ww}) and (\ref{4.23qq})  that for any $x_0>0$ there exists an $\ep>0$ 
such that,   
\begin{equation} 
  \sup_{|x-x_{0}|\le \ep}\bigg|\frac{d}{dx}\frac{ 1}{\bb+  \psi(t/x)}\bigg|\in   L^{1}((0,\ff),\sin^2t/2\,dt).
  \label{444.3}
\end{equation}

 Let    $f(t,x)=1/( \bb+ \psi(t/x) )$  
and note that $f_x(t,x)$ is continuous on $(a,b)$ for all $a>0$.    Therefore    it follows from   (\ref{444.3})
 and  Lemma \ref{lem-8.4}   that,
 \begin{equation}
\nu_{\bb}(x)= \frac{4}{\pi x} \int_{0}^{\ff}\sin^2 t/2 \(\frac{d}{dx} \frac{ 1}{\bb+  \psi(t/x)}\)\,d t\label{4.58as}
 \end{equation}
 and is continuous on $x>0$. Writing   $(\si^{\bb}) ^2 (x)'$ as in (\ref{4.58a}) we see 
it too is continuous. Therefore, $(\si^{\bb})^{2}\in C^{1}\(R^1- \{0 \}\)$. 

  Using (\ref{4.14qq}) we see that (\ref{4.24ww}) is positive.
Therefore,
\begin{equation} \label{4.31nn}
  \frac{4}{\pi x}\int_{0}^{\ff}\sin^2t/2\bigg|\frac{d}{dx}\frac{ 1}{\bb+  \psi(t/x)}\bigg|\,d t\le \frac{4\ga_1}{\pi x^2}\int_{0}^{\ff} \frac{  \sin^2t/2}{\bb+  \psi(t/x)} \,d t=\frac{\ga_1(\si^{\bb}) ^2 (x)}{x}.
\end{equation}
Together with (\ref{4.58a}) and (\ref{4.58as}) this gives (\ref{4.58aqe}).

\medskip We now show that for all   $\bb\ge 0$,  $(\si^{\bb})^{2}(x)\in C^{2}\(R^1- \{0 \}\)$.   
 Using (\ref{4.58a}) we see that,
\be   
 (\si^{\bb})^2(x)''  
   = -\frac{  (\si^{\bb})^2(x)'}{x}+  \frac{  (\si^{\bb})^2(x) }{x^2}+
   \frac{d}{dx}\frac{4 }{\pi x}\int_{0}^{\ff}\sin^2t/2\(\frac{d}{dx}\frac{ 1}{ \bb+ \psi(t/x)}\)\,dt.\label{4.61dd} 
\ee
The last term is equal to
\bea \label{8.31}
&& -\frac{4 }{\pi x^2}\int_{0}^{\ff}\sin^2t/2\(\frac{d}{dx}\frac{ 1}{ \bb+ \psi(t/x)}\)\,dt\\
&&\qquad+\frac{4 }{\pi x}\frac{d}{dx}\int_{0}^{\ff}\sin^2t/2\(\frac{d}{dx}\frac{ 1}{ \bb+ \psi(t/x)}\)\,dt\nn
\eea
The material in the beginning of this proof on differentiating under the integral sign  allows us to   write the first term in (\ref{8.31}) as,
\bea 
  && -\frac{4 }{\pi x^2}\frac{d}{dx}\int_{0}^{\ff}\sin^2t/2\(\frac{ 1}{ \bb+ \psi(t/x)}\)\,dt= -\frac{1}{  x^2}\frac{d}{dx}\(x(\si^{\bb}) ^2 (x)\)\nn\\
  &&\qquad=-\frac{(\si^{\bb}) ^2 (x)}{x^2}-\frac{(\si^{\bb}) ^2 (x)'}{x}
\eea
Consequently.  
\begin{equation} \label{4.61dda}
  (\si^{\bb})^2(x)''  
   = -\frac{2(\si^{\bb}) ^2 (x)'}{x}+\frac{4 }{\pi x}\frac{d}{dx}\int_{0}^{\ff}\sin^2t/2\(\frac{d}{dx}\frac{ 1}{ \bb+ \psi(t/x)}\)\,dt.
\end{equation}

%  Differentiating under the integral sign, which we justify below, we have,
%\be   
% (\si^{\bb})^2(x)''  
%   = -\frac{  (\si^{\bb})^2(x)'}{x}+\frac{4 }{\pi x} \int_{0}^{\ff}\sin^2s/2\(\frac{d^2}{dx^2}\frac{ 1}{\bb+  \psi(s/x)}\)\,ds.\label{4.61dd} 
%\ee
  Set 
\be \label{4.24wws}
g(t,x)= \frac d{dx}\frac{1}{  (\bb+ \psi(t/x))  }   
  =  \frac{(t/x^2)\psi'(t/x)}{ (\bb+ \psi(t/x))^2 }. \nn
 \ee 
 Then  
\bea  
 g_{x}(t,x) := \frac {d^2}{dx^2}\frac{1}{  (\bb+ \psi(t/x))  }&=&-\frac{2(t/x^2)\psi'(t/x)}{ x(\bb+ \psi(t/x))^2 } - \frac{(t/x^2)^2\psi''(t/x)}{ (\bb+ \psi(t/x))^2 }\nn \\
   &+&  \frac{2\((t/x^2)\psi'(t/x)\)^2}{ (\bb+ \psi(t/x))^3 }:=I+II+III,\nn 
\eea
  and   is continuous on $(a,b)$ for all $a>0$.
Using (\ref{4.14qq}) we see that,   
\bea 
  |I|&\le& \frac{2\ga_{1}}{x^2  (\bb+  \psi (t/x)) }\\
   |II|&\le&\frac{ \ga _{1}(\ga _{1}-1)}{x^2  (\bb+  \psi (t/x)) }   \nn\\
   |III|&\le& \frac{2\ga^2_{1}}{x^2  (\bb+  \psi (t/x)) }\nn.
\eea
 
  It follows from  (\ref{4.23qq})
  that for any $x_0>0$ there exists an $\ep>0$ 
such that, 
\begin{equation} 
  \sup_{|x-x_{0}|\le \ep}\big| g_{x}(t,x)\big|\in L^{1}([0,\ff), \sin^2\(t/2 \)  \,dt).
  \label{444.3q}
\end{equation}
   It then  follows from Lemma \ref{lem-8.4} that the last term  in (\ref{4.61dda})  is continuous in $x>0$. In addition, we show above that $(\si^{\bb})^{2}(x)\in C^{1}\(R^1- \{0 \}\)$. Therefore, by (\ref{4.61dda}),  $(\si^{\bb})^{2}(x)\in C^{2}\(R^1- \{0 \}\)$. \qed

 \begin{remark} {\rm 
Lemmas   \ref{lem-4.1}--\ref{lem-7.3} require that $\psi$ is the characteristic exponent of a stable mixture. We want to know what functions    $ (\si^{\bb})^2 $ we can have when $\psi$ has this property. This information is given in  \cite[Corollary 9.6.5]{book} which gives the following lemma:

   \begin{lemma} \label{lem-4.6}
  Let $h$ be any function that is regularly varying at infinity with positive index or is slowly varying at infinity and increasing. Then, for any $1<\ga<2$, there exists a L\'evy process that is a stable mixture   for which,
  \begin{equation} 
  (\si^{\bb})^2(x)\sim |x|^{\ga-1}h(\log 1/|x|),\qquad\text{as}\quad x\to 0,\quad\forall\,\bb\ge 0.
  \end{equation}
  If, in addition,
  \begin{equation} \label{4.55}
  \int_1
^\ff\frac{dx}{h(x)}<\ff,
\end{equation}
the above statement is also true when $\ga=2$.
\end{lemma}
}\end{remark}

(The statement in \cite[(9.175)]{book} is given only for  $ (\si^{0})^2$ and includes the statement that it is concave.  If we remove this, as we do in Lemma 
\ref{lem-4.6}, it is easy to see from the material in \cite[Section 9.6]{book},
that what remains holds for   $ (\si^{\bb})^2$ for all $\bb\ge 0$. Note also that
there is a typo in the statement of \cite[Corollary 9.6.5]{book}. What we have in (\ref{4.55}) is given as the integral from 0 to 1. In the proof of \cite[(9.175)]{book} we state that it follows from \cite[(9.178)]{book}) since $\hat \rho(x)\sim 1/h(x)$ at infinity. Clearly this gives (\ref{4.55}).)

 \begin{remark} {\rm \label{rem-sm} In \cite[Lemma 9.6.1]{book} we define a stable mixture to be a L\'evy process with a characteristic exponent $\psi\(\la\)$ that can be written as in (\ref{4.13}) where $\mu$ is a finite positive measure on $(1,2]$ such that, \begin{equation} \label{4.15}
 \int_{1}^{2} \,\frac{d\mu(s)}{2-s}<\ff.
\end{equation}
  With this definition $\mu$ can not have positive mass at $s=2$. Therefore a L\'evy process with this characteristic exponent  does not have a Gaussian component.  We now refer to such processes as stable mixtures without a Gaussian component.  
}\end{remark}

 \begin{remark} {\rm \label{rem-8.3nn} 
Consider (\ref{4.6zzz}) in which $\tau$ is a regularly varying function at 0 with index $3$. Set \begin{equation} 
  \tau(u)=u^3g(u),
\end{equation}
so that  $g(u)$ is a slowly varying function at 0. The function $\psi_1(\la)$
is asymptotic to a regularly varying functions at infinity with index 2 if and only if,
\begin{equation} 
  \int_0^{x}\frac{u^2}{\tau(u)}\,du = \int_0^{x}\frac{1}{ ug(u) }\,du:=h(x),
\end{equation}
 is a slowly varying function at 0. Therefore, when  
 \begin{equation} \label{ald}
  \frac{1}{ xg(x) }=h'(x)
\end{equation}
   we have,
 \begin{equation} \label{4.6znk}
   \psi_{1}(\la)\sim  2\la^2h(1/\la),\qquad\text{as $\la\to\ff$}.
\end{equation}
Since $g(x)$ is a slowly varying function at 0, it follows from (\ref{ald}) that this occurs  whenever $h'$ is a regularly varying function  at zero with index $-1$.  
 
  Here is a simple criteria for determining when this is the case: Without loss of generality we can describe $h$ near 0 by,
 \begin{equation} \label{rem-8.3mm}
  h(x)=\exp\(\int_0^x\frac{\ep(u)}{u}\,du\),
\end{equation}
where $\lim_{x\to 0}\ep(x)=0$. Note that,\begin{equation} 
   h'(x)=\frac{\ep(x)h(x)}{x}.
\end{equation}
  Therefore when $h$ is a slowly varying function at zero, $h'$ is a regularly varying function at zero  with index $-1$ if and only if $\ep(x)$ is a slowly varying function at zero. 
   
}\end{remark}

    \section{ \label{sec-ex1} 
  Excessive functions   that   satisfy the hypotheses of Theorems  \ref{theo-1.2zz}--\ref{theo-1.5}}

The innovation in this paper is to consider permanental processes that have kernels that are not symmetric. We construct these kernels by adding excessive functions to   symmetric potentials as in (\ref{1.3qssj}).   
    In   Part I of Theorems \ref{theo-1.2zz}--\ref{theo-1.5}
    we want to show that for different numbers $d$ there are excessive functions   $f, g\in C^{1}(\De_d(\de))$,  for some $\de>0$.    We also    want to show that  $g(x)\ne Cf(x)$ in $\De_{d}(\de)$   for any constant $C$. This implies that (\ref{1.3qssj}) is not symmetric on $\De_{d}(\de)$.
 In    Part II of Theorems \ref{theo-1.2zz}--\ref{theo-1.5}
    we want to show that for different numbers $d$ there are excessive functions   $f, g\in C^{2}(\De_d(\de))$  for some $\de>0$ such that $f'(d)=g'(d)=0$.     As above we also    want to show that  $g(x)\ne Cf(x)$ in $\De_{d}(\de)$   for any constant $C$   for any $\de>0$.  
    
    In both of these cases we can simplify the condition that $g(x)\ne Cf(x)$ in $\De_{d}(\de)$   for any constant $C$    for any $\de>0$ to 
   \begin{equation} \label{9.1nn}
   g(x)\ne  f(x) \text{ in $\De_{d}(\de)$ for any $\de>0$.}  
\end{equation}  
To see this suppose that $  f$ and $   g$
 are two  excessive functions. Then,
 \begin{equation} \label{9.2nn}
 \wt f(x)= f(x)+  g(d),\qquad\text{and}\qquad  \wt g(x)=  g(x)+  f(d)
\end{equation}
 are  also excessive functions. Since $\wt f(d)=\wt g(d) $ the condition $\wt g(x)\ne C\wt f(x)$ in $\De_{d}(\de)$   for any constant $C$ follows from (\ref{9.1nn}).
  
\medskip The many examples of excessive functions that satisfy the criteria of Part I of Theorems \ref{theo-1.2zz}--\ref{theo-1.5} make it obvious that we can satisfy (\ref{9.1nn}). The situation is more difficult for Part II of Theorems \ref{theo-1.2zz}--\ref{theo-1.5} because of the added requirement that $f'(d)=g'(d)=0$. We concentrate on this case.

\subsection{Theorem \ref{theo-1.2zz}, I } \label{sec-8.1}
 
\medskip In Theorem  \ref{theo-1.2zz}   we consider potential densities of the form 
 \begin{equation} \label{1.21nnq}
u^{\bb}(x,y)=\frac{1}{2\pi}\int_{-\ff}^{\ff}\frac{ \cos\la (x-y)}{\bb+ \psi(\la)}\,d\la,\qquad x,y\in R^{1}, \,\,\bb>0;  
\end{equation}
see (\ref{1.21nn}).   Throughout the paper we take  $\psi$ even.

Set $u^{\bb}(x,y)=u^{\bb}(x-y)$.
Then   when $h\in L^{1}$,
\bea \label{4.19mmq}
 f(x)&=&\int_{-\ff}^\ff  u^{\bb}(x-y) h(y) \,dy\\
 &=& \frac{1}{ 2\pi}\int_{-\ff}^{\ff}e^{i \la x}\frac{\wh h(\la) }{\bb+ \psi(\la)}\,d\la,\nn
\eea
where $\wh h$ is the Fourier transform of $h$.  Clearly $f$ is bounded but, more significantly,  it follows from this that $f^{(k)}$ the $k-$th derivative of $f$ is the  Fourier transform of 
\begin{equation} \label{4.15a}
  C_k \frac{\la^k\wh h(\la) }{\bb+ \psi(\la)},
\end{equation}
where $C_k$ is a constant depending on $k.$ Since $\wh h(0)=\|h\|_1$, and we can find examples of $\wh h$ with any rate of decay as $\la\to\ff$, we can find potentials $f$ with any number of derivatives, including $f\in C^{\ff}$. This shows that the condition   in Theorem \ref{theo-1.2zz}, I. that  $g$ and $ f$ are   $C^{1}(R^1)$ excessive functions for   $Y$ can be easily satisfied for many functions $g$ and $ f$.

\begin{example} {\rm \label{ex-9.1mm}  
  Let   $   h_{ a,b }(x)= I_{[a,b]}(x)$, $a<b$     and,
\begin{equation} \label{9.5nn}
   f_{a,b}(x)=\int_{-\ff}^\ff  u^{\bb}(x-y) h_{a,b}(y) \,dy.
\end{equation}
To show that $f_{a,b}\in C^1(R^1)$ it suffices to show $f_{-1,1}\in C^1(R_1)$, since the other cases are just translations and rescaling of this one. In this case,
\begin{equation} \label{4.15as}
  C_1 \frac{\la\wh h(\la) }{\bb+ \psi(\la)}= C'_1 \frac{\sin \la  }{\bb+ \psi(\la)},
\end{equation}
which is in $L^{1}$. Therefore $f_{-1,1}\in C^1(R_1)$.

 However it is easier to show this by a direct calculation, since
\begin{equation} \label{ex-1}
   f_{a,b}(x)=\int_{a}^b  u^{\bb}(x-y) \,dy=\int_{x-b}^{x-a} u^{\bb}(y) \,dy.
\end{equation}
Therefore,
\begin{equation} \label{ex-2}
  f'_{a,b}(x)=  u^{\bb}(x-a)-  u^{\bb}(x-b) ,
\end{equation}
which shows that $f_{a,b}\in C^1(R_1)$.}

\end{example}

\subsection {Theorem \ref{theo-1.2zz}, II }

It follows from (\ref{4.19mmq}) that when $ (\la\wh h(\la) )/(\bb+ \psi(\la))\in L^1$,
\begin{equation} \label{4.91z}
 f'(x)=\frac{i}{2\pi}\int_{0}^{\ff} e^{i \la x} \frac{\la\wh h(\la) }{\bb+ \psi(\la)}\,d\la.
\end{equation} 
We get the following simple lemma:  
  
\begin{lemma} \label{8.7mm} When $ (\la\wh h(\la) )/(\bb+ \psi(\la))\in L^1$, 
\be
  \lim_{x\to\ff}f'(x)=0.\label{4.911}
 \ee 
  When $h$ is symmetric around $d$, i.e.,    when $h_{d}(y)=h (d+y )$ is an even  function,   
 \be
 f'(d)=0. \label{4.912}
 \ee
 In particular, when $h$ is an even function  
 $ f'(0)=0$.
 
\end{lemma}

\Proof  The limit at infinity follows because   $f'$ is the Fourier transform of an $L^1$ function.  
Furthermore,
\begin{equation}
f'(0)=\frac{i}{2\pi}\int_{0}^{\ff}  \frac{\la\wh h(\la) }{\bb+ \psi(\la)}\,d\la.
\end{equation}
When   $h(x)$ is an even function,    $\la\wh h(\la)$ is an odd function so that $f'(0)=0$. 
Now consider,
\begin{equation} 
  f(x+d)=\int_{-\ff}^\ff  u^{\bb}(x+d-y) h(y) \,dy=\int_{-\ff}^\ff  u^{\bb}(x -y) h_d(y) \,dy:=f_d(x)
\end{equation}
Since $h_d(y)$ is an even function, $f'_d(0)=0$.
Consequently,   $f'(d)=0$.\qed

  Using (\ref{4.15a}) we can find many examples of   potential functions $f\in C^2(R^1)$.

\medskip To simplify the notation in what follows we substitute $x_{0}$ for $d$.

   \begin{lemma} \label{lem-9.2}
  Let  $   h_{ a,b }(x)= I_{[a,b]}(x)$, $a<b$     and, 
\begin{equation} \label{9.12nnr}
   f_{a,b}(x)=\int_{-\ff}^\ff  u^{\bb}(x-y) h_{a,b}(y) \,dy.
\end{equation}
Then  
\begin{equation}\label{9.12nn}
f'_{a,b}(x_0)=0 \,\text{ when }\, x_0=\frac{a+b}{2}.
\end{equation}
Furthermore, when  $u^\bb(x)$ is the potential   density of a L\'evy process with a   characteristic exponent    that is a stable mixture,  $ f_{a,b}\in C^2(R^{1}-\{a,b\})$.  \end{lemma}

\Proof  The statement in (\ref{9.12nn}) follows from Lemma \ref{8.7mm} since $h_{a,b}$ is symmetric about $x_0$.  Furthermore, it follows from (\ref{ex-2}) that
 \be   
   f_{a,b}'(x)=  \frac{(\si^{\bb})^2(x-a)-(\si^{\bb})^2(x-b)}2.\label{9.16}
\ee 
  Therefore the  final sentence in the statement of this lemma  follows from  Lemma \ref{lem-7.3}. \qed

 \begin{lemma} \label{lem-9.3nn} When  $u^\bb(x,y)$ in (\ref{1.21nnq}) is the potential density of a L\'evy process with a   characteristic exponent  $\psi(\la)$  that is a stable mixture,   for any $x_{0}$ there exist   potential functions   $  f,   g \in C^{2}(\De_{x_0}(\de) )$ for some $\de>0$, such that  $  f'(x_0)=  g'(x_0)=0$ and  $  g(x)\ne C f(x)$ for any constant  $C$    for $x\in \De_{x_0}(\de) $,  for   any $\de>0$.  
  \end{lemma}

\Proof   Recall that $u^\bb(x,y)$ is a function of $x-y$.   We often write it as $u^\bb(x-y)$. Let $f(x)=f_{a,b}(x)$   as given in (\ref{9.12nnr}) with  $( a+ b)/2=x_0$   and define, \begin{equation} \label{pole.54}
   g(x)=\int u^{\bb}(x-y)h_{\ov a, \ov b} (y)\,dy ,
\end{equation} 
where $\ov a<a<b<\ov b$ and  $(\ov a+\ov b)/2=x_0$.  Using (\ref{ex-1}) and (\ref{ex-2})    we see that, \begin{equation}  
   g'(x)=u^{\bb}(x-\ov a)-u^{\bb}(x-\ov b)\label{pole.2ww},
\end{equation}
and $g'(x_0) =f'(x_0) =0$.

 Let $\wt C$ be   the constant  for which,
   \be g (x_0) = \wt Cf (x_0).
\end{equation}
  To prove this lemma it suffices to show that   $ g(x)\neq  \wt C f(x)$      for $x\in \De_{x_0}(\de) $,  for any $\de>0$.

 Consider
  \begin{equation} \label{9.15mm}
  l (x)  =g (x)-\wt C f(x).
\end{equation}
Since $u^\bb$ is an even function, 
\be   \label{9.15nn}
  l' (x) =\( u^\bb(x-\ov a)-u^\bb( \ov b-x)\)-\wt C\(u^\bb(x- a)-u^\bb(  b-x)\). 
\ee 
 Set      $c=(b-a)/2$   and $a-\ov a=\ov b-b=q $.   It follows that (\ref{9.15nn}) is, 
\bea 
l'(x) &=&u^{\bb}(c+q+x-x_{0})-u^{\bb}(c+q-( x-x_{0}) ) \\
&&\qquad-\wt C\( u^{\bb}(c+x-x_{0})-u^{\bb}(c-( x-x_{0}))\).\nn
\eea 
Obviously, $l'(x_0)=0$.    Note also that,
\be \label{9.24x}
l''(x_0 )=2\((u^{\bb})'(c+q )-\wt C(u^{\bb})'(c  )  \).   
\ee

  We now show that there exist numbers
 $c$ and $q$ such that  $(u^{\bb})'(c+q )\ne \wt C(u^{\bb})'(c  )$,  so that for these   $c$ and $q$,  $|l''(x_0 )|\ne 0$. To see that   such $c$ and $q$ exist, note that  $\lim_{x\to \ff}u^\bb(x)=0$ as $x\to \ff$  by (\ref{4.911}), which is impossible if $(u^\bb)'(x)$ is a constant   not equal to 0.  One can show that   $u^\bb(x)= 0$,   for all $x$ sufficiently large, is not possible. Since  $|l''(x_0 )|\ne 0$ we see that \begin{equation} \label{9.1nnm}
   f(x)\ne  \wt C g(x) \text{ in $\De_{d}(\de)$ for any $\de>0$.}  
\end{equation}     \qed

  \medskip We continue provide examples of excessive functions $f$ and $g$ which satisfy the conditions of  Theorem \ref{theo-1.2zzk} I and II, Theorem \ref{theo-1.2aa} I and II and Theorem  \ref{theo-1.5} { a}, { b} and { c}. The proof for Theorem \ref{theo-1.2zzk}  uses some results from  the proof for Theorem \ref{theo-1.2aa} so we   consider Theorem \ref{theo-1.2aa} first.

   \subsection {Theorem \ref{theo-1.2aa}, I }
   
%   {\bf  Jay: I am a little confused about he issues here. Let's discuss them.}
% 
%\textbf{Mike: What follows gives a little more information than what we had.}
%{\bf I don't like condition (\ref{9.29mm}). This shows it can be eliminated in some cases. I also think (\ref{9.28nns}) is an interesting formula which may be useful. Also this presentation continues to use sines and cosines as we do throughout the paper.}

  Consider potentials for $u^{(0)}(x,y)$ of the form,
\begin{equation} \label{4.19mme}
 f(x)=\int_{-\ff}^\ff  u^{(0)}(x,y) h(y) \,dy.
\end{equation}
Using obvious trigonometric identities we see that  
 \be  \label{9.28nqn}
  u^{(0)}(x,y)=\frac{1}\pi\int_{0}^{\ff}\frac{ (1-\cos\la x)(1-\cos\la y) +\sin\la x\sin\la y}{\psi(\la)}\,d\la. 
\ee 
Replace $u^{(0)}(x,y)$ in (\ref{4.19mme}) by the integral in (\ref{9.28nqn}) and consider the resulting double integral.
 Let $h_e(y)=(h(y)+h(-y))/2$ and $h_o(y)=(h(y)-h(-y))/2$. When $(|\la|\wedge 1)/\psi(\la)\in L^1$   we can interchange the order of integration   and get 
\be  \label{9.28nns}
  f(x)=\frac{1}\pi\int_{0}^{\ff}\frac{ (1-\cos\la x)\(\wh h(0)-\wh h_e(\la)\)+\sin\la x\,\wh h_o(\la)}{\psi(\la)}\,d\la, 
\ee
since $\wh h_e(0)=\wh h(0)$. We write this as,
\be  \label{9.28nnsaa}
  f(x)=\phi(x)\|h\|_1-\frac{1}\pi\int_{0}^{\ff}\frac{ (1-\cos\la x)\,\wh h_e(\la) -\sin\la x\,\wh h_o(\la)}{\psi(\la)}\,d\la. 
\ee

 \begin{lemma} \label{prop-8.3qx}   Let $u^{(0)}(x,y)$  be the potential   density of a L\'evy process with a   characteristic exponent   $\psi(\la) $   that is a stable mixture   killed the first time it hits 0. Then for all positive  even    functions  $h\in L^1$, or  for all positive      functions  $h\in L^1$ satisfying  
\begin{equation} \label{9.29mm}
   \frac { | \la \wh h(\la) | } {\psi(\la)}\in L^1((0,\ff),d\la),
\end{equation}  the corresponding  potential  $f \in  C^1(R^1-\{0 \})$.
  \end{lemma}  

 \Proof By Lemma  \ref{lem-7.3},  since    $2\phi(x)=(\si^{(0)})^2(x)$, $\phi'(x)$ is continuous  on $R^1-\{0 \}$. Also, when $h$ is even or when (\ref{9.29mm})
 holds it follows from Lemma \ref{lem-8.4} that for all  $x \in   R^1-\{0 \}$,    
\be  \label{9.28nnsaaq}
  f'(x)=\phi'(x)\|h\|_1+\frac{2}\pi\int_{0}^{\ff}\frac{  \la\sin\la x \, \wh h_e(\la) +\la\cos\la x\,\wh h_o(\la)}{\psi(\la)}\,d\la. 
\ee
  Since  for all $0<a<b$,
\begin{equation} 
  \sup_{x\in (a,b)} \frac{|\la\sin\la x \,  \wh h_e(\la)  |}{\psi(\la)} \in L^1((0,\ff),d\la),
\end{equation}and \begin{equation} 
  \sup_{x\in (0,\ff)} \frac{|\la\sin\la x \,  \wh h_e(\la) +\la\cos\la x\,\wh h_o(\la)|}{\psi(\la)}\le \frac{2|\la \wh h(\la)|}{\psi(\la)},
\end{equation}
 the integral in (\ref{9.28nns}) is continuous in $x \in   R^1-\{0 \}$.     \qed

\begin{example} {\rm    Let  $  h(x)= h_{ a,b }(x)= I_{[a,b]}(x)$, $a<b$ in (\ref{4.19mme}). Then,
\begin{equation} 
  \wh h(x)=\frac{e^{i\la b}-e^{i\la a}}{\la}
\end{equation}
and   \begin{equation} 
  \frac{\la|\wh h(\la)| }{\psi(\la)}\le \frac{2(\la|a-b|\wedge 1)}{\psi(\la)}.
\end{equation} Therefore when,
\begin{equation} 
  \int_{0}^\de\frac{\la}{\psi(\la)}\,d\la<\ff,\qquad \text{for some $\de>0$,}
\end{equation} 
the corresponding potential $f$ in (\ref{4.19mme}) is in $C^1(R^1-\{0 \})$.
}\end{example}

   \subsection {Theorem \ref{theo-1.2aa}, II }

   We begin with an extension of \cite[Lemma 9.6.3]{book}.

\begin{lemma} \label{lem-8.3} When $\psi(\la)$ is the characteristic exponent of a stable mixture and is not equal to $C\la^2$ for some $C>0$, 
\begin{equation} 
  (\si^{(0)})^2(x)=\frac{2}{\pi}\int_{0}^{\ff}\frac{ 1-\cos\la x}{\psi(\la)}\,d\la\in C^{2}\(R^1- \{0 \}\)
\end{equation}
and is strictly concave on $[0,\ff)$. In particular  $((\si^{(0)})^2(x))''<0$ for all $x>0$. 
\end{lemma}

\Proof The fact that   $ (\si^{(0)})^2(x) \in C^{2}\(R^1- \{0 \}\)$ is given in  Lemma \ref{lem-7.3} .
For concavity, refer to the proof of  \cite[Lemma 9.6.3]{book}. The proof was written under the assumption that $\psi(\la)=o(\la^2)$, as $\la\to\ff$. It only needs a minor change  to give this lemma. In \cite[Lemma 9.6.3]{book} we introduce the function $g(x)$ and show that $g''(x)\le 0$. To prove this lemma we must sharpen this to $g''(x)< 0$. That is, we need strictly less than in  \cite[(9.167)]{book}. We get this by showing that the Schwartz inequality applied in \cite[Lemma 9.6.3]{book} is strict. It follows from the second line of \cite[(9.169)]{book} that for the Schwartz inequality to be an equality we must have,
\begin{equation} 
  \frac{s-1}{s}=C'(s-1)s,\qquad \forall\, s\in (1,2]
\end{equation}
for some constant $C'>0$, which is only possible if $\mu(s)$, in \cite[Lemma 9.6.3]{book}, is concentrated at  single point. Other parts of the proof of \cite[Lemma 9.6.3]{book} show that this point has to be 2. \qed

  \begin{lemma} \label{prop-8.3q} Let $u^{(0)}(x,y)$  be the potential   density of a L\'evy process with a   characteristic exponent   $\psi(\la) $   that is a stable mixture  and  $\psi(\la)\ne C\la^2$ for any $C>0$.  Then for any $x_{0}$ there exist excessive functions   $\wt f,\wt g$ for $u^{(0)}(x,y)$ such that  $\wt f,\wt g  \in C^{2}(\De_{x_0}(\de) )$ for some $\de>0$,   $\wt f'(x_0)=\wt g'(x_0)=0$ and  $\wt g(x)\ne C\wt f(x)$ for $x\in \De_{x_0}(\de) $, for any constant  $C$.      
  \end{lemma}

 \Proof  
  Let  $f_{a,b}$ be as given in (\ref{ex-1}) where $h_{a,b}(x)=I_{[a,b]}(x)/(b-a)$, $0<a<b$. Then, 
  \bea
f_{a,b}(x) 
&=&\phi (x)  +  \int_{-\ff }^{ \ff}  {\phi (y)}h_{a,b} (y) \,dy-  \int_{a}^b   {\phi (y-x)} h_{a,b} (y) \,dy \label{9.45mm}\\
&=&\phi (x)  +  \int_{a}^b \frac{\phi (y)}{b-a}  \,dy-  \int_{a}^b   \frac{\phi (y-x)}{b-a} \,dy.\nn
\eea 
 Using  (\ref{ex-2})  we see that, 
\be\label{9.43r}
f_{a,b}'(x)
=\phi' (x) -  \frac{ \phi (x-a)-\phi ( -b+x)
}{b-a}:=\phi' (x)-m_{-b+x,x-a}.
\ee
Note that $m_{-b+x,x-a}$ is   the slope of the line  segment from  $(-b+x , \phi( -b+x))$  to $(x-a, \phi(x-a))$.
  
  Since $\phi(x)$ is an non--negative continuous even function with $\phi(0)=0$  and  $\phi(x)$  increasing and strictly concave for $x> 0$,
 for any $x_{0}>0$, we can find  $0<a<x_{0}<b$ such that
the slope of the line segment from  $(x_0-b, \phi( x_0-b))$ to $( x_0-a, \phi(  x_0-a))$ is equal to $ \phi' (x_0)$, i.e.,  
\be\label{9.43ss}
 \phi' (x_0) = \frac{ \phi (x_0-a)-\phi ( -b+x_0)
}{b-a}=m_{ x_0-b, x_0-a}.
\ee 
It follows from (\ref{9.43r}) that $f'_{a,b}(x_0)=0$. 

  For the following argument it is useful 
to draw the graph of of the line segment between $(x_0-b, \phi( x_0-b))$ to $( x_0-a, \phi(  x_0-a))$. (It is important to note that     $x_0-b<0,  x_0-a>0$.)

It is easy to see that there exist   $b'<b$ and $a'>a$ such that the line segment joining   $\(-b'+x_0, \phi(-b'+x_0)\) $ and $\(x_0-a', \phi(  x_0-a')\)$ have the same slope, i.e. that,
 \begin{equation} 
  m_{ x_0-b, x_0-a}=m_{ x_0-b', x_0-a'}.
\end{equation} 
Define $f_{a',b'}(x)$ as in (\ref{9.45mm}). Clearly  
\be 
f'_{a',b'}(x_0)=f'_{a,b }(x_0)
\ee
Consequently,   as we point out in the paragraph containing (\ref{9.1nn}), to complete this proof we need only show that $f _{a',b'}(x)\ne f _{a,b }(x )$ in $\De_{x_0}(\de)$, for  all $\de>0$. We do this by showing that    $   f'_{a',b'}\(x\)<  f'_{a,b}\(x\)$ for $x\in \De_{x_0}(\de) $ for some $\de>0$.  
By (\ref{9.43r}) this is equivalent to showing that for   all $\de>0$, sufficiently small,\be
  \frac{ \phi (x_0+\de-a)-\phi ( -b+x_0+\de)
}{b-a}<\frac{ \phi (x_0+\de-a')-\phi ( -b'+x_0+\de)
}{b'-a'}.
\ee
 This follows because
%This is because for   all $x-x_0$ sufficiently small $x -  a'<x -  a$ and 
%$  - b+x< -  b'+x$. Therefore, since $\phi'(|\cdot|)$ is decreasing we have,  
%{\bf  the point is that $\phi '(-  b+x  )<0$! }
\be 
\phi '(x_0 -  a')>\phi' (x_0 -  a)\qquad\text{and}\qquad -\phi '(-  b'+x_0  )>-\phi' (-  b +x_0 ) ,\label{99.1}
\ee
since $\phi'(x)$ is decreasing when $x>0$ and increasing when $x<0$,  i.e., takes on smaller negative values, and $b'-a'<b-a$.  

\medskip It follows from Lemma \ref{lem-9.2} that  $    f_{a,b}\(x\)$ and $   f_{a',b'}\(x\)$ are in $C^{2}(\De_{x_0}(\de) )$. 
 \qed

   The case when $u^{(0)}(x,y)$    is the  potential   density of a L\'evy process with     characteristic exponent $\psi(\la)=C\la^2$, in which case the process is a constant time Brownian motion,  in covered in Example  \ref{ex-9.4}.   
   
     \subsection {Theorem \ref{theo-1.2zzk}, I }
In this theorem,   \be 
 v^{\bb} (x,y) =u^\bb(x-y)-\frac{u^\bb(x)u^\bb(y)}{u^\bb(0)},\nn
 \ee
with potentials,
\bea \label{4.19qq}
\wt f(x)&=&\int_{-\ff}^\ff  v^{\bb}(x,y) h(y) \,dy\\
  &=&\int_{-\ff}^\ff  u^{\bb}(x-y) h(y)\,dy-\frac{u^\bb(x) }{u^\bb(0)}f(0)\nn \\
   &  =&f(x)-\frac{u^\bb(x) }{u^\bb(0)}f(0)\nn 
\eea
where $f$ is the potential in (\ref{4.19mmq}).

As we pointed out in Section \ref{sec-8.1}  we can find examples of $\wh h$ with any rate of decay as $\la\to\ff$.  Consequently we can find  functions  $f$ with any number of derivatives, including $f\in C^{\ff}$.  Note that,
\be  \label{4.19qqx}
\wt f'(x)  = f'(x)+\((\si^\bb)^2(x)\)'\frac{f(0) }{2u^\bb(0)}.  
\ee

By Lemma \ref{lem-7.3}, $(\si^{\bb})^{2}(x)\in C^{1}\(R^1- \{0 \}\)$. Therefore, when  $u^\bb(x)$ is the potential of a L\'evy process with a   characteristic exponent    that is a stable mixture    $ \wt f(x)\in C^1(R^1-\{0 \})$.

 \subsection { Theorem \ref{theo-1.2zzk}, II }
\begin{lemma} \label{lem-9.7mm}
 When  $ v^{\bb}(x,y)$    is the potential  density of a L\'evy process with a   characteristic exponent $\psi(\la) $   that is a stable mixture and $ \psi(\la)\ne C\la^2$ for any $C>0$, then for any   $x_{0}\ne 0$ there exist excessive functions $\wt f,\wt g$ for $u^{(0)}(x,y)$ such that  $\wt f,\wt g  \in C^{2}(\De_{x_0}(\de) )$ for some $\de>0$,     $\wt f'(x_0)=\wt g'(x_0)=0$ and  $\wt g(x)\ne Cf(x)$ for $x\in \De_{x_0}(\de) $, for any constant  $C$. 
\end{lemma}

\Proof    Let $Z $ be a symmetric L\'evy process. It follows from  \cite[Example 4.5.4]{book} that  for any $\bb>0$, 
the $\bb$--potential density $u^{(\bb)}(x,y)  $ of the process obtained by killing  $Z$ the first time it hits $0$ is the same as   $ v^{\bb}(x,y)$, the $0$--potential density of the process obtained by first killing $Z$ at the end of an independent exponential 
time with mean $1/\bb$, and then killing it  the first time it hits $0$.  That is,
\begin{equation} \label{obv.1}
 u^{(\bb)}(x,y)= v^\bb(x,y), \qquad \forall \bb>0.
\end{equation}
  Consider the  resolvent equation for $\{u^{(\bb)}(x,y),\, \bb\geq 0\}$ where we substitute  $ v^\bb(x,y)$ for $u^{(\bb)}(x,y)$ when $\bb>0$:
  \begin{equation}
\int_{-\ff}^\ff u^{(0)}(x,y)h(y)\,dy-\int_{-\ff}^\ff v^{\bb}(x,y)h(y)\,dy=\bb\int_{-\ff}^\ff v^{\bb}(x,y) U^{(0)}h(y) \,dy,\label{resv.10}
\end{equation}
 where  \be
 U^{(0)}h(z)=\int_{-\ff}^\ff u^{(0)}(z,y)h(y)\,dy. 
 \ee
 Combining these two equations we see that, 
\begin{equation}
\int_{-\ff}^\ff u^{(0)}(x,y)h(y)\,dy= U^{(0)}h(x)= \int_{-\ff}^\ff v^{\bb}(x,y)\( h+\bb U^{(0)}h\) (y)\,dy.     \label{resv.11}
\end{equation}

Let $ h_{a,b}(x)=I_{[a,b]}(x)/(b-a)$, $0<a<b$ and $  h_{a', b'}(x)=I_{[a', b']}(x)/(   b'-  a')$, $0<a'<b'$, be the functions in (\ref{9.45mm}) that determine $f_{a , b }(x)$ and $f_{a', b'}(x)$ in Lemma \ref{prop-8.3q}. I.e.,
\begin{equation}
U^{(0)} h_{a,b}(x)=f_{a,b}(x), \quad\mbox{ and }\quad U^{(0)}h_{a', b'}(x)=f_{a', b'}(x).\label{resv.12}
\end{equation}
 Consider the two potentials,
\begin{equation}
k_{a,b}(x) =\int v^{\bb}(x,y)\(   h_{a,b}+\bb f_{a,b}\) (y)\,dy, \label{resv.13}
\end{equation}
and\begin{equation}
k_{a', b'}(x) =\int v^{\bb}(x,y)\(   h_{a', b'}+\bb f_{a', b'}\) (y)\,dy, \label{resv.15}
\end{equation}
for $v^{\bb}(x,y)$.
By (\ref{resv.11})  
 \begin{equation}
k_{a,b}(x) =f_{a,b}(x) \hspace{.2in}\mbox{ and } \hspace{.2in} k_{a', b'}(x)=f_{a', b'}(x).   
\end{equation} 
Following the proof of Lemma \ref{prop-8.3q} we see that 
$k_{a,b}'(x_{0}) 
=k_{a', b'}'(x_{0})=0$ and $k_{a,b}'(x ) 
-k_{a', b'}'(x )>0$   for $x\in \De_{x_{0} }(\de)$ for some $\de>0$. %Set,
%\begin{equation} 
%  \wt k_{a,b}(x) =k_{a,b}(x)+k_{a',b'}(x_0)\quad\text{and} \quad \wt k_{a',b'}(x) =k_{a',b'}(x)+k_{a ,b }(x_0).
%\end{equation}
%Consider (\ref{9.1nn}) and (\ref{9.2nn}) with $f_{a,b}$ and $\wt f_{a,b}$ replaced by $k_{a,b}$ and $\wt k_{a,b}$ and similarly for $f_{a',b'}$ and $\wt f_{a',b'}$. 
 
The rest of the proof follows  the end of the proof of Lemma \ref{prop-8.3q}.  \qed

   The case when $v^{\bb}(x,y)$    is the $\bb$ potential   density of a L\'evy process with     characteristic exponent $\psi(\la)=C\la^2$, in which case the process is a diffusion,  in covered in Example  \ref{ex-9.3}.

  \subsection{Theorem \ref{theo-1.5} }\label{sec-9.9a}

 In Parts a. b. and c. of  
Theorem \ref{theo-1.5} we require that  the   excessive functions  have a derivative equal to 0 at 
some point in their state space as in Part b. of Theorems \ref{theo-1.2zz}--\ref{theo-1.2aa}. This is because   $\si^{2}(d+t,d)$ in (\ref{1.34sz}) is equal to a constant times $|x|$ in each of a. b. and c. Therefore $(i)$
of Theorem \ref{theo-1.2z} is not satisfied. 

Theorem \ref{theo-1.5} is for diffusions. Sophisticated tools exist for studying them. We describe their potentials in much greater generality than those in Part b. of Theorems \ref{theo-1.2zz}--\ref{theo-1.2aa}.
 We treat each of the three parts of Theorem \ref{theo-1.5} separately.

 In the discussion preceding Theorem \ref{theo-1.5} we define $\cal Z$ to a diffusion in $R^{1}$ which is regular and without traps.   In the examples we give     $\cal Z$ is  symmetric with respect to Lebesgue measure and   has generator,
\begin{equation}
L=\frac{1}{2 }\frac{d}{dx} b(x)\frac{d}{dx}=\frac{1}{2 }b(x)\frac{d^{2}}{dx^{2}}+\frac{1}{2 }b'(x)\frac{d}{dx},\label{gen.8}
\end{equation}
where $b\in C^{1}(R^{1})$ and is strictly positive.   (This is the most general  form of a symmetric generator with continuous coefficients, of a diffusion in $R^1$.)

\subsubsection{Theorem \ref{theo-1.5}, \textbf{a.}} \label{sec-9.7.1}
 
   We define  $Z$ to be   $\cal Z$  killed  at the end of an independent exponential time with mean $1/\bb$, $\bb>0$. The potential density for $Z$   has the form,
 \begin{equation} \label{gen.1m}
\wt u^\bb (x,y)= \left\{
 \begin{array} {cc}
 p(x)q(y),& \quad x\leq y  
 \\
 q(x)p(y),& \quad y\leq x   
\end{array}  \right. ,
\end{equation}
where $p,q\in C^{2}(R^{1})$  are    positive    strictly convex functions  with the properties that $p$ is  increasing and  
$q$ is decreasing.

\begin{lemma} \label{lem-9.9}
  For any   $x_{0 }\in R^1$ there exist potentials $\wt f,\wt g$ for $\wt u^{\bb}(x,y)$ such that  $\wt f,\wt g  \in C^{2}(\De_{x_0}(\de) )$ for some $\de>0$,     $\wt f'(x_0)=\wt g'(x_0)=0$ and  $\wt g(x)\ne Cf(x)$ for $x\in \De_{x_0}(\de) $, for any constant  $C$. 
\end{lemma}

\Proof   Let  $h$ be a continuous non--negative integrable  function and set,
\bea \label{9.58ww}
  f(x)&=&\int^{\ff}_{-\ff}\wt u^\bb(x,y)h(y)\,dy\\
  &=&  q(x)\int^{x}_{-\ff}p(y)h(y)\,dy  +  p(x) \int^{\ff }_{x}q(y)h(y)\,dy  .\nonumber 
\eea
 Using   the fact that  $h$ is continuous we have, 
 \begin{eqnarray}  
 \label{gen.3}  f'(x)  &  = &\frac{d}{dx}\( q(x)\int^{x}_{-\ff}p(y)h(y)\,dy \)+ \frac{d}{dx}\(p(x) \int^{\ff }_{x}q(y)h(y)\,dy\) \nonumber\\
 &= &  q'(x)\int^{x}_{-\ff}p(y)h(y)\,dy  +  p'(x) \int^{\ff }_{x}q(y)h(y)\,dy.
 \end{eqnarray}
Consequently,
 \be
   f' (x_0)=   q'(x_0)\int^{x_0}_{-\ff}p(y)h(y)\,dy  +  p'(x_0) \int^{\ff }_{x_0}q(y)h(y)\,dy.  \label{gen.3ar}
 \ee
Let  $h_{1}(y)=1_{\{(-\ff,x_0]\}}h(y)$ and  $h_{2}(y)=1_{\{[x_0,\ff)\}}h(y)$.  Clearly, $f' (x_0)=0$
when
 \be
\int^{x_0}_{-\ff}p(y)h_{1}(y)\,dy=-\frac{p'(x_0)}{   q'(x_0)}   \int^{\ff }_{x_0}q(y)h_{2}(y)\,dy.  \label{gen.3asm}
 \ee
Note that $p,q$ and $-{p'(x_0)}/{   q'(x_0)}$ are all strictly positive. 
It is obvious that for any positive continuous function $h_{2}(y)$, for which  the right--hand side of (\ref{gen.3asm}) is finite, there exist an infinite number of different  positive continuous functions $ h_{2}(y)$ with $h_{1}(x_0)=h_{2}(x_0)$,  for which  (\ref{gen.3asm}) is satisfied.

Let $( h_1(y), h_2(y))$ and $(\wt h_1(y),\wt h_2(y))$ be two such pairs and consider,
\bea
h (y) &=& \wt h_1(y)1_{\{(-\ff,x_0]\}}+\wt h_2(y)I_{\{(x_0,\ff \}}\qquad\text{and,}\\
 \wt h(y) &=& \wt h_1(y)1_{\{(-\ff,x_0]\}}+\wt h_2(y)I_{\{(x_0,\ff \}}.\nn
\eea
Considering the infinite number of way we can choose these functions it is clear that we can choose  a $\wt h(y)$ that is not equal to a constant multiple of $ h(y)$. It now follows from
Lemma \ref{cor-wron1}  below
that the corresponding potentials are not multiples of each other in any neighborhood of $x_0$.\qed
\begin{lemma}\label{cor-wron1}    When 
\begin{equation}
\int^{\ff }_{-\ff}\wt u^\bb(x,y)h_{1}(y)\,dy=C\int^{\ff }_{-\ff}\wt u^\bb(x,y)(x,y)h_{2}(y)\,dy\label{gen.11},
\ee
for some $C$ in some interval $I$,   $h_{1}(y)=Ch_{2}(y)$ on $I$. 
\end{lemma}

This lemma is an immediate corollary of the following theorem:
  
\bt\label{lem-wron} There is a constant $c_{p,q}>0$ such that for any continuous non-negative integrable  function $h$
\begin{equation}
\( L-\bb\) \int^{\ff }_{-\ff}\wt u^\bb(x,y)h(y)\,dy=-c_{p,q}h(x).\label{gen.10r}
\end{equation}
\et

\noindent{\bf  Proof of Lemma \ref{lem-wron}}
Using (\ref{gen.3})   we see that,
  \begin{eqnarray}
 f''(x)   &=& \frac{d}{dx}  \( q'(x)\int^{x}_{-\ff}p(y)h(y)\,dy\)   +\frac{d}{dx}  \(  p'(x) \int^{\ff }_{x}q(y)h(y)\,dy\) \nonumber \\&  = &  q''(x)\int^{x}_{-\ff}p(y)h(y)\,dy  +  p''(x) \int^{\ff }_{x}q(y)h(y)\,dy    \label{gen.3b}\\
 &&\qquad    +\(q'(x)p(x)-q(x)p'(x) \) h(x). \nonumber
 \end{eqnarray} 
   We have,
  \begin{equation} \label{gen.8a}
  \( Lf\) (x)=\frac{1}{2}b(x)f''(x)+\frac{1}{2}b'(x)f'(x)
\end{equation}
and similarly for $ \( Lp\) (x)$ and $ \( Lq\) (x)$.
Consequently, using    (\ref{gen.3}) (\ref{gen.3b}) and   (\ref{gen.8a}) we see that,
 \begin{eqnarray}
 ( Lf)(x)
 \label{gen.9} & = &  \( Lq\) (x)\int^{x}_{-\ff}p(y)h(y)\,dy  +   \( Lp\)(x) \int^{\ff }_{x}q(y)h(y)\,dy  \\
 &&  \qquad+\frac{1}{2 }b(x)\(q'(x)p(x)-q(x)p'(x) \) h(x). \nonumber
 \end{eqnarray}  

Note that,\begin{equation} \label{9.57}
  Lq=\bb q\qquad\text{and}\qquad Lp=\bb p.
\end{equation} 
%\textbf{  I suggest: This is a well known result but we can not find a direct proof. 
 %Breiman, \cite{Breiman} gives this as Problem 19 at the end of   Section 16. %Apparently it is a consequence of Theorem 16.75 which he does not prove. He refers %the reader to,  \cite[Ito and McKean page 149 ff.]{IM}. We leave this to the ambitious %reader.}
 See  the end of Section 16.12 in \cite{Breiman}, which refers to \cite{IM} for proofs.

  Using (\ref{gen.9}) and (\ref{9.57})   we can finally write,\begin{equation}
\( L-\bb\) \int^{\ff }_{-\ff}\wt u^\bb(x,y)h(y)\,dy=\frac{1}{2 }b(x)\(q'(x)p(x)-q(x)p'(x) \)h(x).\label{gen.10}
\end{equation}
Clearly the Wronskian, 
$W(x)=  ({1}/{2 })b(x)\(q'(x)p(x)-q(x)p'(x) \),$  is  negative. We now show that it is   a constant.   This follows because  
\begin{eqnarray}
 W'(x) 
&=&  \frac{1}{2 } \( b'(x)\(q'(x)p(x)-q(x)p'(x) \)\) \nonumber\\
&&\quad+ \frac{1}{2 } \( b(x)\(q''(x)p(x)-q(x)p''(x) \)\)\\&=&\( Lq\)(x)\,p(x)-\,\( Lp\)(x)q(x)=0, \nonumber
\end{eqnarray}
by (\ref{9.57}).  
We  label it   $-c_{p,q}$.  Using this in (\ref{gen.10}) gives (\ref{gen.10r}).\qed

\subsubsection{Theorem \ref{theo-1.5}, \textbf{b.}}\label{sec-9.9b}

 We consider    $Z'$ to be the  process with state space  $T=(0,\ff) $ which is obtained by starting   $  Z$ in $T$ and then  killing it  the first time it hits 0. The potential of $Z'$ is   
\bea
\wt v^\bb  (x,y)&=&\wt u^\bb (x,y)-\frac{\wt u^\bb (x,0)\wt u^\bb (0,y)}{\wt u^\bb(0,0)},\qquad  x,y>0, \label{gen.12}\\
  &=&\wt u^\bb (x,y)-\frac{p(0)}{q(0)}q(x)q(y),\qquad  x,y>0. \nn
  \eea
  
  \begin{lemma} \label{lem-9.11}
  For any   $x_{0 }\in R^1$ there exist potentials $\wt f,\wt g$ for $\wt v^{\bb}(x,y)$ such that  $\wt f,\wt g  \in C^{2}(\De_{x_0}(\de) )$ for some $\de>0$,     $\wt f'(x_0)=\wt g'(x_0)=0$ and  $\wt g(x)\ne Cf(x)$ for $x\in \De_{x_0}(\de) $, for any constant  $C$. 
\end{lemma}

\Proof The proof is a minor modification of the proof of Lemma \ref{lem-9.9}. We have,
  
 \begin{eqnarray}
 &&\( L-\bb\) \int^{\ff }_{0}\wt v^\bb(x,y)h(y)\,dy
 \\
 &&\qquad=\( L-\bb\) \int^{\ff }_{0}\wt u^\bb(x,y)h(y)\,dy-\frac{p(0)}{q(0)}\(\( L-\bb\)q \) (x) \int^{\ff }_{0}  q(y) h(y)\,dy  \nonumber\\
 &&\qquad =\( L-\bb\) \int^{\ff }_{0}\wt u^\bb(x,y)h(y)\,dy=-c_{p,q}h(x),  \nonumber
 \end{eqnarray}
  where we use (\ref{9.57}) to see that $Lq=\bb q$ followed by Lemma \ref{lem-wron}. Therefore, as we  point in Lemma \ref{cor-wron1}, if
\begin{equation}
\int^{\ff }_{-\ff}\wt v^\bb(x,y)h_{1}(y)\,dy=C\int^{\ff }_{-\ff}\wt v^\bb(x,y)(x,y)h_{2}(y)\,dy\label{gen.11q}
\end{equation}
for some $C$ in some interval $I$, then $h_{1}(y)=Ch_{2}(y)$ on $I$.

Fix $ x_{0}>0$. We now show how to find  potentials
\begin{equation} 
   f (x)=\int^{\ff }_{0}\wt v^\bb(x,y)h (y)\,dy
\end{equation}
 for which $f' (x_0)=0$.  Using  (\ref{9.58ww})  we have,
\bea
\label{18.3} 
 f'(x)& =&p'(x)\int_{x}^\ff  q( y)h(y)\,dy+q'(x)\int_{0}^{x}  p( y)h(y)\,dy   \nonumber\\
&&\qquad-\frac{p(0)}{q(0)}q'(x)\int_{0}^\ff q(y) h(y)\,dy        \nonumber\\
&=& \( p(x)-\frac{p(0)}{q(0)}q(x)\)' \int_{x}^\ff  q( y)h(y)\,dy   \nonumber\\
&&\qquad +q'(x)  \int_{0}^{x}  \(p( y)-\frac{p(0)}{q(0)} q(y)  \) h(y)\,dy.       \nonumber
\end{eqnarray}
Set   $\rho (x)=p(x)- \({p(0)}/{q(0)}\)q(x)$ and write this as
\begin{equation}
  f'(x)= \rho'(x)\int_{x}^\ff  q( y)h(y)\,dy+q'(x)\int_{0}^{x} \rho(y) h(y)\,dy.  \label{18.18.}
\end{equation}

Since $ \rho(0)=0$ and $\rho'(x) >0$ for $x>0$ we see that $\rho(y)>0$ for all $y>0$.  In addition $q'(x)<0$. Consequently $f'(x_0)=0$ for all positive continuous functions $h$ for which, 
 \begin{equation}
\int_{0}^{x_0}  \rho(y) h(y)\,dy=\frac{- \rho'(x_0)}{q'(x_0)}\int_{x_0}^\ff  q( y)h(y)\,dy.  \label{18.5}
\end{equation}
This is similar to the situation  in the proof of Lemma \ref{lem-9.9}. Using the argument given there we see that there exist an infinite number of functions satisfying the hypotheses of this lemma.\qed

 \begin{example} {\rm  \label{ex-9.3} 
In (\ref{exe.13ja}) we show that  for $\bb>0$,
  \begin{equation}      u^{\bb}_C(x-y) :=\frac{1}{ \pi}\int_{0}^{\ff}\frac{ \cos\la (x-y)}{\bb+ C\la^2}\,d\la= {e^{- ({ \bb   C})^{1/2}\,| x-y|}\over    2(\bb C) ^{1/2} },\quad x,y\in  R^{ 1}.  \label{exe.13jab}
\end{equation} 
We can write this as in (\ref{gen.1m}) with,
\begin{equation} 
  p(x)= {e^{({ \bb  /C})^{1/2}\,   x }\over    \sqrt 2(C\bb )^{1/4} },\qquad\text{and}\qquad  q(y)= {e^{-({ \bb  /C})^{1/2}\,   y }\over \sqrt 2(C\bb )^{1/4}  }.
\end{equation}

 We   consider $v^{\bb}(x,y)$ for these values of $p$ and $q$ and can use Lemma \ref{lem-9.11} to augment Lemma \ref{lem-9.7mm} when $ \psi(\la)\ne C\la^2$. 
 
 Note that,
 \be
 \tau(d) =\frac{1}{2( \bb   C)^{1/2}  
}\frac{d}{dx} \log\( e^{2( \bb /  C)^{1/2}x} \)\bigg |_{x=d}=\frac{1}{C}.
  \ee
Using (\ref{1.27mm}), or more generally Lemma \ref{lem-4.9}, we see  that
 \begin{equation} 
   (\si_C^{\bb})^2\(x\) \sim \tau(d)x,\qquad \text{as $x\downarrow 0$}.
\end{equation}
 This verifies  that the denominators in (\ref{1.34szqa}) and  (\ref{1.34szq}) are the same.
}\end{example}

 \subsubsection{Theorem \ref{theo-1.5}, \textbf{c.}}\label{sec-9.9}
 
In the discussion preceding Theorem \ref{theo-1.5} and in Section \ref{sec-9.7.1} we define the diffusion $\cal Z$ in $R^{1}$ which has the  generator
\begin{equation}
L=\frac{1}{2 }\frac{d}{dx} b(x)\frac{d}{dx}, \label{gen.8qm}
\end{equation}
where $b\in C^{1}(R^{1})$ and is strictly positive.   Let  $s(x)\in C^{2}([0,\ff))$ be a positive strictly increasing function with $s(0)=0$ and $\lim_{x\to\ff}s(x)=\ff$. Let $b(x)=1/s'(x)$, for $x\in [0,\ff)$.

\bt\label{theo-diffs}  Let $\ov Z$ be the process with state space  $ (0,\ff)$ that is obtained by starting  $\cal Z$ in $ (0,\ff)$ and then  killing it  the first time it hits $0$. Then    $\ov Z$ has potential density  
\begin{equation}
u(x,y)=2  \(s(x)\wedge s(y)\),\qquad  x,y>0,\label{d2.1}
\end{equation}
with respect to Lebesgue measure. 
\et

 \Proof
  It follows from \cite[VII. Exercise 3.20]{RY} that $s(x)$ is the scale function for the diffusion in $R^{1}$ with generator $L$ defined in   (\ref{gen.8qm}) with $b(x)=1/s'(x)$.  In addition the speed measure of this diffusion is two times Lebesgue measure.  Taking the limit $a\to 0$ and  $b\to\ff$ in \cite[VII. Corollary 3.8]{RY} shows that the potential density of  $\ov Z$,   the diffusion  in $(0,\ff)$ with  infinitesimal generator $L$ killed the first time it hits $0$, is given by  (\ref{d2.1}).  \qed

Note that for functions $f\in C_b^2(0,\ff)$,
\begin{equation}
L f(s(x))=\frac{1}{2 }\frac{d}{dx} \frac1{s' (x)}\frac{d}{dx}f(s(x))=\frac1{2}\frac{d}{dx}f'(s(x))=\frac1{2}f''(s(x))s'(x). \label{gen.8q}
\end{equation}
Note that when $s(x)=x$  we have $L f(x))= \frac1{2}f''(x)$.

\bt\label{theo-diffsq} A function   
    $  f_s\in C^{2}(0,\ff)$ is excessive  for $\ov Z $ if and only if  
    \be   f_s(x)=f(s(x)) \text{ for some  $0\leq f\in C^{2}(0,\ff)$ with  $f''(x)\leq 0$.}
    \ee  
\et
  It is helpful to note that the conditions on $f$ imply that it is an increasing positive concave function.  

\medskip\Proof  
 We first show that a positive function $h$ is   excessive for $\ov Z $ if and only if  $L  h(x)\le 0$. Recall that a   function $h\geq 0$ is   excessive for $\ov Z $ if  
\bea\label{9.96}
 E^{x}\(h( \ov Z_{t} )\) &\leq& h(x), \forall x,t\in (0,\ff)\quad\text{and}\label{9.97a}\\ E^{x}\(h( \ov Z_{t} )\) &\to& h(x) \text{ as  $t\downarrow 0$, $\forall$  $x\in  (0,\ff)$.}  \label{9.97}
  \eea
  The   condition in (\ref{9.97}) holds for all  $h\in C_{b}(0,\ff)$ because the paths of $ \ov Z $ are continuous.   
 
 It follows from \cite[Chapter VII, Proposition 2.2]{RY}  that when   $h\in C_{b}^{2}(0,\ff)$, 
\begin{equation}
E^{x}\(h( \ov Z_{t} )\)-h(x)=E^{x}\(\int_{0}^{t}L  h( \ov Z_{s} )\,ds\).\label{ry.ex}
\end{equation}
Consequently, if $L h(z)\leq 0$ for all $z\in(0,\ff) $ then $ E^{x}\(h(  \ov Z_{t} \)\leq  h(x)$ for all $x,t\in(0,\ff) $ and, as we just pointed out, (\ref{9.97}) holds.  
  Conversely, (\ref{ry.ex}) implies that \be
\lim_{t\to 0}\frac{E^{x}\(h( \ov Z_{t} )\)-h(x)}{t}=L h(x),
\ee 
so that if $h$ is excessive for   $\ov Z $,   which implies    that (\ref{9.96}) holds,  we must have $L h(x)\leq 0$ for all $x$.

  Since $s'(x)>0$ it follows from (\ref{gen.8q}) that $ L\(f(s(x))\)\leq 0$, that is, $f(s(x))$ is excessive  for $\ov Z$, if and only if $f''\leq 0$.  Finally,   note that  any function  $ f_s\in C_{b}^{2}(0,\ff)$ can be written in the form $  f_s(x)=f(s(x))$ where $f_s(x)= f (s^{-1}(x))\in C_{b}^{2}(0,\ff)$.  This shows that the criteria in this theorem  describes  all excessive functions of   $\ov Z$. \qed

% \frac{In Theorem \ref{theo-1.5}, \textbf{c.} we require that $  f_s'(x_{0})=0$ for some $x_0\in (0,\ff)$.  Since $  f_s'(x)=f'(s(x)) s'(x)$ and $s'(x)>0$,   we see that when    $  f_s'(x_{0})=0$ for some $x_0\in (0,\ff)$  we must have  $f'(s(x_{0}))=0$. Consequently, the  condition that $f''(x)\leq 0$,   implies that   $f'(x)=0 $  for all $x\geq s(x_{0})$ \nc and that $f'_s(x)=0 $  for all $x\geq  x_{0} $. All such  functions   $f\in C^{2}(0,\ff)$ can be written as, \textbf{this notation is very confusing. We are using $  f_s(x)=f(s(x))$ but here it seems that $f_{  s(x_{0})}(x)$ is function of $x$. How about:}

  In Theorem \ref{theo-1.5}, \textbf{c.} we require that $  f_s'(x_{0})=0$ for some $x_0\in (0,\ff)$.  Since $  f_s'(x)=f'(s(x)) s'(x)$ and $s'(x)>0$,   we see that when    $  f_s'(x_{0})=0$ for some $x_0\in (0,\ff)$  we must have  $f'(s(x_{0}))=0$. Consequently,  the  condition that $f''(x)\leq 0$,   implies that   $f'(x)=0 $  for all $x\geq s(x_{0})$. Let $\wh  f $ denote the functions $f$ with these additional properties. All such  functions   $\wh f\in C^{2}(0,\ff)$ can be written as, 
 \be
  \label{9.79q} \wh f (y)=\left\{
\begin{array} {cc}
 \displaystyle    C-h\( {s(x_{0})-y} \)  & \qquad y\in (0, s(x_{0} )] \\\\
    \displaystyle   C & y>  s(x_{0}) 
\end{array}\right.,
\ee 
  for any constant $C>0$ and increasing  convex function $h\in   C^{2}(0,s(x_0)]$ with $h(0)=h'(0)=h''(0)=0$,  and $h(s ({x_0}))\le C$.   

Here are   specific examples, 
\be
  \label{9.79} \wh f^{(p)}(y)=\left\{
\begin{array} {cc}
 \displaystyle    1-\(\frac{s(x_{0})-y}{s(x_{0})}\)^{p} & \qquad y\in (0, s(x_{0}) ] \\\\
    \displaystyle  1 & y>s(x_{0})
\end{array}\right.,
\ee 
 for any $p>2$.   
 
 \begin{lemma} \label{lem-9.11q}
For any   $x_{0 }>0$ there exist excessive functions $  f_s$ and $ g_s$ for $\ov Z$ such that  $ f_s$ and $  g_s  \in C^{2}(0, \ff )$,     $  f_s'(x_0)= g_s'(x_0)=0$ and   $  g_s(x)\ne C  f_s(x)$ for $x\in \De_{x_0}(\de) $, for any $\de>0$ and constant  $C'$. 
\end{lemma}

\Proof   By (\ref{9.79q})   we can take any two   functions $h$, that satisfy the conditions in (\ref{9.79q}), that are not equal on $[s(x_0)-\de,s(x_0)]$ for any $\de>0$. To be more specific let    $f_s (x)= \wh f^{(p)}(s(x))$ and $g_s (x)=\wh f^{(q)}(s(x))$     with $p\ne q>2$. It is easy to see that $  g_s(x)\ne C'  f_s(x)$ for $x\in \De_{x_0}(\de) $, for any $\de>0$ and constant  $C'$.\qed

   The next theorem  gives a  different complete description 
   of all the excessive functions considered    Theorem \ref{theo-diffsq} and shows that they are potentials.

    \begin{theorem} \label{theo-9.4} Let  $0\leq f\in C^{2}(0,\ff)$  with    $f''(x)\leq 0$ for all $x\in (0,\ff)$,
    and   $f'(x)=0$ for all $x\ge s(x_0)$ for some $x_0\in (0,\ff)$. Let $f_s( x) =f(s(x))$ for all $x\in (0,\ff)$. Then 
   \bea \label{9.94mm}   f_s( x) &=&\int^{s(x_0)  }_{0}2(s(x)\wedge y )\(-  \frac{1}{2} f''(y)\)\,dy\\
          &=& \int^{  x_0  }_{0}2(s(x)\wedge s(y) )\(-L    f_s(y)\)\,dy\nn.
           \eea
  \end{theorem} 
  
  \Proof{\rm \label{rem-9.1} 
It follows from the Riesz decomposition theorem, \cite[VI, (2.11),  (2.19)]{BG}, that $f$ has a unique representation of the form,
\begin{equation}
f(x)=C_{0}x+\int_{0}^{\ff}(x\wedge y) \,d\mu(y),\label{RD.1}
\end{equation}
for some constant $C_{0}\geq 0 $ and finite positive measure $\mu$.  Since for any finite measure $\nu$,
\begin{equation}
\lim_{x\to\ff}\frac{1}{x}\int_{0}^{\ff}(x\wedge y )\,d\nu(y)=0,\label{RD.0}
\end{equation}
(see, e.g.,     \cite[(5.36)]{MRejp}), and  $f'(x)=0$ for all $x\geq s(x_{0})$,    $C_{0}=0$.

Clearly,   $-f''(y)\,dy$ is a finite, positive   measure with compact support. We now show that $\mu(y)=-f''(y)\,dy$.
 To see this, note that
 \begin{eqnarray}
-\frac{d}{dx}\int^{\ff}_{0}(x\wedge y) \, \,f''(y)\,dy
 \label{d5.104a}& =&-\frac{d}{dx}\int^{x}_{0}y \,f''(y)\,dy-\frac{d}{dx}\(x \int^{\ff}_{x} \,f''(y)\,dy\) \nonumber\\
&=& - \int^{\ff }_{x}\,f''(y)\,dy=f'(x).
 \end{eqnarray}
Consequently, 
\begin{equation}
f(x)=C_{1}x+\int^{\ff}_{0}(x\wedge y) (-  f''(y))\,dy,\label{RD.2}
\end{equation} 
for some constant $C_{1}$.
Using (\ref{RD.0}) again and the fact that  $f'(x)=0$ for all $x\geq x_{0}$,  we see that $C_{1}=0$.  Therefore,
\begin{equation}
f(x)= \int^{s(x_{0})}_{0}(x\wedge y )\(-f''(y)\)\,dy.
\end{equation}
It follows that
\begin{eqnarray}
f_{s}(x)=f(s(x))&=& \int^{s(x_{0}) }_{0}(s(x)\wedge y )\(-f''(y)\)\,dy\\
&=&   \int^{ x_{0} }_{0}(s(x)\wedge s(z) )\(-f''(s(z))\)s'(z)\,dz.\nonumber
\end{eqnarray}
The first integral is first integral in (\ref{9.94mm}). The second integral and   (\ref{gen.8q})  gives the    second integral in (\ref{9.94mm}).\qed

 \begin{remark}{\rm  \label{rem-9.1}
  
   It follows from (\ref{potdef}) that for any $h\geq 0$ the potential function
\begin{equation}
G(x)= \int^{\ff}_{0}2(s(x)\wedge s(y) )h(y)\,dy\label{}
\end{equation}
is an excessive function for a diffusion with potential density given in (\ref{d2.1}). In Theorem \ref{theo-1.5},  {c.} we also require that   $G\in C^{2}(0,\ff)$ and $G'(x_{0})=0.$ We   now show that this is the case when  $h\in C(0,\ff)$ and 
$h(y)=0$ for $y\geq x_{0}$. In this case, 
\begin{equation}
G(x)= \int^{x_{0}}_{0}(s(x)\wedge s(y) )h(y)\,dy.\label{}
\end{equation}
We have,  
 \begin{eqnarray}
G'(x)&=& \frac{d}{dx}\int^{x_{0}}_{0}2(s(x)\wedge s(y) )h(y)\,dy
 \label{dover.m}\\
 & =& \frac{d}{dx}\int^{x\wedge x_{0}}_{0}2s(y)\, \,h(y)\,dy+\frac{d}{dx}\( 2s(x) \int^{x_{0}}_{x\wedge x_{0}} \,h(y)\,dy\) \nonumber\\
&=&  2 s(x)\, \,h(x)1_{\{x\leq  x_{0}\}}+2s'(x) \int^{x_{0}}_{x\wedge x_{0}} \, h(y)\,dy-2s(x)\, \,h(x)1_{\{x\leq  x_{0}\}}\nn\\
&=& 2s'(x) \int^{x_{0}}_{x\wedge x_{0}} \,h(y)\,dy.
 \end{eqnarray}
This shows that that $G'(x_{0})=0$. 
In addition, 
 \begin{eqnarray}
G''(x)&=&\frac{d}{dx}2s'(x) \int^{x_{0}}_{ x\wedge x_{0}} \,h(y)\,dy
 \label{dover.m}\\
 & =& 2s''(x) \int^{x_{0}}_{x\wedge x_{0}} \,h(y)\,dy -2s'(x)\, \,h(x)1_{\{x\leq  x_{0}\}} \nonumber ,
 \end{eqnarray}
which shows that $G\in C^{2}(0,\ff)$. 

  Note that $G'(x_{0})=G''(x_{0})=0$ and 
\begin{equation}
G(x_{0})= \int^{x_{0}}_{0}2(s(x_{0})\wedge s(y) )h(y)\,dy=\int_0^{x_0} 2s(y) h(y)\,dy.\label{mlr.1}
\end{equation}

}\end{remark}

\begin{example}{\rm  \label{ex-9.4} Using  (\ref{1.44mm}) we see that when $\psi(\la)=C\la^2$,
  \begin{equation} 
   (\si^{(0)})^2(x)=2\phi(x)= \frac{2}{ \pi}\int_{0}^{\ff}\frac{1-\cos\la x}{C\la^2}\,d\la =\frac{2}{C} |x|
\end{equation}
Using (\ref{1.44mm}) we see that for $x,y>0$,
\begin{equation} 
  u^{(0)}(x,y)=\frac{1}{C} (2x\wedge 2y).
\end{equation}
It follows from Theorem \ref{theo-diffsq} that the excessive functions for $u^{(0)}(x,y)$ are all functions $0\leq f\in C^{2}(0,\ff)$ with $f''(x)\le 0$.  
  }\end{example}

 \end{document}